\documentclass[12pt]{article}
\usepackage[round,colon,authoryear]{natbib}
\usepackage{bibentry}
\RequirePackage[OT1]{fontenc}
\RequirePackage{amsthm,amsmath,amsfonts,amssymb,{mathtools}}
\RequirePackage[colorlinks]{hyperref}
\usepackage{bbm}
\hypersetup{
 colorlinks=false,
 citecolor=Violet,
 linkcolor=Red,
 urlcolor=Blue}
\RequirePackage{epsfig, graphicx, color}
\RequirePackage{latexsym}
\RequirePackage{psfrag}

\allowdisplaybreaks

\usepackage{fullpage}
\usepackage{epsf,epsfig,subfigure,caption}
\usepackage{graphics,graphicx}
\usepackage{epstopdf}
\usepackage{algorithm,algorithmicx,algpseudocode}

\newtheorem{theorem}{Theorem}
\newtheorem{lemma}{Lemma}

\newtheorem{assump}{Assumption}

\newcommand{\bfm}[1]{\ensuremath{\mathbf{#1}}}

     \def\EE{\mathbb{E}}

     \def\OO{\mathbb{O}}
     \def\PP{\mathbb{P}}
     
     \def\RR{\mathbb{R}}
\def\bs{\bfm s}

     \def\ZZ{\mathbb{Z}}

\def\calE{{\cal  E}} 
 
\def\calG{{\cal  G}} 
 
\def\calI{{\cal  I}}

\def\calN{{\cal  N}} 
\def\calO{{\cal  O}} 
\def\calP{{\cal  P}}

\def\calS{{\cal  S}}

\DeclareMathOperator{\tr}{tr}

\def\newpage{\vfill\eject}

\newdimen\biblioindent    \biblioindent=30pt

 at 8truept

\def\tran{\mathsf{T}}
\def\what{\widehat}
\def\wtilde{\widetilde}

\newcommand{\beq}{\begin{equation}}
  \newcommand{\eeq}{\end{equation}}
\newcommand{\beqn}{\begin{eqnarray}}
  \newcommand{\eeqn}{\end{eqnarray}}
\newcommand{\beqnn}{\begin{eqnarray*}}
  \newcommand{\eeqnn}{\end{eqnarray*}}

\numberwithin{equation}{section}

\usepackage {tikz}
\usetikzlibrary {positioning}

\usepackage{mathtools}
\DeclarePairedDelimiter\ceil{\lceil}{\rceil}

\newcounter{CondCounter}

\def\mumax{\mu_{\scriptscriptstyle \sf max}}

\def\supinit{\scriptscriptstyle \sf init}
\def\supunbiased{\scriptscriptstyle \sf unbs}
\def\submaxx{\scriptscriptstyle \sf 2,max}
\def\submax{\scriptscriptstyle \sf max}
\def\subf{\scriptscriptstyle \sf F}
\def\subo{\scriptscriptstyle \sf O}
\def\kmax{k_{\scriptscriptstyle \sf max}}
\def\supklt{\scriptscriptstyle \sf KLT}
\def\supkmo{\scriptscriptstyle \sf KMO}
\def\supmwc{\scriptscriptstyle \sf MWC}

\begin{document}
\title{Statistical Inferences of Linear Forms for Noisy Matrix Completion$^{\ast}$}
\author{Dong Xia$^{\dagger}$ and Ming Yuan$^{\sharp}$\\
$^{\dagger}${ Hong Kong University of Science and Technology}\\
$^{\sharp}${ Columbia University}}
\date{}

\maketitle

\footnotetext[1]{
Dong Xia's research is partially supported by Hong Kong RGC Grant ECS 26302019. Ming Yuan's research is supported in part by NSF Grant DMS-1803450.}

\begin{abstract}
We introduce a flexible framework for making inferences about general linear forms of a large matrix based on noisy observations of a subset of its entries. In particular, under mild regularity conditions, we develop a universal procedure to construct asymptotically normal estimators of its linear forms through double-sample debiasing and low-rank projection whenever an entry-wise consistent estimator of the matrix is available. These estimators allow us to subsequently construct confidence intervals for and test hypotheses about the linear forms. Our proposal was motivated by a careful perturbation analysis of the empirical singular spaces under the noisy matrix completion model which might be of independent interest. The practical merits of our proposed inference procedure are demonstrated on both simulated and real-world data examples.
\end{abstract}

\section{Introduction}\label{sec:intro}
Noisy matrix completion (NMC) refers to the reconstruction of a low rank matrix $M\in\RR^{d_1\times d_2}$ after observing a small subset of $M$'s entries with random noise. Problems of this nature arise naturally in various applications. For the sake of generality, we shall cast it in the framework of trace regression where each observation is a random pair $(X,Y)$ with $X\in\RR^{d_1\times d_2}$ and $Y\in\RR$. The random matrix $X$ is sampled uniformly from the orthonormal basis $\mathfrak{E}=\{e_{j_1}e_{j_2}^{\tran}: j_1\in[d_1], j_2\in[d_2]\}$ where $[d]=\{1,\cdots,d\}$ and $\{e_{j_1}\}_{j_1\in[d_1]}$ and $\{e_{j_2}\}_{j_2\in[d_2]}$ are the canonical basis vectors in $\RR^{d_1}$ and $\RR^{d_2}$, respectively. It is worth pointing out that, while we shall focus on the canonical basis in this work, our framework can be easily extended to general product basis where $\{e_{j_1}\}_{j_1\in[d_1]}$ and $\{e_{j_2}\}_{j_2\in[d_2]}$ are arbitrary orthonormal basis in $\RR^{d_1}$ and $\RR^{d_2}$, respectively.  Without loss of generality, we shall assume that $d_1\geq d_2$ and denote $\alpha_d=d_1/d_2$ the aspect ratio of $M$. The response variable $Y$ is related to $X$ via
\begin{equation}\label{eq:modelXY}
Y=\langle  M, X\rangle +\xi
\end{equation}
where $\langle M, X\rangle=\tr(M^{\tran}X)$, and the independent measurement error $\xi$ is assumed to be a centered sub-Gaussian random variable. Our goal is to infer $M$ from $n$ i.i.d. copies $\{(X_i,Y_i)\}_{i=1}^n$ obeying (\ref{eq:modelXY}) when, in particular, $M$ is of (approximately) low rank and $n$ is much smaller than $d_1d_2$. 

In the absence of measurement error (e.g., $\xi=0$), \cite{candes2009exact} first discovered that exact matrix completion can be solved efficiently by relaxing the non-convex and non-smooth rank constraint of a matrix to its nuclear norm. Following the pioneering work, nuclear-norm penalized least squares estimators \citep{negahban2011estimation,rohde2011estimation,cai2010singular, cai2016matrix, candes2009power, candes2010matrix,gross2011recovering} and numerous other variants \citep{koltchinskii2011nuclear,klopp2014noisy,liu2011universal,recht2010guaranteed,sun2012calibrated,cai2015rop, gao2016optimal} have been studied. It is now understood, from these developments, that even when the observations are contaminated with noise, statistically optimal convergence rates are attainable by efficiently computable convex methods. For instance, \cite{koltchinskii2011nuclear} proved that a modified matrix LASSO estimator, denoted by $\what M^{\supklt}$, achieves the convergence rate:
\begin{equation}\label{eq:KLT_rate}
\|\what M^{\supklt}-M\|_{\subf}^2=O_P\Big((\sigma_{\xi}+ \|M\|_{\submax})^2\cdot\frac{rd_1^2d_2\log d_1}{n}\Big)
\end{equation}
as long as $n\gg d_1\log d_1$ where $r$ is the rank of $M$ and $\sigma_{\xi}^2$ is the variance of $\xi$. Here, $\|\cdot\|_{\subf}$ denotes the matrix Frobenius norm and $\|\cdot\|_{\submax}$ denotes the max-norm defined as $\|A\|_{\submax}=\max_{j_1\in[d_1], j_2\in[d_2]}|A(j_1,j_2)|$. It is worth noting that (\ref{eq:KLT_rate}) was established without additional assumptions on $M$. As a result, the rate given on the righthand side of (\ref{eq:KLT_rate}) depends on $\|M\|_{\submax}$ and does not vanish even when $\sigma_{\xi}=0$.

In addition to convex methods, non-convex approaches such as those based on matrix-factorization have also been developed. For instance, \cite{keshavan2010matrix_b} proposed a non-convex estimator based on the thin SVD, denoted by $\what M^{\supkmo}$, and show that
\begin{equation}\label{eq:KMO_rate}
\|\what M^{\supkmo}-M\|_{\subf}^2=O_P\Big(\sigma_{\xi}^2\cdot\frac{rd_1^2d_2\log d_1}{n}\Big)
\end{equation}
assuming that $n\gg rd_1(r+\log d_1)$ and $M$ satisfies the so-called incoherent condition. See also, e.g., \cite{zhao2015nonconvex,chen2015fast, cai2016structured} and references therein.
The rate (\ref{eq:KMO_rate}) is optimal up to the logarithmic factors, see, e.g., \cite{koltchinskii2011nuclear} and \cite{ma2015volume}, for a comparable minimax lower bound. More recently, an alternative scheme of matrix factorization attracted much attention. See, e.g., \cite{wang2016unified,ge2016matrix,zheng2016convergence,chen2019inference,chen2019noisy,ma2017implicit,chen2019nonconvex}. In particular, \cite{ma2017implicit} showed this approach yields an estimator, denoted by $\what M^{\supmwc}$, that is statistically optimal not only in matrix Frobenius norm but also in entry-wise max-norm, i.e., 
\begin{equation}\label{eq:MWC_rate}
\|\what M^{\supmwc}-M\|_{\submax}^2=O_P\Big(\sigma_{\xi}^2\cdot\frac{rd_1\log d_1}{n}\Big)
\end{equation}
provided that $n\gg r^3d_1\log^3 d_1$. 

While there is a rich literature on statistical estimation for NMC, results about its statistical inferences are relatively scarce. In \cite{carpentier2015uncertainty}, a debiasing procedure, based on sample splitting, was proposed for the nuclear norm penalized least squares estimator which enables constructing confidence region for $M$ with respect to matrix Frobenius norm when $n\gg rd_1\log d_1$. Their technique, however, cannot be directly used to make inferences about individual entries or linear forms as confidence regions for $M$ with respect to matrix Frobenius norm can be too wide for such purposes. To this end, \cite{carpentier2018adaptive} proposed another procedure to construct entrywise confidence intervals. However their procedure requires that the design, namely the underlying distribution of $X$ satisfy the so-called {\it restricted isometry property} which is violated when $X$ is sampled uniformly from $\mathfrak{E}$. Another proposal introduced by \cite{cai2016geometric} can be used to construct confidence intervals for $M$'s entries. However, it requires that the sample size $n\gg d_1d_2$ which is significantly larger than the optimal sample size requirement for estimation. In addition, during the preparation of the current work, \cite{chen2019inference} announced a different approach to constructing confidence intervals for the entries of $M$.

The present article aims to further expand this line of research by introducing a flexible framework for constructing confidence intervals and testing hypotheses about general linear forms of $M$, with its entries as special cases, under optimal sample size requirement. In a nutshell, we develop a procedure that, given any entry-wise consistent estimator $\what M^{\supinit}$ in that $\|\what M^{\supinit}-M\|_{\submax}=o_P(\sigma_{\xi})$, can yield valid statistical inferences for $m_T:=\tr(M^{\tran}T)$ under mild regularity conditions. More specifically, we show that, through double-sample debiasing and spectral projection, we can obtain from the initial estimator a new one, denoted by $\what M$, so that
\begin{equation}\label{eq:mproj_sigma}
\frac{\tr(\what M^{\tran}T)-\tr(M^{\tran}T)}{\sigma_{\xi}(\|U^{\tran}T\|_{\subf}^2+\|TV\|_{\subf}^2)^{1/2}\sqrt{d_1d_2/n}}\overset{{\rm d}}{\longrightarrow}\calN(0,1),
\end{equation}
provided that
\begin{align*}
\|U^{\tran}T&\|_{\subf}+\|TV\|_{\subf}
\gg \|T\|_{\ell_1}\sqrt{\frac{r}{d_1}}\cdot \max\bigg\{\sqrt{\frac{r\log d_1}{d_2}},\   \frac{\sigma_{\xi}}{\lambda_r}\sqrt{\frac{\alpha_d rd_1^2d_2\log^2 d_1}{n}}\bigg\}
\end{align*} 
where $U, V$ are $M$'s left and right singular vectors and $\lambda_r$ is its $r$-th singular value, and $\|\cdot\|_{\ell_1}$ stands for the vectorized $\ell_1$ norm. We not only show that (\ref{eq:mproj_sigma}) holds under optimal sample size (independent of $T$) but also derive its non-asymptotic convergence rate explicitly. Note that condition for $\|U^{\tran}T\|_{\subf}+\|TV\|_{\subf}$ in a certain sense is necessary to avoid non-regular asymptotic behavior when $\|U^{\tran}T\|_{\subf}+\|TV\|_{\subf}=0$. Moreover, we show that under similar conditions, \eqref{eq:mproj_sigma} continues to hold when we replace $\sigma_{\xi}$, $\|U^{\tran}T\|_{\subf}$ and $\|TV\|_{\subf}$ by suitable estimates, denoted by $\what\sigma_{\xi}$, $\|\what U^{\tran}T\|_{\subf}$ and $\|T\what V\|_{\subf}$ respectively:
\begin{equation}\label{eq:mproj_hatsigma}
\frac{\tr(\what M^{\tran}T)-\tr(M^{\tran}T)}{\what\sigma_{\xi}(\|\what U^{\tran}T\|_{\subf}^2+\|T\what V\|_{\subf}^2)^{1/2}\sqrt{d_1d_2/n}}\overset{{\rm d}}{\longrightarrow}\calN(0,1).
\end{equation}
The statistic on the lefthand side is now readily applicable for making inferences about the linear form $\tr(M^{\tran}T)$.

Our proposal greatly generalizes the scope of earlier works on inferences for entries of $M$ in several crucial aspects. Firstly, unlike earlier approaches that focus on a specific estimator of $M$, our procedure can be applied to any entry-wise consistent estimator. This not only brings potential practical benefits but also helps us better understand the fundamental differences between estimation and testing in the context of NMC. For instance, our results suggest that, perhaps surprisingly, when it comes to make valid inferences with optimal sample sizes, the rate of convergence of the initial estimate is irrelevant as long as it is consistent; therefore a suboptimal estimator may be used for making optimal inferences.

Secondly, our approach can be applied in general when $T$ is sparse, and depending on its alignment with the singular spaces of $M$, even to cases where it is dense and $\|T\|_{\ell_1}^2/\|T\|_{\subf}^2$ is of the order $O(d_2)$. Entry-wise inferences correspond to the special case when $T$ takes the form $e_{i}e_{j}^{\tran}$. Extensions to more general linear forms could prove useful in many applications. For example, in recommender systems, it may be of interest to decide between items $j_1$ and $j_2$ which should we recommend to user $i$. This can obviously be formulated as a testing problem:
\begin{equation}\label{eq:HT}
H_0: M(i,j_1)=M(i,j_2)\quad {\rm v.s.}\quad H_1: M(i,j_1)>M(i,j_2),
\end{equation}
which can be easily solved within our framework by taking $T=e_{i}e_{j_1}^{\tran}-e_{i}e_{j_2}^{\tran}$. More generally, if the target is a group of users $\calG\subset [d_1]$, we might take a linear form $T=\sum_{i\in\calG}e_{i}(e_{j_1}-e_{j_2})^{\tran}$. At a technical level, inferences about general linear forms as opposed to entries of $M$ present nontrivial challenges because of the complex dependence structure among the estimated entries. As our theoretical analysis shows, the variance of the plug-in estimator for the linear form depends on the alignment of the linear form with respect to the singular space of $M$ rather than the sparsity of the linear form.

An essential part of our technical development is the characterization of the distribution of the empirical singular vectors for NMC where we take advantage of the recently developed spectral representation for empirical singular vectors. Similar tools have been used earlier to derive confidence regions for singular subspaces with respect to $\ell_2$-norm for low-rank matrix regression (LMR) when the linear measurement matrix $X$s are Gaussian \citep{xia2018confidence}, and the planted low rank matrix (PLM) model where every entry of $M$ is observed with i.i.d. Gaussian noise \citep{xia2019normal}. In both cases, Gaussian assumption plays a critical role and furthermore, it was observed that first order approximation may lead to suboptimal performances. In absence of the Gaussian assumption, the treatment of NMC is technically more challenging and requires us to derive  sharp bounds for the $(2,\max)$-norm for the higher order perturbation terms. Interestingly, it turns out that, unlike LMR or PLM, a first order approximation actually suffices for NMC.

Even though our framework applies to any max-norm consistent matrix estimator, for concreteness, we introduce a novel rotation calibrated gradient descent algorithm on Grassmannians that yields such an initial estimator. The rotation calibration promotes fast convergence on Grassmannians so that constant stepsize can be selected to guarantee geometric  convergence. We note that existing results on max-norm convergence rates are established for sampling without replacement \citep{ma2017implicit}. It is plausible that \eqref{eq:MWC_rate} may continue to hold under our assumption of independent sampling given the close connection between the two sampling schemes, but an actual proof is likely much more involved and therefore we opted for the proposed alternative for illustration as it is more amenable for analysis.

The rest of our paper is organized as follows. In next section, we present a general framework for estimating $m_T=\tr(M^{\tran}T)$ given an initial estimator through double-sample-debiasing and spectral projection. In Section~\ref{sec:theory}, we establish the asymptotic normality of the estimate obtained. In Section~\ref{sec:conf_int}, we propose data-driven estimates for the noise variance and the true singular vectors, based on which confidence intervals of $m_T$ are constructed. In Section~\ref{sec:init_ncx}, we introduce a rotation calibrated gradient descent algorithm on Grassmannians, which, under mild conditions, provides the initial estimator $\what M^{\supinit}$ so that $\|\what M^{\supinit}-M\|_{\submax}=o_P(\sigma_{\xi})$.  Numerical experiments on both synthetic and real world datasets presented in Section~\ref{sec:numerical} further demonstrate the merits of the proposed methodology. All proofs are presented in the online supplement.

\section{Estimating Linear Forms}\label{sec:method}
We are interested in making inferences about $m_T=\tr(M^{\tran}T)$ for a given $T$ based on observations $\mathfrak{D}=\{(X_i,Y_i): 1\leq i\leq n\}$ satisfying model (\ref{eq:modelXY}), assuming that $M$ has low rank. To this end, we first need to construct an appropriate estimate of $m_T$ which we shall do in this section. 

Without loss of generality, we assume $n$ is an even number with $n=2n_0$, and split $\mathfrak{D}$ into two sub-samples:
$$
\mathfrak{D}_1=\big\{(X_i, Y_i)\big\}_{i=1}^{n_0}\quad {\rm and}\quad \mathfrak{D}_2=\big\{(X_i, Y_i)\big\}_{i=n_0+1}^{n}.
$$
In what follows, we shall denote $M$'s thin singular value decomposition (SVD) by $M=U\Lambda V^{\tran}$ where $U\in\OO^{d_1\times r}, V\in\OO^{d_2\times r}$ and $\Lambda={\rm diag}(\lambda_1,\cdots,\lambda_r)$ represent $M$'s singular vectors and singular values, respectively. The Stiefel manifold $\OO^{d\times r}$ is defined as 
$
\OO^{d\times r}:=\big\{A\in \RR^{d\times r}: A^{\tran}A=I \big\}. 
$
We arrange $M$'s positive singular values non-increasingly, i.e., $\lambda_1\geq\cdots\geq \lambda_r>0$. 

Assuming the availability of an initial estimator, our procedure consists of four steps as follows: 

\begin{itemize}
\item {\it Step 1 (Initialization)}: By utilizing the first and second data sub-sample $\mathfrak{D}_1, \mathfrak{D}_2$ separately, we apply the initial estimating procedure on noisy matrix completion to yield initial (biased in general) estimates $\widehat M_1^{\supinit}$ and $\what M_2^{\supinit}$, respectively. 

\item{\it Step 2 (Debiasing)}:  Using the second data sub-sample $\mathfrak{D}_2$, we debias $\what M_1^{\supinit}$:
$$
\what M_1^{\supunbiased}=\what M_1^{\supinit}+\frac{d_1d_2}{n_0}\sum_{i=n_0+1}^{n}\big(Y_i-\langle \what M_1^{\supinit}, X_i\rangle\big)X_i.
$$
Similarly, we use the first data sub-sample $\mathfrak{D_1}$ to debias $\what M_2^{\supinit}$ and obtain
$$
\what M_2^{\supunbiased}=\what M_2^{\supinit}+\frac{d_1d_2}{n_0}\sum_{i=1}^{n_0}\big(Y_i-\langle \what M_2^{\supinit}, X_i\rangle\big)X_i.
$$
\item{\it Step 3 (Projection)}:  Compute the top-$r$ left and right singular vectors of $\what{M}_1^{\supunbiased}$, denoted by $\what U_1$ and $\what V_1$. Similarly, compute the top-$r$ left and right singular vectors of $\what{M}_2^{\supunbiased}$, denoted by $\what U_2$ and $\what V_2$. 
Then, we calculate the (averaged) projection estimate 
$$
\what{M}=\frac{1}{2}\what U_1\what U_1^{\tran}\what M_1^{\supunbiased} \what V_1\what V_1^{\tran}+\frac{1}{2}\what U_2\what U_2^{\tran}\what M_2^{\supunbiased} \what V_2\what V_2^{\tran}.
$$

\item {\it Step 4 (Plug-in)}: Finally, we estimate $m_T$ by $\what m_T=\tr(\what{M}^{\tran}T)$.
\end{itemize}

We now discuss each of the steps in further details.

\paragraph{Initialization.}
Apparently, our final estimate depends on the initial estimates $\what M_1^{\supinit}, \what M_2^{\supinit}$. However, as we shall show in the next section, such dependence is fairly weak and the resulting estimate $\what m_T$ is asymptotically equivalent as long as the estimation error of $\what M_1^{\supinit}$ and $\what M_2^{\supinit}$, in terms of max-norm, is of a smaller order than $\sigma_{\xi}$. More specifically, we shall assume that
\begin{assump}\label{assump:init_entry}
There exists a sequence $\gamma_{n,d_1,d_2}\to 0$ as $n, d_1,d_2\to\infty$ so that with probability at least $1-d_1^{-2}$,
\begin{eqnarray}\label{eq:M1M2_rate}
\|\what M_1^{\supinit}-M\|_{\scriptscriptstyle \sf max}+\|\what M_2^{\supinit}-M\|_{\submax}\leq C\gamma_{n,d_1,d_2}\cdot \sigma_{\xi}
\end{eqnarray}
for an absolute constant $C>0$.
\end{assump}

In particular, bounds similar to \eqref{eq:M1M2_rate} have recently been established by \cite{ma2017implicit, chen2019inference}. See eq. (\ref{eq:MWC_rate}). Assumption \ref{assump:init_entry} was motivated by their results. However, as noted earlier, (\ref{eq:MWC_rate}) was obtained under sampling without replacement and for positively semi-definite matrices. While it is plausible that it also holds under independent sampling as considered here, an actual proof is lacking at this point. For concreteness, we shall present a simple algorithm in Section~\ref{sec:init_ncx} capable of producing an initial estimate that satisfies Assumption~\ref{assump:init_entry}.

\paragraph{Debiasing.} The initial estimate is only assumed to be consistent. It may not necessarily be unbiased or optimal. To ensure good quality of our final estimate $\what m_T$, it is important that we first debias it which allows for sharp spectral perturbation analysis. Debiasing is an essential technique in statistical inferences of high-dimensional sparse linear regression \citep[see, e.g.,][]{zhang2014confidence,javanmard2014confidence, van2014asymptotically, cai2017confidence} and low-rank matrix regression \citep[see, e.g.,][]{cai2016geometric, carpentier2015iterative, carpentier2018adaptive,xia2018confidence}. Oftentimes, debiasing is done in absence of the knowledge of $\EE {\rm vec}(X){\rm vec}(X)^{\top}$ and a crucial step is to construct an appropriate decorrelating matrix. In our setting, it is clear that $\EE {\rm vec}(X){\rm vec}(X)^{\top}= (d_1d_2)^{-1}I_{d_1d_2}$. This allows for a much simplified treatment via sample splitting, in the same spirit as earlier works including \cite{carpentier2015uncertainty, xia2018confidence}, among others. The particular double-sample-splitting technique we employ was first proposed by \cite{chernozhukov2018double} and avoids the loss of statistical efficiency associated with the sample splitting. It is worth noting that if the entries are not sampled uniformly, the debiasing procedure needs to be calibrated accordingly.

In addition to reducing possible bias of the initial estimate, the sample splitting also enables us to extend the recently developed spectral representation for empirical singular vectors under Gaussian assumptions to general sub-Gaussian distributions.




\paragraph{Spectral Projection.}  Since $M$ have low rank, it is natural to apply spectral truncation to a matrix estimate to yield an improved estimate. To this end, we project $\what M_1^{\supunbiased}$ and $\what{M}_2^{\supunbiased}$ to their respective leading singular subspaces. Note that, while $\what M_1^{\supunbiased}, \what{M}_2^{\supunbiased}$ are unbiased, their empirical singular vectors $\what U_1, \what U_2, \what V_1$ and $\what V_2$ are typically not. The spectral projection serves the purpose of reducing entry-wise variances at the cost of negligible biases.

It is worth noting that the estimate $\what M$ may not be of rank $r$. If an exact rank-$r$ estimator is desired, it suffices to obtain the best rank-$r$ approximation of  $\what M$ via singular value decomposition and all our development remains valid under such a modification.  In general, getting the initial estimates is the most computational expensive step as the other steps involving fairly standard operation without incurring any challenging optimization. This is noteworthy because it suggests that as long as we can compute a good estimate, it does not cost much more computationally to make inferences.

\section{Asymptotic Normality of $\what m_T$}\label{sec:theory_normal}\label{sec:theory}

We now show the estimate $\what m_T$ we derived in the previous section is indeed suitable for inferences about $m_T$ by establishing its asymptotic normality.

\subsection{General results}
For brevity, let $e_{j}$ denote the $j$-th canonical basis in $\RR^{d}$ where $d$ might be $d_1$ or $d_2$ or $d_1+d_2$ at different appearances. With slight abuse of notation, denote by $\|\cdot\|$ the matrix operator norm or vector $\ell_2$-norm depending on the dimension of its argument. Denote the condition number of $M$ by
\begin{equation}\label{eq:cond_num}
\kappa(M)=\lambda_1(M)/\lambda_r(M)=\lambda_1/\lambda_r.
\end{equation}
As is conventional in the literature, we shall assume implicitly that rank $r$ is known with $r\ll d_2$ and $M$ is well-conditioned so that $\kappa(M)\leq \kappa_0$. 
In practice, $r$ is usually not known in advance and needs to be estimated from the data. Our experience with numerical experiments such as those reported in Section \ref{sec:numerical} suggests that our procedure is generally robust to reasonable estimate of $r$. Although a more rigorous justification of such a phenomenon has thus far eluded us, these promising empirical observations nonetheless indicate a more careful future investigation is warranted.

In addition, we shall assume that $U$ and $V$ are incoherent, a standard condition for matrix completion.
\begin{assump}\label{assump:incoh} Let $\|U\|_{\submaxx}=\max_{j\in[d_1]}\|e_{j}^{\tran}U\|$ and there exists $\mumax>0$ so that 
$$
\max\Big\{\sqrt{\frac{d_1}{r}}\|U\|_{\submaxx},\sqrt{\frac{d_2}{r}}\|V\|_{\submaxx}\Big\}\leq \mumax.
$$
\end{assump}

We also assume that the noise $\xi$ is independent with $X$ and sub-Gaussian such that
\begin{assump}\label{assump:noise}
The noise $\xi$ is independent with $X$ and 
\begin{equation}\label{eq:noiseassump}
\EE \xi=0,\quad \EE\xi^2=\sigma_{\xi}^2,\quad {\rm and}\quad  \EE e^{s\xi}\leq e^{s^2\sigma_{\xi}^2},\  \forall s\in\RR
\end{equation}
Let $\alpha_d=d_1/d_2$. There exists a large enough absolute constant $C_1>0$ so that 
\begin{equation}\label{eq:SNR}
\lambda_r\geq C_1\mumax \kappa_0^2\sigma_{\xi}\sqrt{\frac{\alpha_d rd_1^2d_2\log^2d_1}{n}}.
\end{equation}
\end{assump}

The SNR condition (\ref{eq:SNR}) is optimal up to the logarithmic factors if $\alpha_d, \mumax, \kappa_0=O(1)$. Indeed, the consistent estimation of singular subspaces  requires $\lambda_r\gg \sigma_{\xi}\sqrt{ rd_1^2d_2/n}$. This condition is common for non-convex methods of NMC. However, when $\alpha_d\gg 1$, i.e., $M$ is highly rectangular, condition (\ref{eq:SNR}) is significantly stronger than the optimal SNR condition even if $\mumax, \kappa_0=O(1)$. It is unclear to us whether this sub-optimality is due to technical issues or reflection of more fundamental differences between statistical estimation and inference. 

To avoid the nonregular asymptotics, we focus on the case when $T$ does not lie entirely in the null space of $M$. More specifically, we assume that
\begin{assump}\label{assump:T}
There exists a constant $\alpha_T>0$ such that
$$
\|U^{\tran}T\|_{\subf}\geq \alpha_T\|T\|_{\subf}\cdot \sqrt{\frac{r}{d_1}}\quad {\rm or}\quad \|TV\|_{\rm F}\geq \alpha_T\|T\|_{\subf}\cdot \sqrt{\frac{r}{d_2}}.
$$
\end{assump}

The alignment parameter $\alpha_T$ in Assumption~\ref{assump:T} is allowed to vanish as $d_1,d_2,n\to\infty$. Indeed, as we show below, the asymptotic normality of $\what m_T-m_T$ only requires that
\begin{align}
\alpha_T\ge \frac{\|T\|_{\ell_1}}{\|T\|_{\subf}}\cdot\max\left\{\mumax^2\sqrt{\frac{r\log d_1}{d_2}},\   \frac{\kappa_0\mumax^2\sigma_{\xi}}{\lambda_r}\sqrt{\frac{\alpha_d rd_1^2d_2\log^2d_1}{n}}\right\}.\label{eq:alpha_T}
\end{align}
We are now in position to establish the asymptotic normality of $\what m_T$.

\begin{theorem}\label{thm:normal_mt}
Under Assumptions~\ref{assump:init_entry}-\ref{assump:T}, there exist absolute constants $C_1, C_2, C_3, C_4, C_5,C_6>0$ so that if $n\geq C_1\mumax^2 rd_1\log d_1$, then
\begin{align*}
\sup_{x\in\RR}\Big|\PP&\Big(\frac{\what m_T-m_T}{\sigma_{\xi}(\|TV\|_{\subf}^2+\|U^{\tran}T\|_{\subf}^2)^{1/2}\cdot\sqrt{d_1d_2/n}}\leq x\Big)-\Phi(x)\Big|\\
&\hspace{0cm}\leq C_2\frac{\mumax^2\|T\|_{\ell_1}}{\alpha_T\|T\|_{\subf}}\sqrt{\frac{\log d_1}{d_2}}+C_3\kappa_0\frac{\mumax^2\|T\|_{\ell_1}}{\alpha_T\|T\|_{\subf}}\cdot \frac{\sigma_{\xi}}{\lambda_r}\sqrt{\frac{\alpha_d rd_1^2d_2\log^2d_1}{n}}\\
&\hspace{0cm}+C_4\frac{\mumax^4\|T\|_{\ell_1}^2}{\alpha_T^2\|T\|_{\subf}^2}\cdot\frac{r\sqrt{\log d_1}}{d_2}+\frac{6\log d_1}{d_1^2}+C_5\gamma_{n,d_1,d_2}\sqrt{\log d_1}+C_6\mumax\sqrt{\frac{rd_1}{n}}.
\end{align*}
where $\Phi(x)$ denotes the c.d.f. of the standard normal distribution. 
\end{theorem}

By Theorem~\ref{thm:normal_mt}, if $ \mumax, \alpha_d,\kappa_0=O(1)$ and
\begin{align}
\max\bigg\{\frac{\|T\|_{\ell_1}}{\alpha_T\|T\|_{\subf}}\sqrt{\frac{r\log d_1}{d_2}}, \frac{\|T\|_{\ell_1}}{\alpha_T\|T\|_{\subf}}\cdot\frac{\sigma_{\xi}}{\lambda_r}\sqrt{\frac{rd_1^2d_2\log^2d_1}{n}},\gamma_{n,d_1,d_2}\sqrt{\log d_1}\bigg\}\to 0,\label{eq:asymp_cond}
\end{align}
then 
$$
\frac{\what m_T-m_T}{\sigma_{\xi}(\|TV\|_{\subf}^2+\|U^{\tran}T\|_{\subf}^2)^{1/2}\cdot\sqrt{d_1d_2/n}}\overset{{\rm d}}{\longrightarrow}\calN(0,1),
$$
as $n, d_1,d_2\to\infty$.

\subsection{Specific examples}\label{sec:examples}

We now consider several specific linear forms to further illustrate the implications of Theorem \ref{thm:normal_mt}.

\paragraph{Example 1:} 
As noted before, among the simplest linear forms are entries of $M$. In particular, $M(i,j)=\langle M, T\rangle$ with $T=e_{i}e_{j}^{\tran}$. It is clear that $\|T\|_{\ell_1}=\|T\|_{\subf}=1$ and Assumption~\ref{assump:T} is equivalent to
\begin{equation}\label{eq:T1_UV}
\|e_{i}^{\tran}U\|+\|e_{j}^{\tran}V\|\geq \alpha_T\sqrt{\frac{r}{d_1}}.
\end{equation}
Theorem~\ref{thm:normal_mt} immediately implies that 
$$
\frac{\what M(i,j)-M(i,j)}{(\|e_{i}^{\tran}U\|^2+\|e_{j}^{\tran}V\|^2)^{1/2}\cdot\sigma_{\xi}\sqrt{d_1d_2/n}} \overset{\rm d}{\longrightarrow}\calN(0,1),
$$
provided that
\begin{align}
\max\bigg\{\frac{\mumax^2}{\alpha_T}\sqrt{\frac{r\log d_1}{d_2}}, \ \frac{\kappa_0\mumax^2}{\alpha_T}\cdot\frac{\sigma_{\xi}}{\lambda_r}\sqrt{\frac{\alpha_drd_1^2d_2\log^2d_1}{n}},\ \gamma_{n,d_1,d_2}\sqrt{\log d_1}\bigg\}\to 0\label{eq:cond_T1}
\end{align}
as $n, d_1,d_2\to \infty$.

We can also infer from the entry-wise asymptotic normality that
\begin{equation}
\EE \|\what M-M\|_{\subf}^2=(1+o(1))\cdot\frac{\sigma_{\xi}^2rd_1d_2(d_1+d_2)}{n}.\label{eq:Mproj_rate}
\end{equation}
The mean squared error on the righthand side is sharply optimal and matches the minimax lower bound in \cite{koltchinskii2011nuclear}.

\paragraph{Example 2:} In the case when we want to compare $M(i,j_1)$ and $M(i,j_2)$, we can take $T=e_{i}e_{j_1}^{\tran}- e_{i}e_{j_2}^{\tran}$. Because $\|T\|_{\ell_1}/\|T\|_{\subf}=\sqrt{2}$, Assumption~\ref{assump:T} then becomes
\begin{equation}\label{eq:T2_UV}
\|T V\|_{\subf}^2+\|U^{\tran}T\|_{\subf}^2=2\| U^{\tran}e_{i}\|^2+\|V^{\tran}(e_{j_1}-e_{j_2})\|^2\geq \frac{2\alpha_T^2r}{d_1}.
\end{equation}
Theorem~\ref{thm:normal_mt} therefore implies that
$$
\frac{\big(\what M(i,j_1)-\what M(i,j_2)\big)-\big(M(i,j_1)-M(i,j_2)\big)}{\big(2\| U^{\tran}e_{i}\|^2+\|V^{\tran}(e_{j_1}-e_{j_2})\|^2\big)^{1/2}\cdot\sigma_{\xi}\sqrt{d_1d_2/n}} \overset{\rm d}{\longrightarrow}\calN(0,1),
$$
provided that
\begin{align}
\max\left\{\frac{\mumax^2}{\alpha_T}\sqrt{\frac{r\log d_1}{d_2}}, \ \frac{\kappa_0\mumax^2}{\alpha_T}\cdot\frac{\sigma_{\xi}}{\lambda_r}\sqrt{\frac{\alpha_drd_1^2d_2\log^2d_1}{n}},\ \gamma_{n,d_1,d_2}\sqrt{\log d_1}\right\}\to 0.\label{eq:cond_T2}
\end{align}

\paragraph{Example 3:}
More generally, we can consider the case when $T$ is sparse in that it has up to $s_0$ nonzero entries. By Cauchy-Schwartz inequality, $\|T\|_{\ell_1}/\|T\|_{\subf}\leq \sqrt{s_0}$ so that Assumption~\ref{assump:T} holds. By Theorem~\ref{thm:normal_mt}, 
$$
\frac{\what m_T-m_T}{\sigma_{\xi}(\|TV\|_{\subf}^2+\|U^{\tran}T\|_{\subf}^2)^{1/2}\cdot\sqrt{d_1d_2/n}}\overset{{\rm d}}{\longrightarrow}\calN(0,1),
$$
as long as
\begin{align}
\max\left\{\frac{\mumax^2}{\alpha_T}\sqrt{\frac{s_0r\log d_1}{d_2}}, \frac{\kappa_0\mumax^2}{\alpha_T}\cdot\frac{\sigma_{\xi}}{\lambda_r}\sqrt{\frac{s_0\alpha_drd_1^2d_2\log^2d_1}{n}},\gamma_{n,d_1,d_2}\sqrt{\log d_1}\right\}\to 0.\label{eq:asymp_spa}
\end{align}

It is of interest to consider the effect of alignment of $T$ with respect to the singular spaces of $M$. Note that
$$
\|T\|_{\subf}^2=\|U^{\tran}T\|_{\subf}^2+\|U^{\tran}_\perp T\|_{\subf}^2=\|TV\|_{\subf}^2+\|TV_\perp\|_{\subf}^2,
$$
where $U_\perp\in \OO^{d_1\times (d_1-r)}$ and $V_\perp\in \OO^{d_2\times (d_2-r)}$ are the basis of the orthogonal complement of $U$ and $V$ respectively. In the case that $T$ is not dominated by its projection onto $U_\perp$ or $V_\perp$ in that $\|U^{\tran}T\|_{\subf}+\|TV\|_{\subf}$ is of the same order as $\|T\|_{\subf}$, we can allow $T$ to have as many as $O(d_2)$ nonzero entries.

\section{Inferences about Linear Forms}\label{sec:conf_int}

The asymptotic normality of $\what m_T$ we showed in the previous section forms the basis for making inferences about $m_T$. To derive confidence intervals of or testing hypotheses about $m_T$, however, we need to also estimate the variance of $\what m_T$. To this end, we shall estimate the noise variance by
\begin{equation}
\what\sigma_{\xi}^2= \frac{1}{2n_0}\sum_{i=n_0+1}^{n}\big(Y_i-\langle \what M_1^{\supinit}, X_i\rangle\big)^2+\frac{1}{2n_0}\sum_{i=1}^{n_0}\big(Y_i-\langle \what M_2^{\supinit}, X_i\rangle \big)^2.
\end{equation}
and $\|TV\|_{\subf}^2+\|U^{\tran}T\|_{\subf}^2$ by
$$\what s_T^2:={1\over 2}\left(\|T\what V_1\|_{\subf}^2+\|\what U^{\tran}_1T\|_{\subf}^2+\|T\what V_2\|_{\subf}^2+\|\what U^{\tran}_2T\|_{\subf}^2\right).$$
The following theorem shows that the asymptotic normality remains valid if we replace the variance of $\what m_T$ with these estimates:
\begin{theorem}\label{thm:conf_int}
Under Assumptions~\ref{assump:init_entry}-\ref{assump:T},  if  $n\geq C\mumax^2rd_1\log d_1$ for some absolute constant $C>0$ and 
$$
\max\bigg\{\frac{\mumax^2\|T\|_{\ell_1}}{\alpha_T\|T\|_{\subf}}\sqrt{\frac{r\log d_1}{d_2}}, \frac{\kappa_0\mumax^2\|T\|_{\ell_1}}{\alpha_T\|T\|_{\subf}}\cdot \frac{\sigma_{\xi}}{\lambda_r}\sqrt{\frac{\alpha_drd_1^2d_2\log d_1^2}{n}}, \gamma_{n,d_1,d_2}\sqrt{\log d_1}\bigg\}\to 0,
$$
then 
$$
\frac{\what m_T-m_T}{\what\sigma_{\xi}\what s_T\cdot\sqrt{d_1d_2/n}} \overset{\rm d}{\longrightarrow}\calN(0,1),
$$
as $n,d_1,d_2\to\infty$. 
\end{theorem}

Theorem \ref{thm:conf_int} immediately allows for constructing confidence intervals for $m_T$. More specifically, we can define the $100(1-\theta)\%$-th confidence interval as
\begin{align}\label{eq:CI_alpha}
\what{{\rm CI}}_{\theta, T}=\bigg[\what m_T-z_{\theta/2}\cdot\what{\sigma}_{\xi}\what s_T\sqrt{\frac{d_1d_2}{n}},\ \ \what m_T+z_{\theta/2}\cdot\what{\sigma}_{\xi}\what s_T\sqrt{\frac{d_1d_2}{n}}\bigg]
\end{align}
for any $\theta\in(0,1)$, where $z_{\theta}=\Phi^{-1}(1-\theta)$ is the upper $\theta$ quantile of the standard normal. In light of Theorem \ref{thm:conf_int}, we have
$$
\lim_{n,d_1,d_2\to\infty}\ \PP\big(m_T\in \what{{\rm CI}}_{\theta,T}\big)=1-\theta,
$$
for any $\theta\in(0,1)$.

Similarly, we can also consider using Theorem \ref{thm:conf_int} for the purpose of hypothesis test. Consider, for example, testing linear hypothesis
$$
H_0: \langle M, T\rangle =0 \qquad {\rm against }\qquad H_1:\langle M, T\rangle \neq 0.
$$
Then we can proceed to reject $H_0$ if $|\what z|>z_{\theta/2}$ and accept $H_0$ otherwise, where
$$
\what z=\frac{\what m_T}{\what\sigma_{\xi}\what s_T\cdot\sqrt{d_1d_2/n}}.
$$
Following Theorem \ref{thm:conf_int}, this is a test with asymptotic level $\theta$. For example, in the particular case of comparing two entries of $M$:
\begin{equation}\label{eq:H0}
H_0: M(i,j_1)=M(i,j_2)\quad {\rm v.s.}\quad H_1: M(i,j_1)> M(i,j_2),
\end{equation}
the test statistic can be expressed as
\begin{align*}
\what z=\frac{\sqrt{2}(\what M(i,j_1)-\what M (i,j_2))}{\what\sigma_{\xi}\big(\|\what V_1^{\tran}(e_{j_2}-e_{j_1})\|_{\subf}^2+2\|\what U_1^{\tran}e_{i}\|_{\subf}^2+\|\what V_2^{\tran}(e_{j_2}-e_{j_1})\|_{\subf}^2+2\|\what U_2^{\tran}e_{i}\|_{\subf}^2\big)^{1/2}\sqrt{d_1d_2/n}}
\end{align*}
and we shall proceed to reject the null hypothesis if and only if $\what z>z_{\theta}$ to account for the one-sided alternative.

\section{Initial Estimate}\label{sec:init_ncx}
Thus far, our development has assumed a generic max-norm consistent matrix estimate as initial estimator. For concreteness, we now introduce a rotation calibrated gradient descent algorithm on Grassmannians which, under mild conditions, produces such an estimate.

Any rank $r$ matrix of dimension $d_1\times d_2$ can be written as $UGV^{\tran}$ where $U\in\OO^{d_1\times r}$, $V\in\OO^{d_2\times r}$ and $G\in \RR^{r\times r}$. The loss of the triplet $(U, G, V)$ over $\mathfrak{D}$ is given by
\begin{align}\label{eq:loss_L}
L\big(\mathfrak{D}, (U,G,V)\big)=\sum_{(X,Y)\in\mathfrak{D}}\big(Y-\langle UGV^{\tran},X\rangle\big)^2.
\end{align}
Given $(U, V)$, we can easily minimize (\ref{eq:loss_L}) to solve for $G$. This allows us to reduce the problem of minimizing \eqref{eq:loss_L} to a minimization over the product space of two Grassmannians ${\sf Gr}(d_1, r)\times{\sf Gr}(d_2, r)$ as ${\sf Gr}(d, r)=\OO^{d_1\times r}/\OO^{r\times r}$. In particular we can do so via a rotation calibrated gradient descent algorithm on Grassmannians as detailed in Algorithm~\ref{algo:GD} where, for simplicity, we resort to data-splitting. It is plausible that a more elaborative analysis via the leave-one-out (LOO) framework introduced by \cite{ma2017implicit} can be applied to show that our algorithm continues to produce estimates of similar quality without data-splitting, as we observe empirically. An actual proof however is likely much more involved under our setting. For brevity, we opted here for data-splitting.



Let $m=C_1\ceil{\log (d_1+d_2)}$ for some positive integer $C_1\geq 1$. We shall partition the data $\mathfrak{D}=\{(X_i, Y_i)\}_{i=1}^n$ into $2m$ subsets:
$$
\mathfrak{D}_t=\big\{(X_j,Y_j)\big\}_{j=(t-1)N_0+1}^{tN_0},\quad \forall\ t=1,\cdots,2m
$$
where, without loss of generality, we assumed $n=2mN_0$ for some positive integer $N_0$. 

\begin{algorithm}[H]
\caption{Rotation Calibrated Gradient descent on Grassmannians}
\label{algo:GD}
\begin{algorithmic}[2]
\State Let $\what U^{(1)}$ and $\what V^{(1)}$ be the top-$r$ left and right singular vectors of $d_1d_2N_0^{-1}\sum_{j\in\mathfrak{D_1}}Y_jX_j$.
\State Compute $\what{G}^{(1)}=\arg\min_{G\in\RR^{r\times r}} L(\mathfrak{D}_2, (\what U^{(1)}, G, \what {V}^{(1)}))$ and its SVD $\what G^{(1)}=\what L_G^{(1)}\what\Lambda^{(1)}\what R_G^{(1)\tran}$.
\For{$t=1,2,3,\cdots,m-1$}
\State Update by rotation calibrated gradient descent
$$
\what U^{(t+0.5)}=\what U^{(t)}\what L_G^{(t)}-\eta\cdot \frac{d_1d_2}{N_0}\sum_{j\in\mathfrak{D}_{2t+1}}\big(\langle \what U^{(t)}\what G^{(t)}\what V^{(t)}, X_j\rangle-Y_j\big)X_j\what V^{(t)}\what R_G^{(t)}(\what \Lambda^{(t)})^{-1}
$$
$$
\what V^{(t+0.5)}=\what V^{(t)}\what R_G^{(t)}-\eta\cdot \frac{d_1d_2}{N_0}\sum_{j\in\mathfrak{D}_{2t+1}}\big(\langle \what U^{(t)}\what G^{(t)}\what V^{(t)}, X_j\rangle-Y_j\big)X_j^{\tran}\what U^{(t)}\what L_G^{(t)}(\what \Lambda^{(t)})^{-1}
$$
\State Compute the top-$r$ left singular vectors
$$
\what U^{(t+1)}={\rm SVD}(\what U^{(t+0.5)})\quad {\rm and}\quad \what V^{(t+1)}={\rm SVD}(\what V^{(t+0.5)})
$$
\State Compute $\what G^{(t+1)}$ by 
$$
\what{G}^{(t+1)}=\underset{G\in\RR^{r\times r}}{\arg\min}\ L\big(\mathfrak{D}_{2t+2}, (\what U^{(t+1)},G,\what V^{(t+1)})\big)
\textrm{ and its SVD }
 \what G^{(t+1)}=\what L_G^{(t+1)}\what\Lambda^{(t+1)}\what R_G^{(t+1)\tran}
$$
\EndFor
\State Output: $(\what U^{(m)}, \what G^{(m)}, \what V^{(m)})$ and $\what M^{(m)}=\what U^{(m)}\what G^{(m)}(\what V^{(m)})^{\tran}$.
\end{algorithmic}
\end{algorithm}

The algorithm presented here is similar in spirit to those developed earlier by \cite{keshavan2010matrix_a, keshavan2010matrix_b, xia2017polynomial}. A key difference is that we introduce an explicit rule of gradient descent update where each iteration on Grassmannians is calibrated with orthogonal rotations. The rotation calibrations are necessary to guarantee the contraction property for the $(2,\max)$-norm accuracy of empirical singular vectors. Indeed, we show that the algorithm converges geometrically with constant stepsizes.

To this end, write
$$
\what O_U^{(1)}=\arg\min_{O\in\OO^{r\times r}}\|\what U^{(1)}-UO\|\quad {\rm and}\quad \what O_V^{(1)}=\arg\min_{O\in\OO^{r\times r}}\|\what V^{(1)}-VO\|
$$
and, for all $t=1,\cdots,m-1$, denote the SVDs
$$
\what U^{(t+0.5)}=\what U^{(t+1)}\what\Sigma_U^{(t+1)}\what K_U^{(t+1)\tran}\quad {\rm and}\quad \what V^{(t+0.5)}=\what V^{(t+1)}\what\Sigma_V^{(t+1)}\what K_V^{(t+1)\tran}.
$$
For all $t=1,\cdots,m-1$, define the orthogonal matrices
$$
\what O_U^{(t+1)}=\what O_U^{(t)}\what L_G^{(t)}\what K_U^{(t+1)}\quad {\rm and}\quad \what O_V^{(t+1)}=\what O_V^{(t)}\what R_G^{(t)}\what K_V^{(t+1)}.
$$
Then we have
\begin{theorem}\label{thm:gd}
Under Assumptions~\ref{assump:incoh} and \ref{assump:noise}, if $\eta\in[0.25, 0.75]$ and
$$
n\geq C_1\alpha_d\kappa_0^6\mumax^6 r^3d_1\log^2d_1\quad {\rm and}\quad C_2\kappa_0^2\mumax\frac{\sigma_{\xi}}{\lambda_r}\cdot \sqrt{\frac{\alpha_d rd_1^2d_2\log^2d_1}{n}}\leq 1
$$
for some large enough constants $C_1,C_2>0$, then for all $t=1,\cdots,m-1$, with probability at least $1-4md_1^{-2}$,
\begin{align*}
\big\|\what U^{(t+1)}-U&\what O_U^{(t+1)}\big\|_{\submaxx}+\big\|\what V^{(t+1)}-V\what O_V^{(t+1)}\big\|_{\submaxx}\leq C_3\eta\frac{\sigma_{\xi}}{\lambda_r}\sqrt{\frac{rd_1d_2\log^2d_1}{n}}\\
& +\Big(1-\frac{2\eta}{3}\Big)\cdot \big(\|\what U^{(t)}-U\what O_U^{(t)}\|_{\submaxx}+\|\what V^{(t)}-V\what O_V^{(t)}\|_{\submaxx}\big)
\end{align*}
where $C_3>0$ is an absolute constant. Moreover, if in addition $\|M\|_{\submax}/\sigma_{\xi}\leq d_1^{C_4}$ for some constant $C_4>0$, then, by setting $m=2\ceil{C_4\log d_1}$ and $\eta=0.75$, with probability at least $1-C_5d_1^{-2}\log d_1$,
$$
\big\|\what M^{(m)}-M \big\|_{\submax}\leq C_6\mumax\kappa_0 \sigma_{\xi}\sqrt{\frac{r^2d_1\log^2d_1}{n}}
$$
for some absolute constants $C_5,C_6>0$. 
\end{theorem}

We can then apply Algorithm~\ref{algo:GD} to produce initial estimates suitable for inferences about linear forms of $M$. With this particular choice of initial estimate, Assumption~\ref{assump:init_entry} is satisfied with 
$$
\gamma_{n,d_1,d_2}=\mumax\kappa_0\sqrt{\frac{r^2d_1\log^2d_1}{n}}
$$ 
when the sample size $n\geq C_1\alpha_d \kappa_0^6\mumax^6r^3d_1\log^2 d_1$. We note that this sample size requirement in general is not optimal and the extra logarithmic factor is due to data splitting. As this is not the main focus of the current work, no attempt is made here to further improve it.

\section{Numerical Experiments}\label{sec:numerical}
We now present several sets of numerical studies to further illustrate the practical merits of the proposed methodology, and complement our theoretical developments.

\subsection{Simulations}
We first consider several sets of simulation studies. Throughout the simulations, the true matrix $M$ has rank $r=3$ and dimension $d_1=d_2=d=2000$. $M$'s singular values were set to be $\lambda_i=d$ for $i=1,2,3$. In addition, $M$'s singular vectors were generated from the SVD of $d\times r$ Rademacher random matrices. The noise standard deviation was set at $\sigma_{\xi}=0.6$. 

\begin{figure}
\centering     
\includegraphics[width=.40\textwidth]{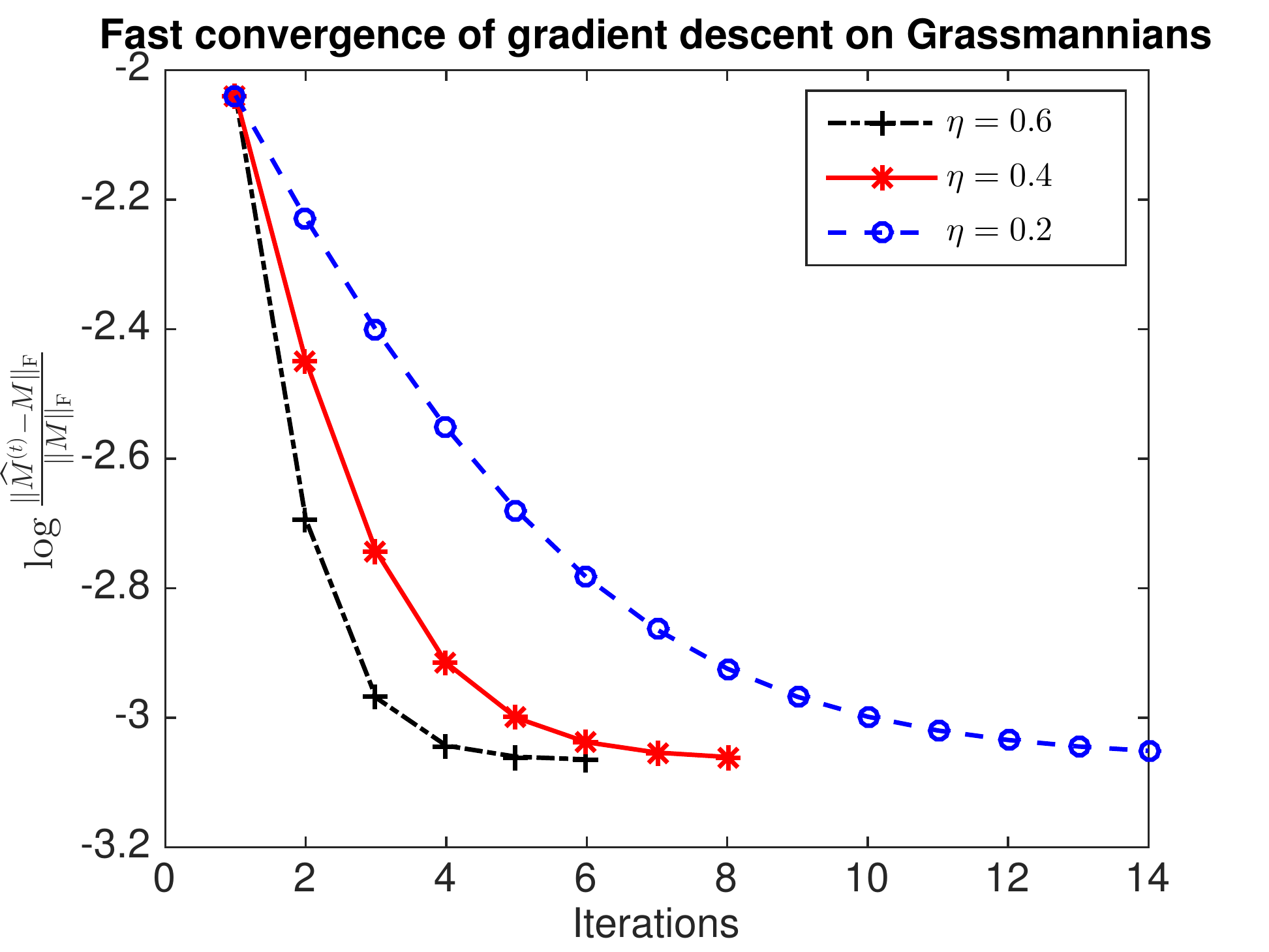}
\includegraphics[width=.40\textwidth]{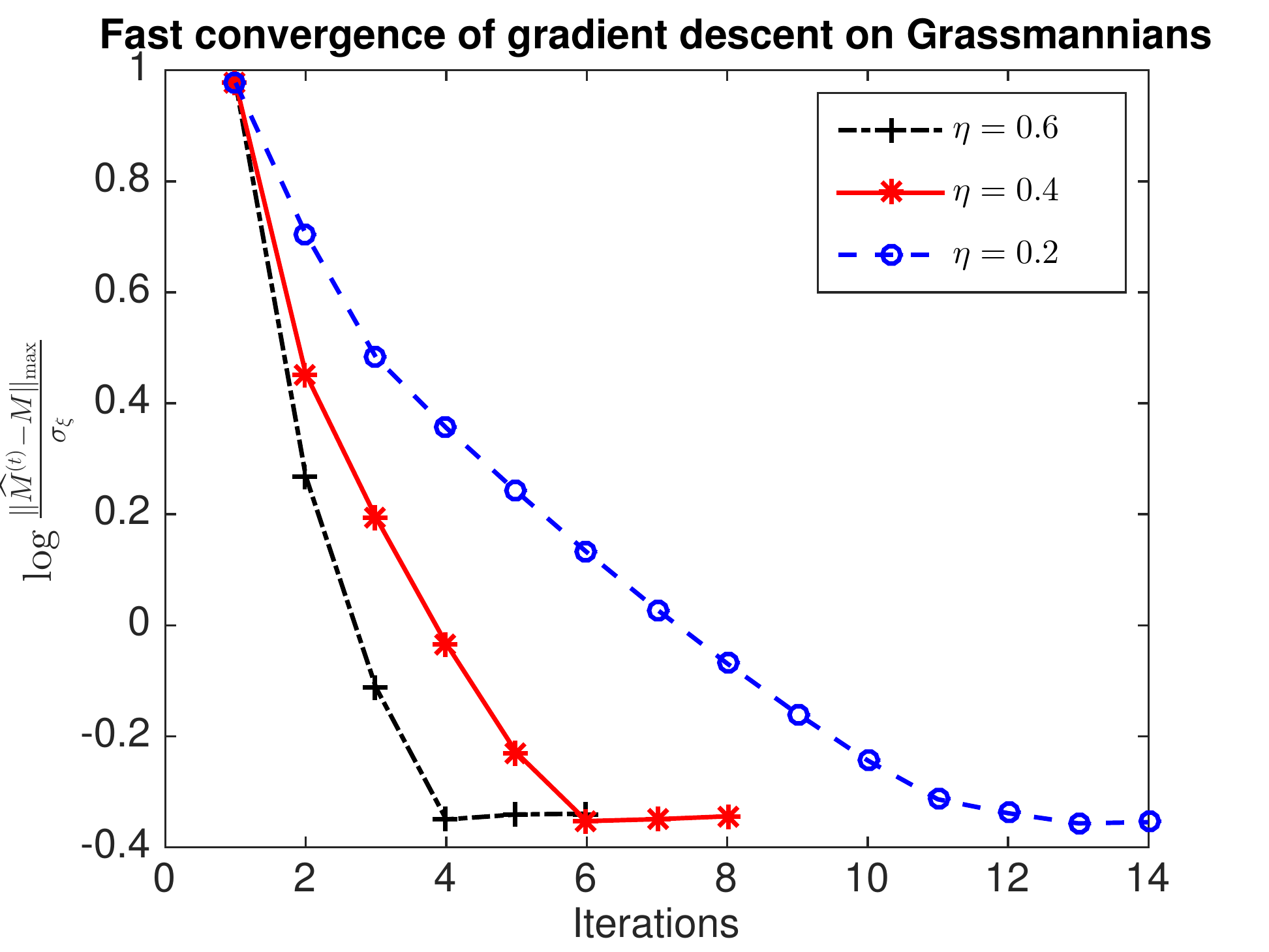}
\caption{Convergence of Algorithm~\ref{algo:GD} in relative matrix Frobenius norm and the max-norm, with respect to step size $\eta$ and the number of iterations. The parameters are $d_1=d_2=d=2000, r=3, \lambda_i=d, \sigma_{\xi}=0.6$ and $U,V$ are generated from the SVD of $d\times r$ Rademacher random matrices. The sample size is  $n=4r^2d\log(d)$ and the noise is Gaussian.}
\label{fig:conv}
\end{figure}

First, we show the convergence performance of the proposed Algorithm~\ref{algo:GD} where both Frobenius norm and max-norm convergence rates are recorded. Even though the algorithm we presented in the previous section uses sample splitting for technical convenience, in the simulation, we did not split the sample. Figure~\ref{fig:conv} shows a typical realization under Gaussian noise, which suggest the fast convergence of Algorithm~\ref{algo:GD}. In particular, $\log\frac{\|\what M^{\supinit}-M\|_{\submax}}{\sigma_{\xi}}$ becomes negative after $3$ iterations when the stepsize is $\eta=0.6$. Recall that our double-sample debiasing approach  requires $\|\what M^{\supinit}-M\|_{\submax}=o_P(\sigma_{\xi})$  for the initial estimates, i.e.,  $\what M_1^{\supinit}, \what M_2^{\supinit}$ in Assumption~\ref{assump:init_entry}. 

\begin{figure}
\centering     
\includegraphics[width=.40\textwidth]{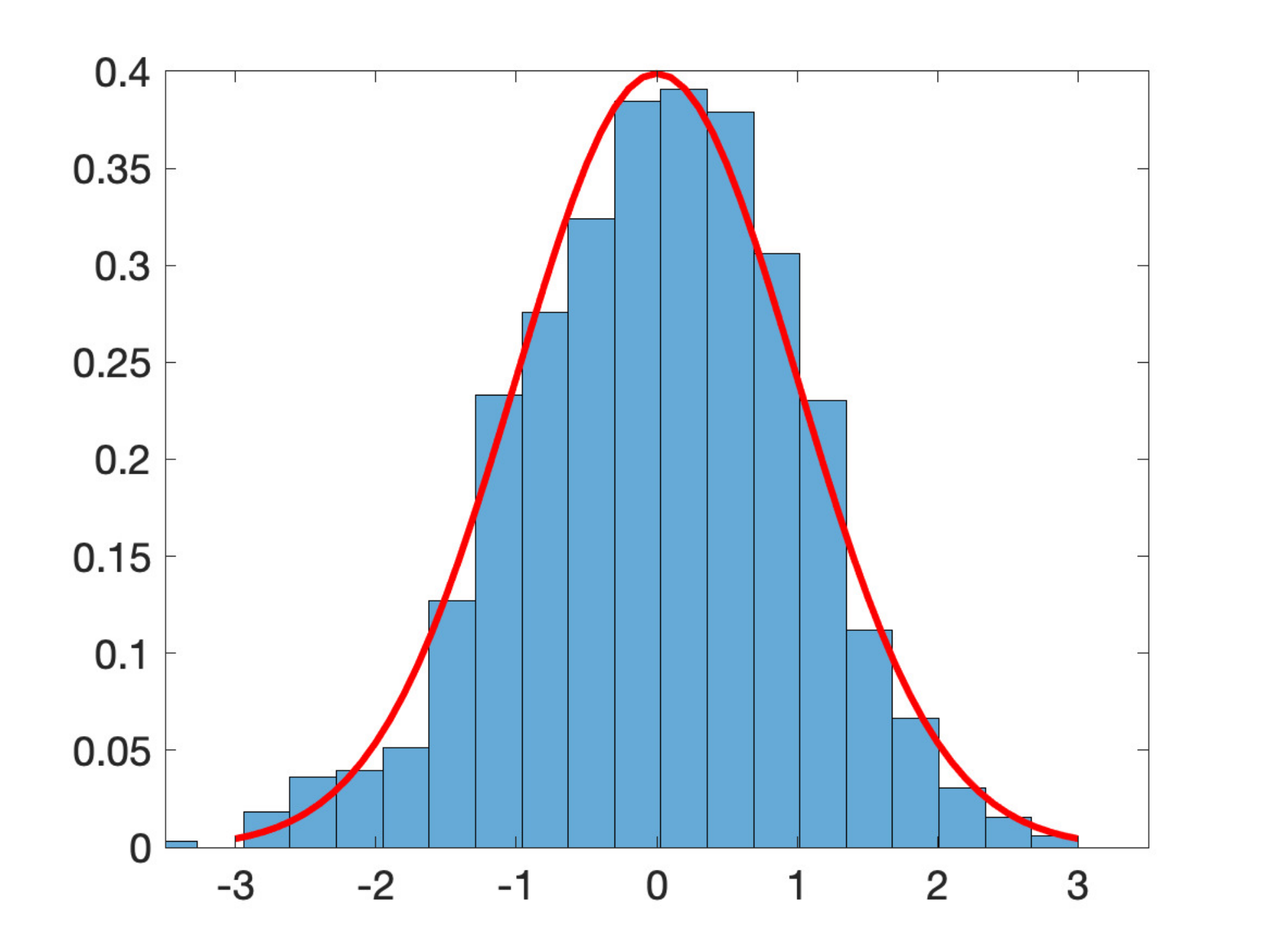}
\includegraphics[width=.40\textwidth]{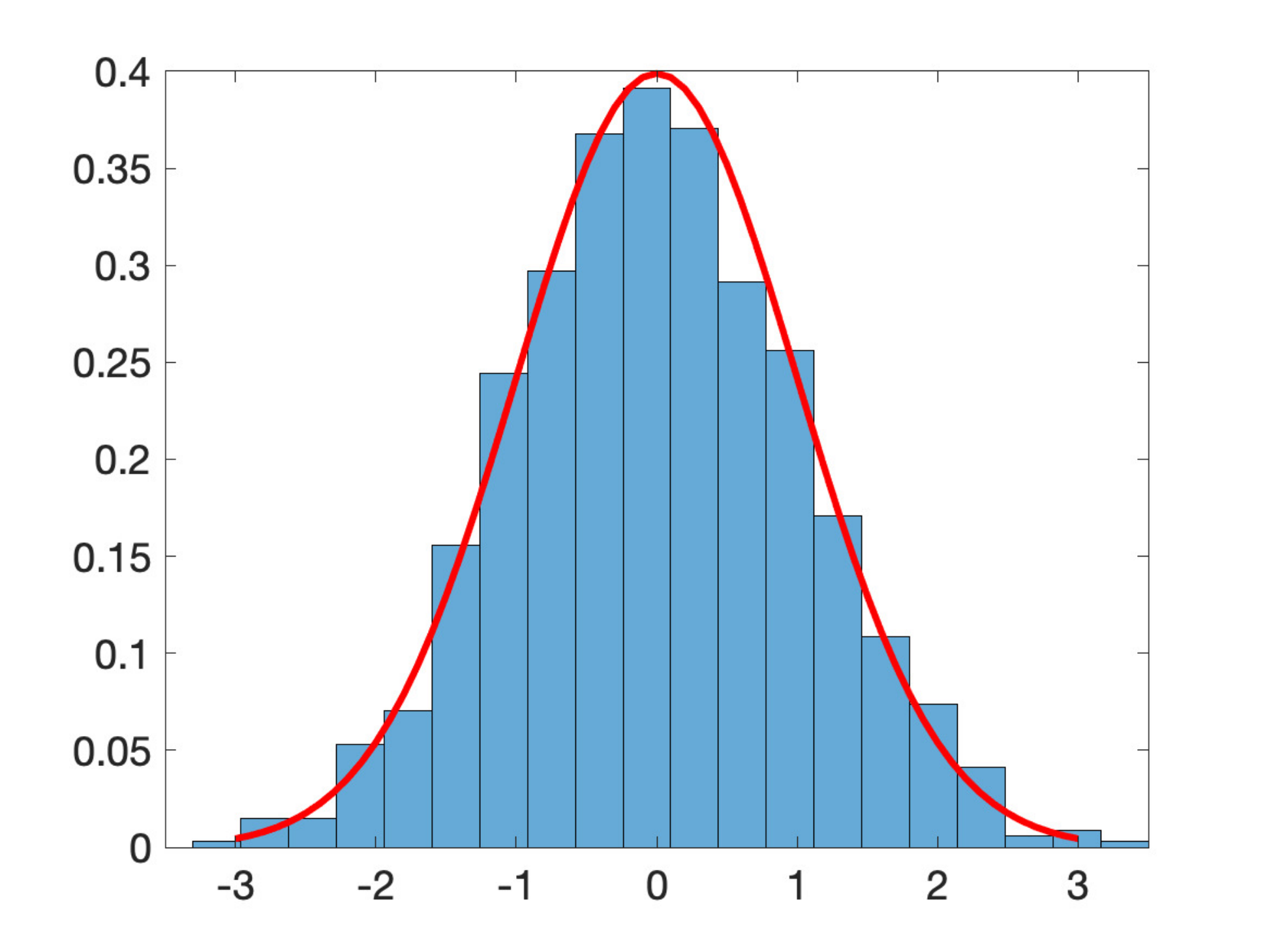}
\includegraphics[width=.40\textwidth]{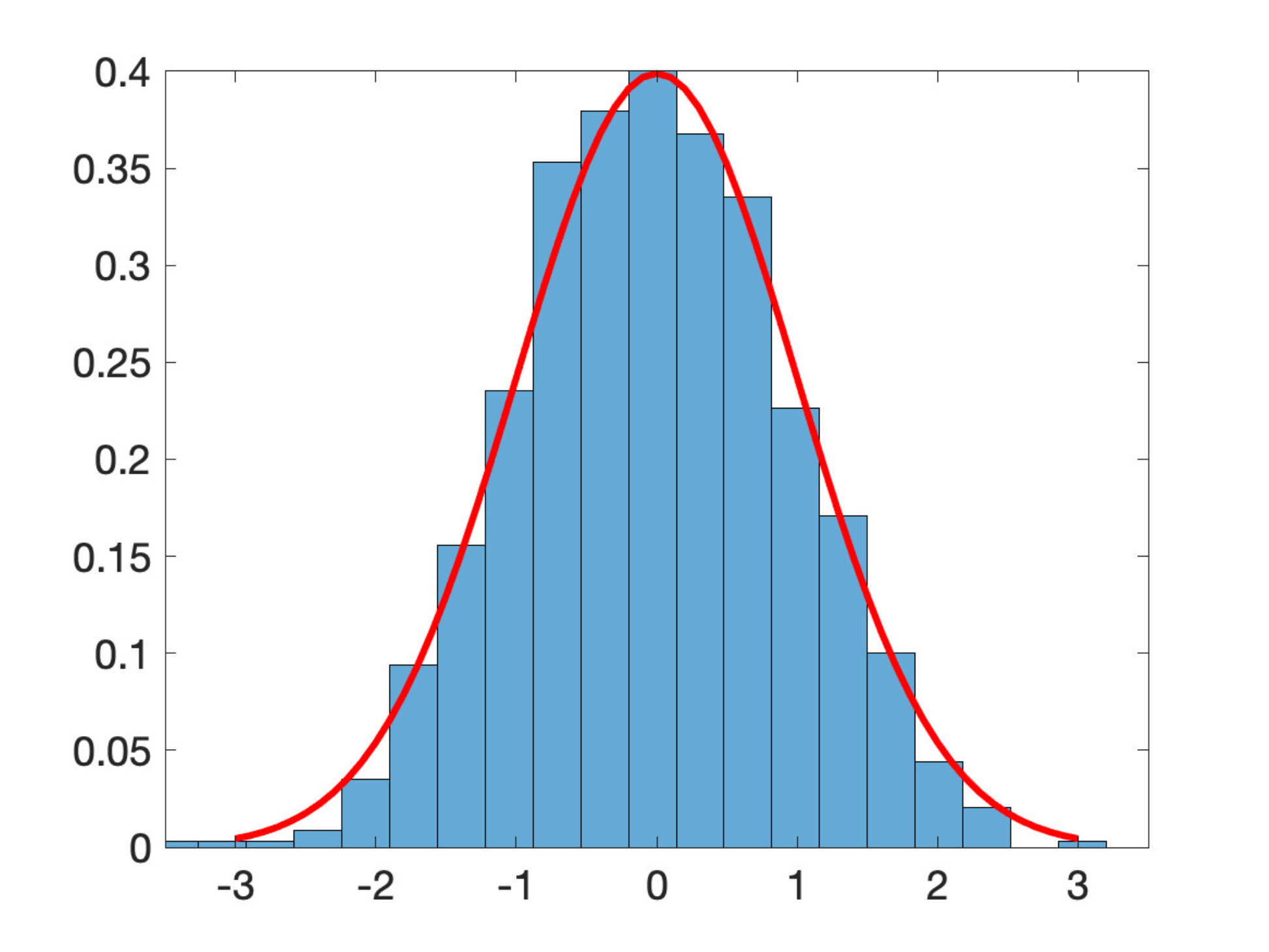}
\includegraphics[width=.40\textwidth]{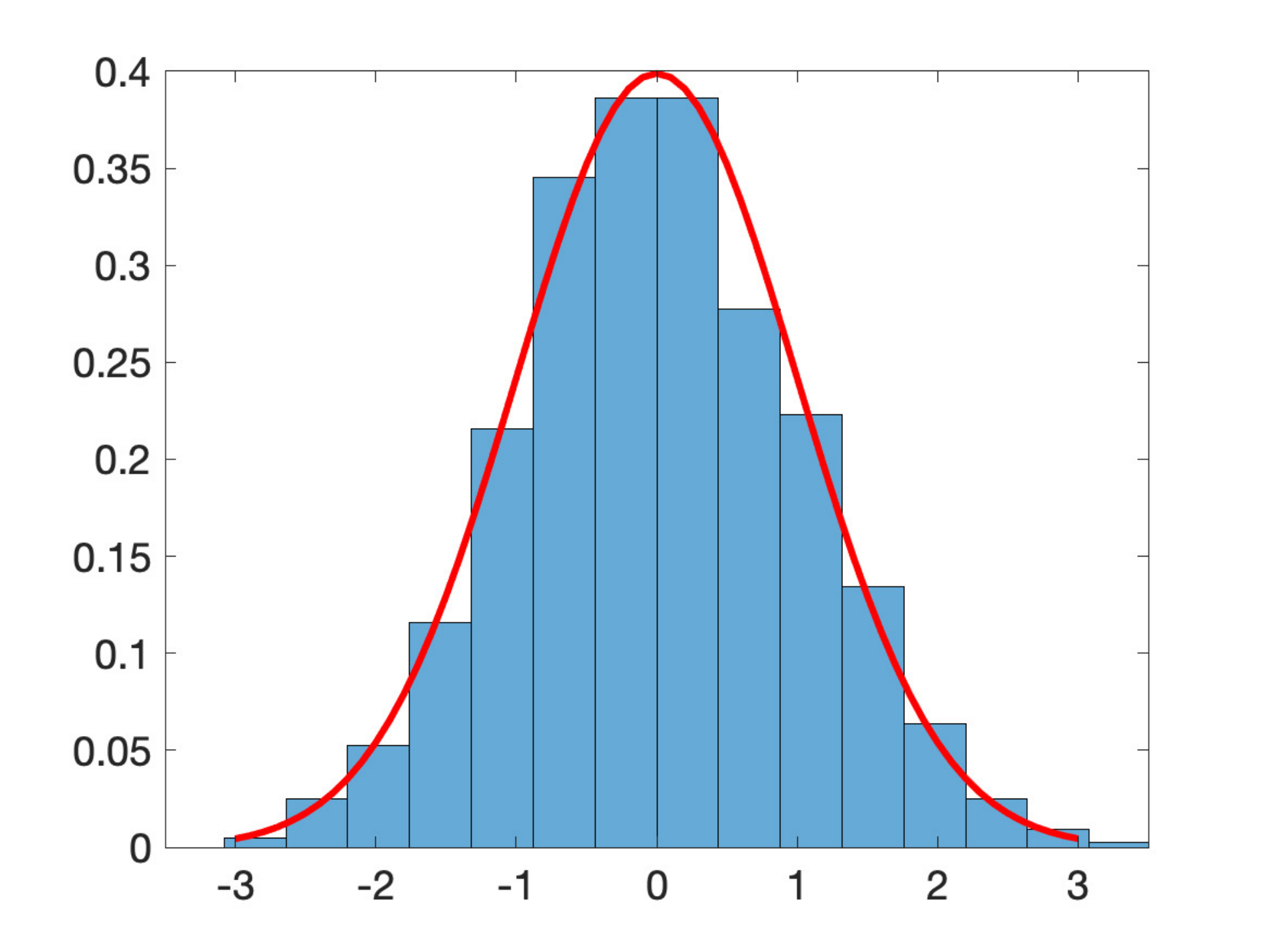}
\caption{Normal approximation of $\frac{\what m_T-m_T}{\what\sigma_{\xi}\what s_T\sqrt{d_1d_2/n}}$. The parameters are $d_1=d_2=d=2000, r=3, \lambda_i=d, \sigma_{\xi}=0.6$ and $U,V$ are generated from the SVD of $d\times r$ Rademacher random matrices. The sample size is  $n=4r^2d\log(d)$ for the top two and $n=5r^2d\log(d)$ for bottom two. The noise is Gaussian. Each density histogram is based on 1000 independent simulations and the red curve represents the p.d.f. of standard normal distributions.  Top left: $T=e_1e_1^{\tran}$, top right: $T=e_1e_1^{\tran}-e_1e_2^{\tran}$. Bottom left: $T=e_1e_1^{\tran}-e_1e_2^{\tran}+e_2e_1^{\tran}$, bottom right: $T=e_1e_1^{\tran}-e_1e_2^{\tran}+2e_{2}e_1^{\tran}+3e_{2}e_2^{\tran}$.}
\label{fig:na}
\end{figure}

Next, we investigate how the proposed inference tools behave under Gaussian noise and for four different linear forms corresponding to $T_1=e_1e_1^{\tran}$, $T_2=e_1e_1^{\tran}-e_1e_2^{\tran}$, $T_3=e_1e_1^{\tran}-e_1e_2^{\tran}+e_{2}e_1^{\tran}$ and 
$$
T_4=e_1e_1^{\tran}-e_1e_2^{\tran}+2e_{2}e_1^{\tran}+3e_{2}e_2^{\tran}.
$$
For each $T$, we drew the density histogram of $(\what m_T-m_T)/\big(\what\sigma_{\xi}\what s_T\sqrt{d_1d_2/n}\big)$ based on $1000$ independent simulation runs. The density histograms are displayed in Figure~\ref{fig:na} where the red curve represents the p.d.f. of standard normal distributions. The sample size was  $n=4r^2d\log(d)$ for $T_1, T_2$ and $n=5r^2d\log(d)$ for $T_3, T_4$. The empirical observation agrees fairly well with our theoretical results.

\begin{figure}
\centering     
\includegraphics[width=.40\textwidth]{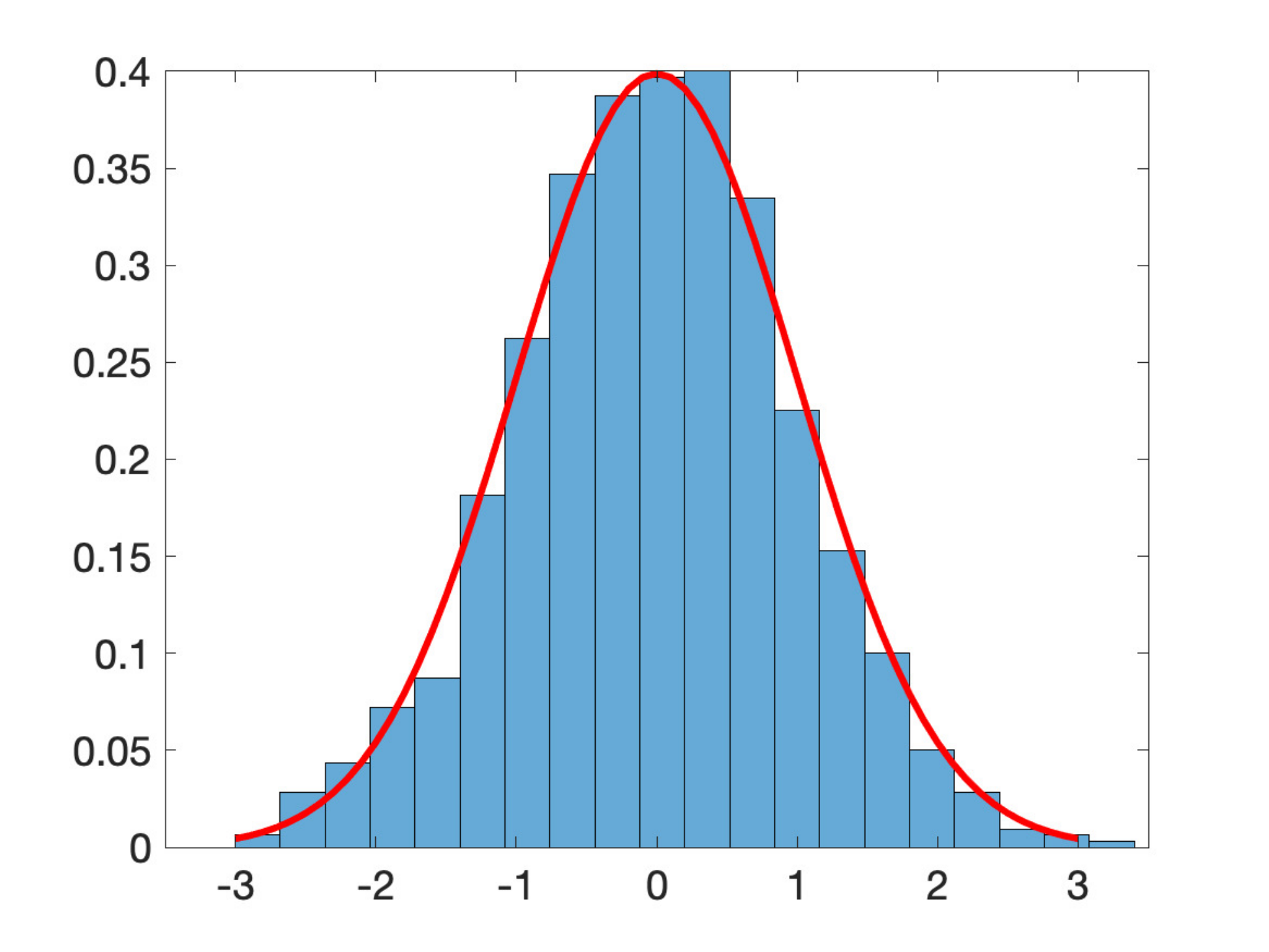}
\includegraphics[width=.40\textwidth]{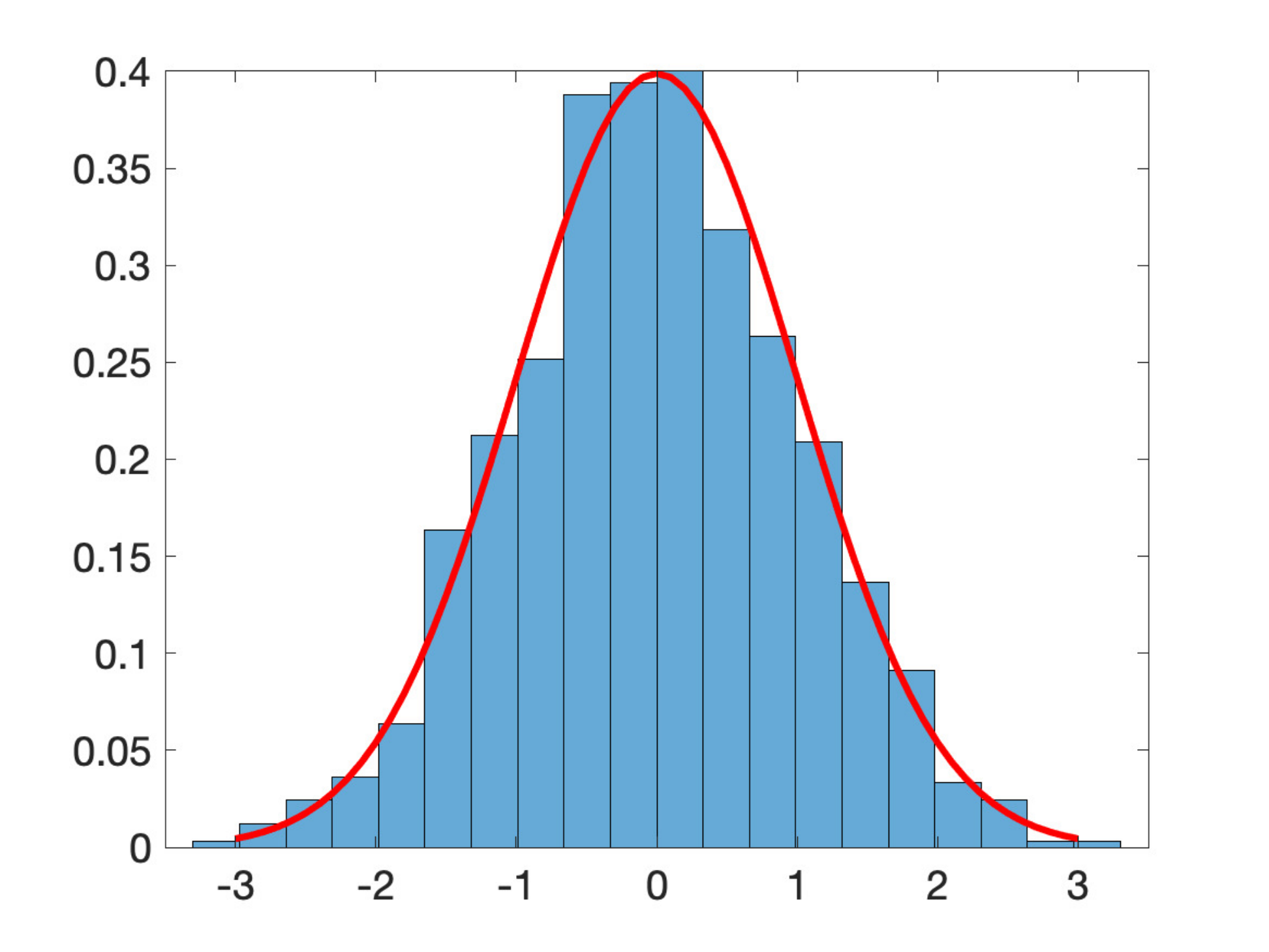}
\includegraphics[width=.40\textwidth]{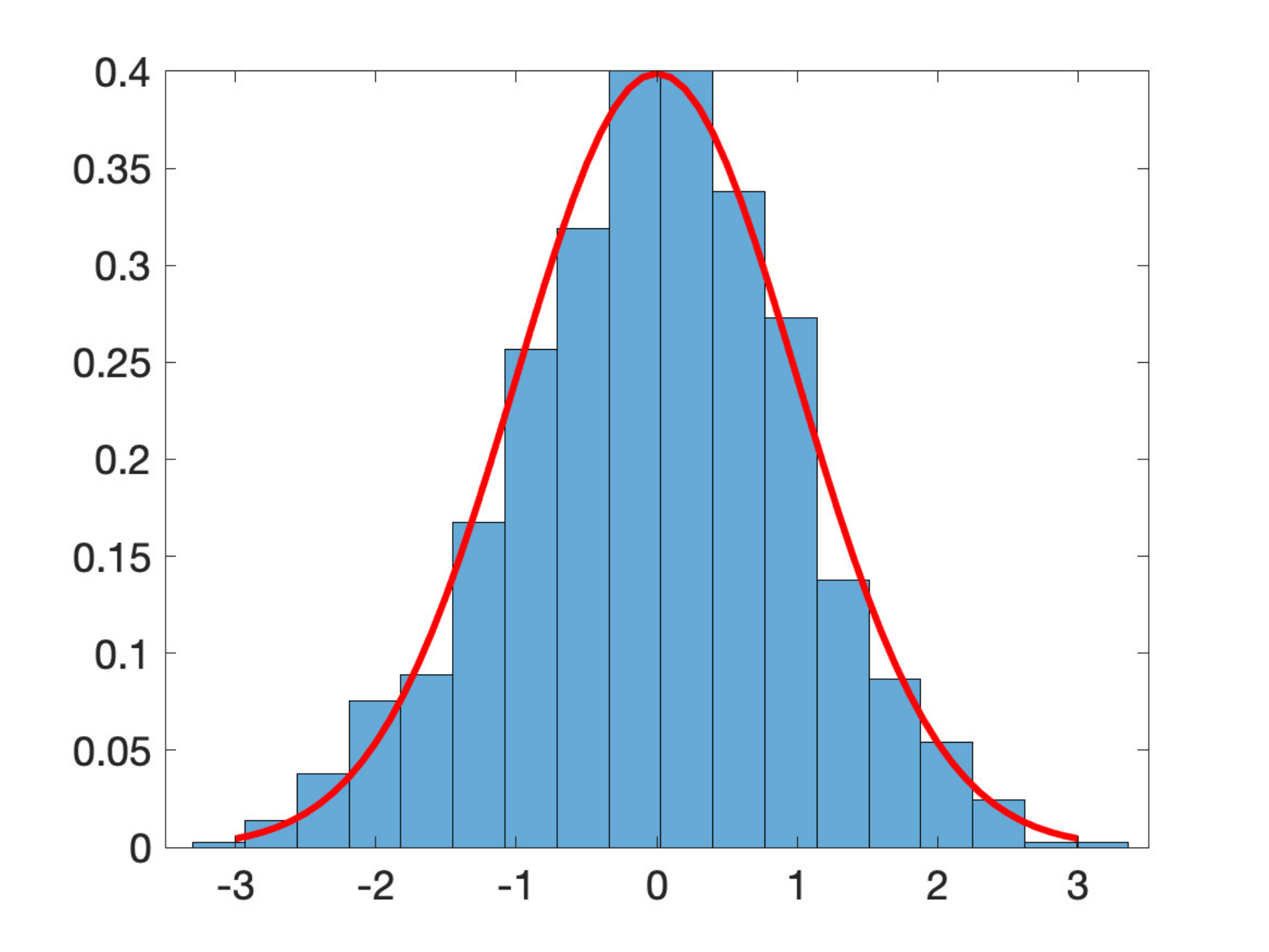}
\includegraphics[width=.40\textwidth]{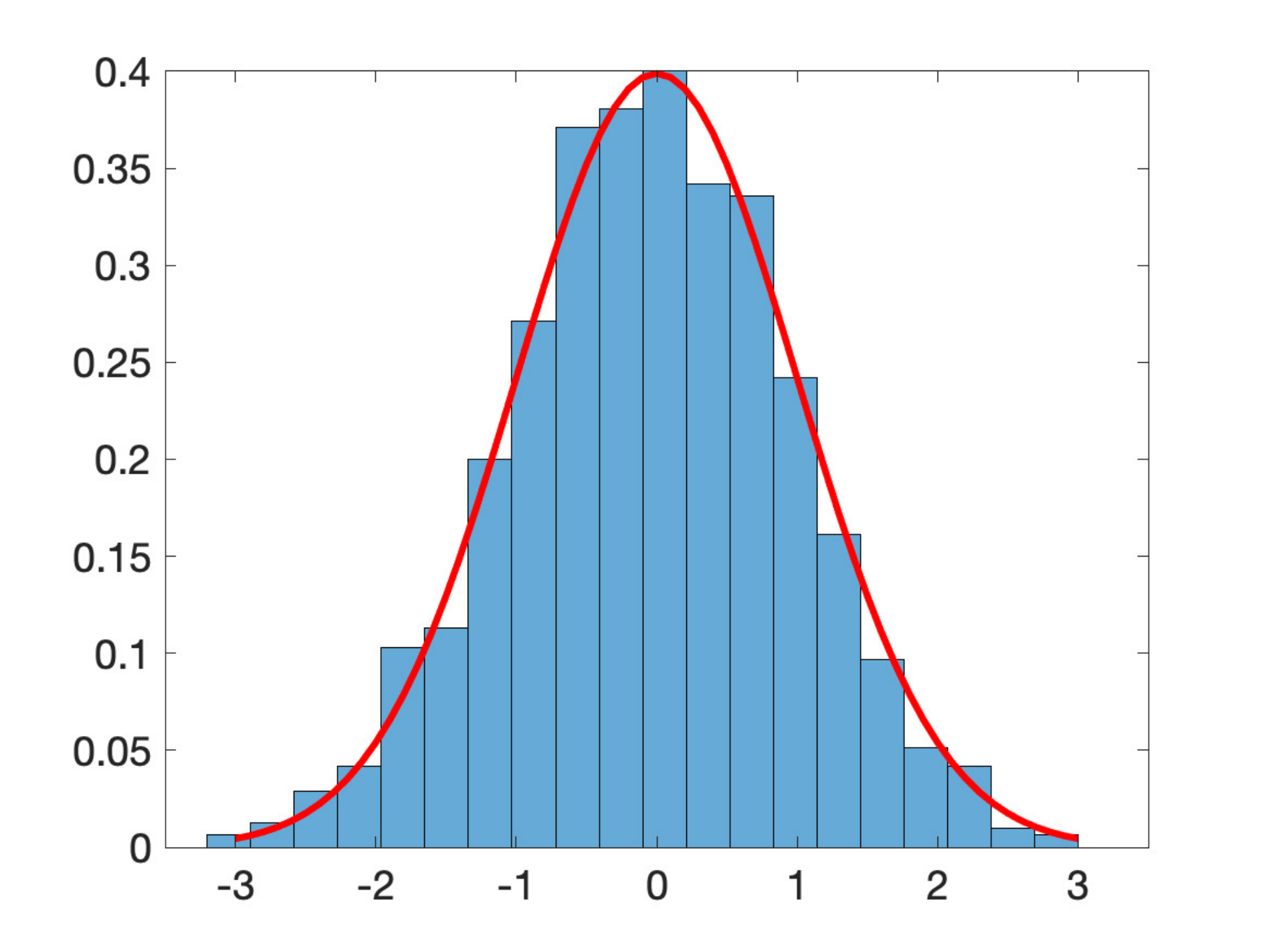}
\caption{Normal approximation of $\frac{\what m_T-m_T}{\what\sigma_{\xi}\what s_T\sqrt{d_1d_2/n}}$. The parameters are $d_1=d_2=d=2000, r=3, \lambda_i=d, \sigma_{\xi}=0.6$ and $U,V$ are generated from the SVD of $d\times r$ Rademacher random matrices. The sample size is  $n=4r^2d\log(d)$ for the top two and $n=5r^2d\log(d)$ for the bottom two. The non-Gaussian noise $(\xi/\sqrt{3}\sigma_{\xi})\in {\rm Unif}([-1,1])$. Each density histogram is based on 1000 independent simulations and the red curve represents the p.d.f. of standard normal distributions. Top left: $T=e_1e_1^{\tran}$, top right: $T=e_1e_1^{\tran}-e_1e_2^{\tran}$. Bottom left: $T=e_1e_1^{\tran}-e_1e_2^{\tran}+e_2e_1^{\tran}$, bottom right: $T=e_1e_1^{\tran}-e_1e_2^{\tran}+2e_{2}e_1^{\tran}+3e_{2}e_2^{\tran}$. }
\label{fig:na_U}
\end{figure}

Finally, we examine the performance of the proposed approach under non-Gaussian noise. In particular, we repeated the last set of experiments with noise $(\xi/\sqrt{3}\sigma_{\xi})\in {\rm Unif}([-1,1])$. The density histograms are displayed in Figure~\ref{fig:na_U} where the red curve represents the p.d.f. of standard normal distributions. Again the empirical evidences support the asymptotic normality of the proposed statistic.

\subsection{Real-world data examples}
We now turn our attention to two real-world data examples -- the Jester and MovieLens datasets. The Jester dataset contains ratings of $100$ jokes from $\sim70K$ users (\cite{goldberg2001eigentaste}). The dataset consists of 3 subsets of data with different characteristics as summarized in Table~\ref{tb:jester-movielens}. For each subset, the numbers of ratings of all users are equal. MovieLens was a recommender system created by GroupLens that recommends movies for its users. We use three datasets released by MovieLens \citep{harper2016movielens} whose details are summarized in Table~\ref{tb:jester-movielens}. In these three datasets, each user rates at least $20$ movies. 


\begin{table}
    \caption{\label{tb:jester-movielens}Summary of Datasets}
   \centering
  \fbox{%
\begin{tabular}{*{5}{c}}
\hline
 Dataset & \#users & \#jokes &  \#ratings per user&rating values    \\
 \cline{1-1}
Jester-1 &  24983  & 100 & $29$&[-10, 10]  \\
Jester-2  & 23500   & 100  & $34$&[-10, 10] \\     
Jester-3 & 24938&100&$ 14$&[-10, 10]\\   
\hline
Dataset & \#users & \#movies & total \#ratings&rating values    \\
\cline{1-1}
ml-100k &  943  & 1682 & $\sim 10^5$&\{1,2,3,4,5\}   \\
ml-1m  & 6040   & 3952  & $\sim10^6$&\{1,2,3,4,5\} \\     
ml-10m& 71567&10681&$\sim10^7$&$\{0.5,1.0,\cdots,4.5,5.0\}$\\   
\hline
\end{tabular}}
\end{table}

\begin{figure}[htbp]
\centering
\includegraphics[width=.32\textwidth]{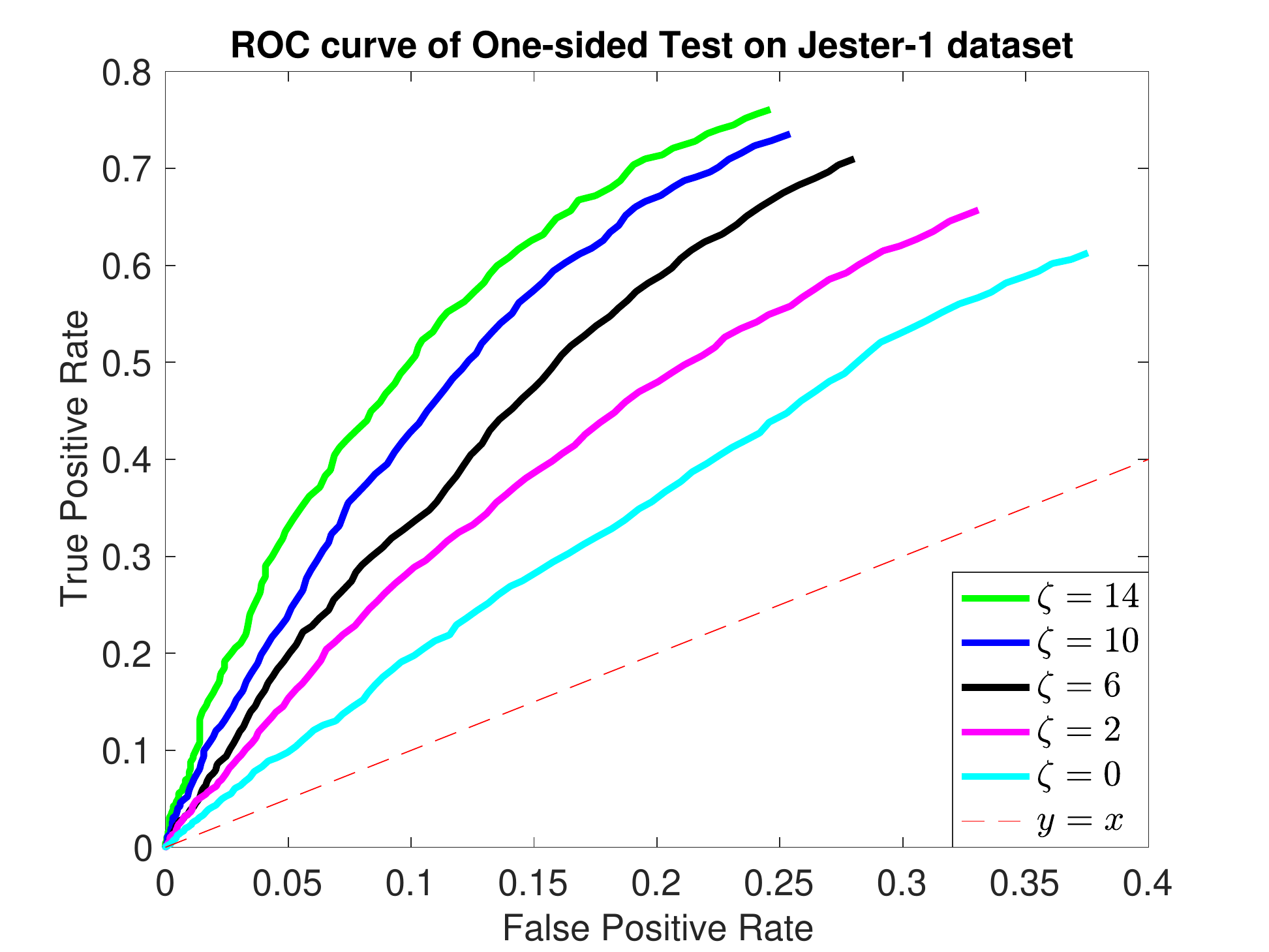}
\includegraphics[width=.32\textwidth]{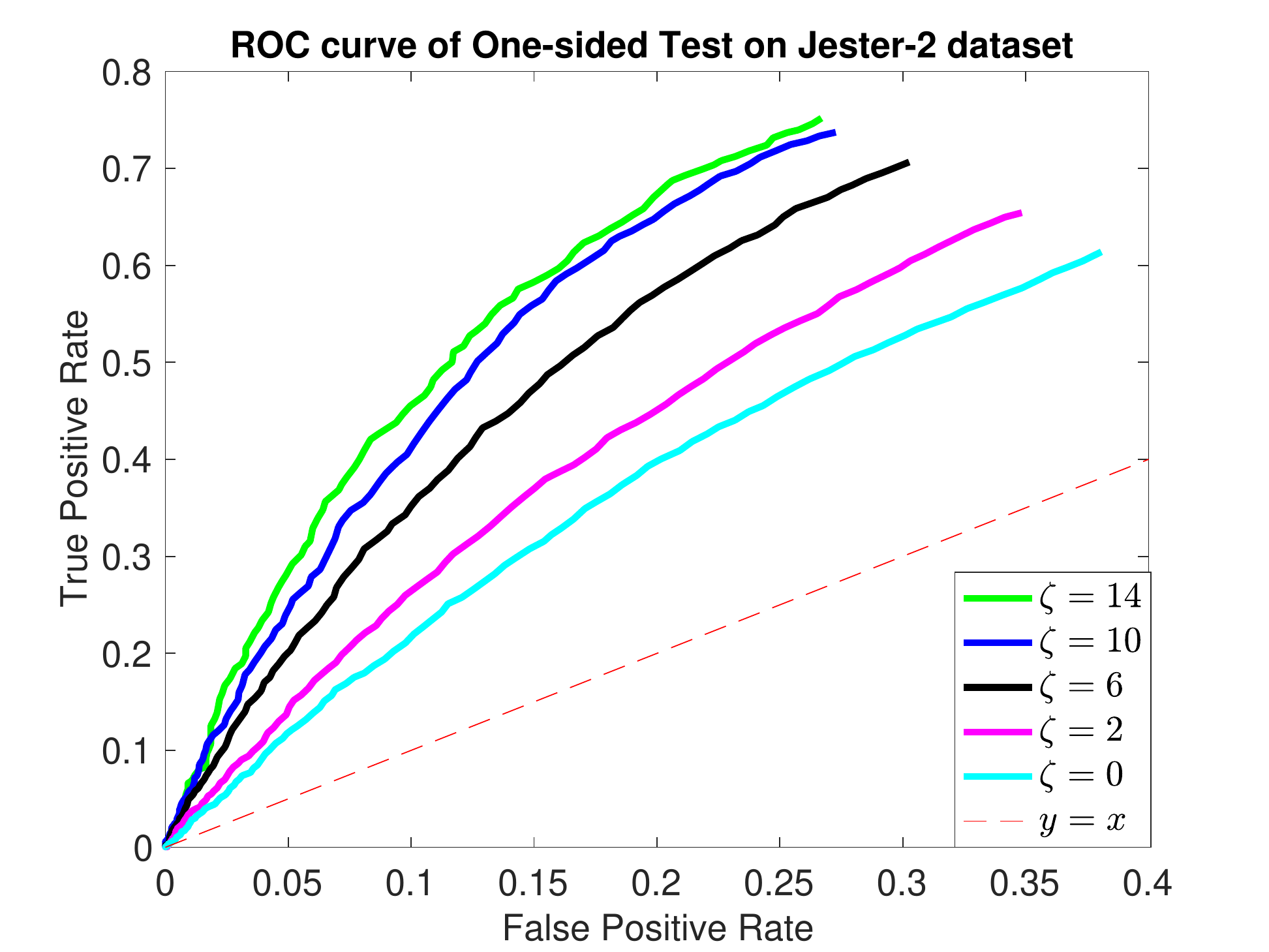}
\includegraphics[width=.32\textwidth]{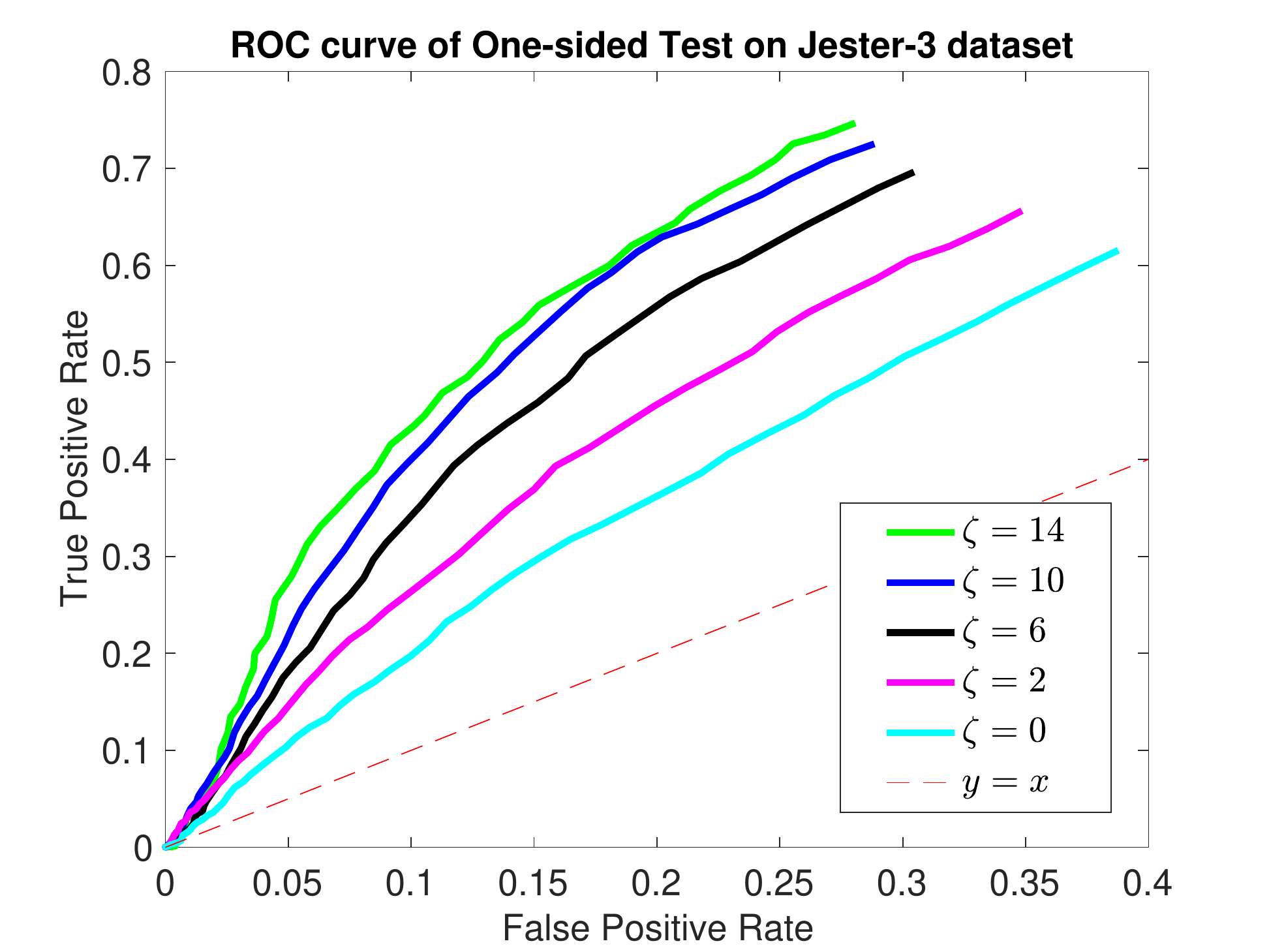}
\caption{ROC curves for one-sided tests $H_0: M(i,j_1)\leq M(i,j_2)\  {\rm v.s.}\  H_1: M(i,j_1)>M(i,j_2)$ on Jester datasets. The testing data are sampled such that $|Y(i,j_1)-Y(i,j_2)|\geq \zeta$. The estimated noise level $\hat\sigma_{\xi}=4.5160$ on Jester-1, $\hat\sigma_{\xi}=4.4843$ on Jester-2, and $\hat\sigma_{\xi}=5.1152$ on Jester-3. The rightmost point of each ROC curve corresponds to the significance level $\theta=0.5$ so that $z_{\theta}=0$.}
\label{fig:jester-ROC-oneside}
\end{figure}

For illustration, we consider the task of recommending jokes or movies to a particular users. Because of the lack of ground truth, we resort to resampling. For the Jester dataset, we randomly sample $\sim2000$ users, and for each user $2$ ratings that at least $\zeta\in\{0,2,6,10,14\}$ apart. We removed these ratings from the training and used the proposed procedure to infer, for each user ($i$), between these two jokes ($j_1$ or $j_2$) with ratings which one should be recommended. This amounts to the following one-sided tests:
$$
H_0: M(i,j_1)\leq M(i,j_2)\quad {\rm v.s.}\quad H_1: M(i,j_1)>M(i,j_2).
$$
We ran the proposed procedure on the training data and evaluate the test statistic $\hat{z}$ for each user from the testing set. In particular, we fixed the rank $r=2$ corresponding to the smallest estimate $\hat\sigma_{\xi}$. Note that we do not know the true value of $M(i,j)$ and only observe $Y(i,j)=M(i,j)+\xi(i,j)$. We therefore use ${\mathbb I}(Y(i,j_1)>Y(i,j_2))$ as a proxy to differentiate between $H_0$ and $H_1$. Assuming that the $\xi$ has a distribution symmetric about 0, then ${\mathbb I}(Y(i,j_1)>Y(i,j_2))$ is more likely to take value $0$ under $H_0$, and $1$ under $H_1$. We shall evaluate the performance of our procedure based on its discriminant power in predicting ${\mathbb I}(Y(i,j_1)>Y(i,j_2))$. In particular, we record the ROC curve of $\hat{z}$ for all users from the testing set. The results, averaged over $10$ simulation runs for each value of $\zeta$, are reported in Figure~\ref{fig:jester-ROC-oneside}. Clearly, we can observe an increase in predictive power as $\zeta$ increases suggesting $\hat{z}$ as a reasonable statistic for testing $H_0$ against $H_1$. 

%

\begin{figure}
\centering
\includegraphics[width=.32\textwidth]{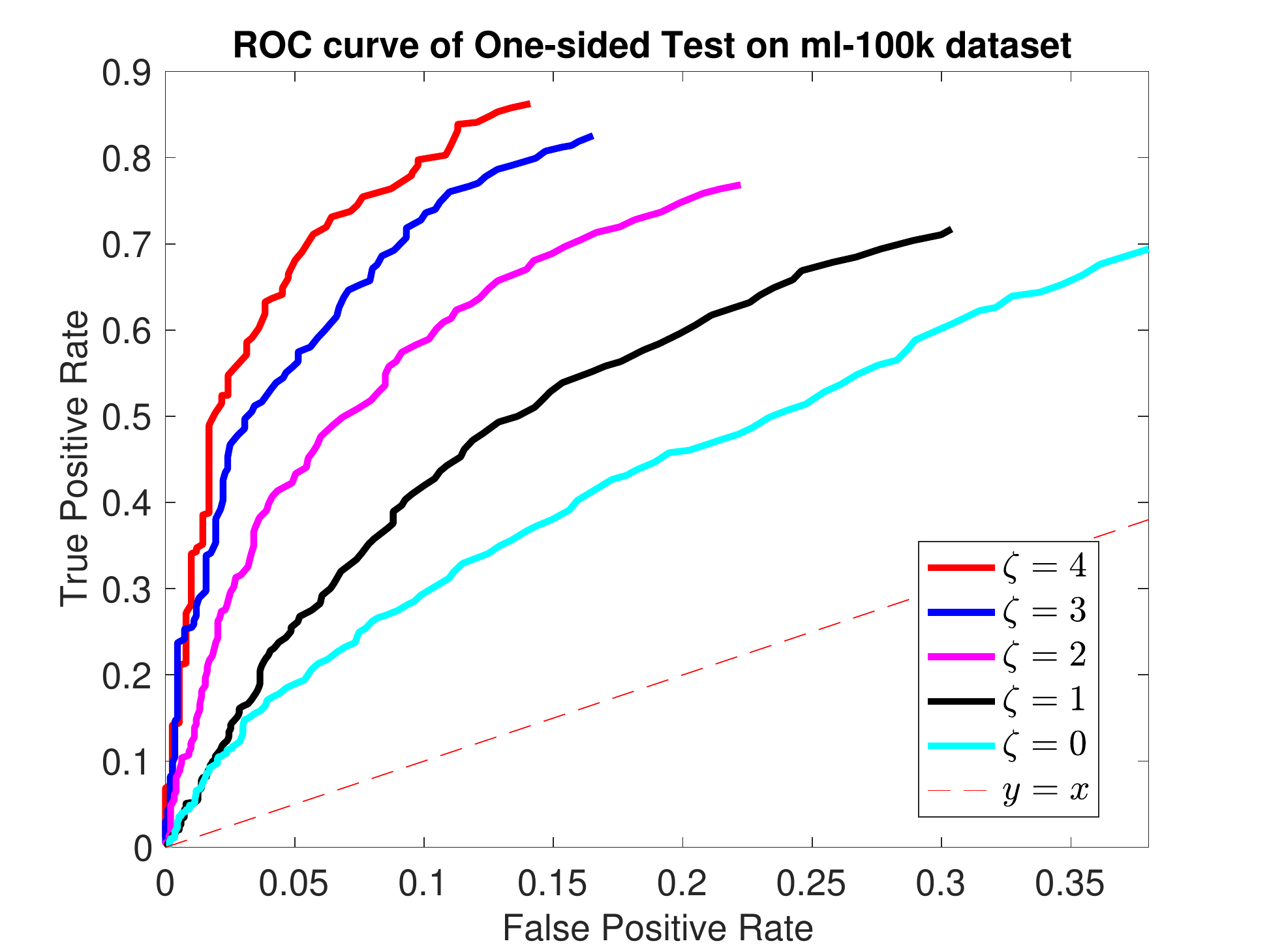}
\includegraphics[width=.32\textwidth]{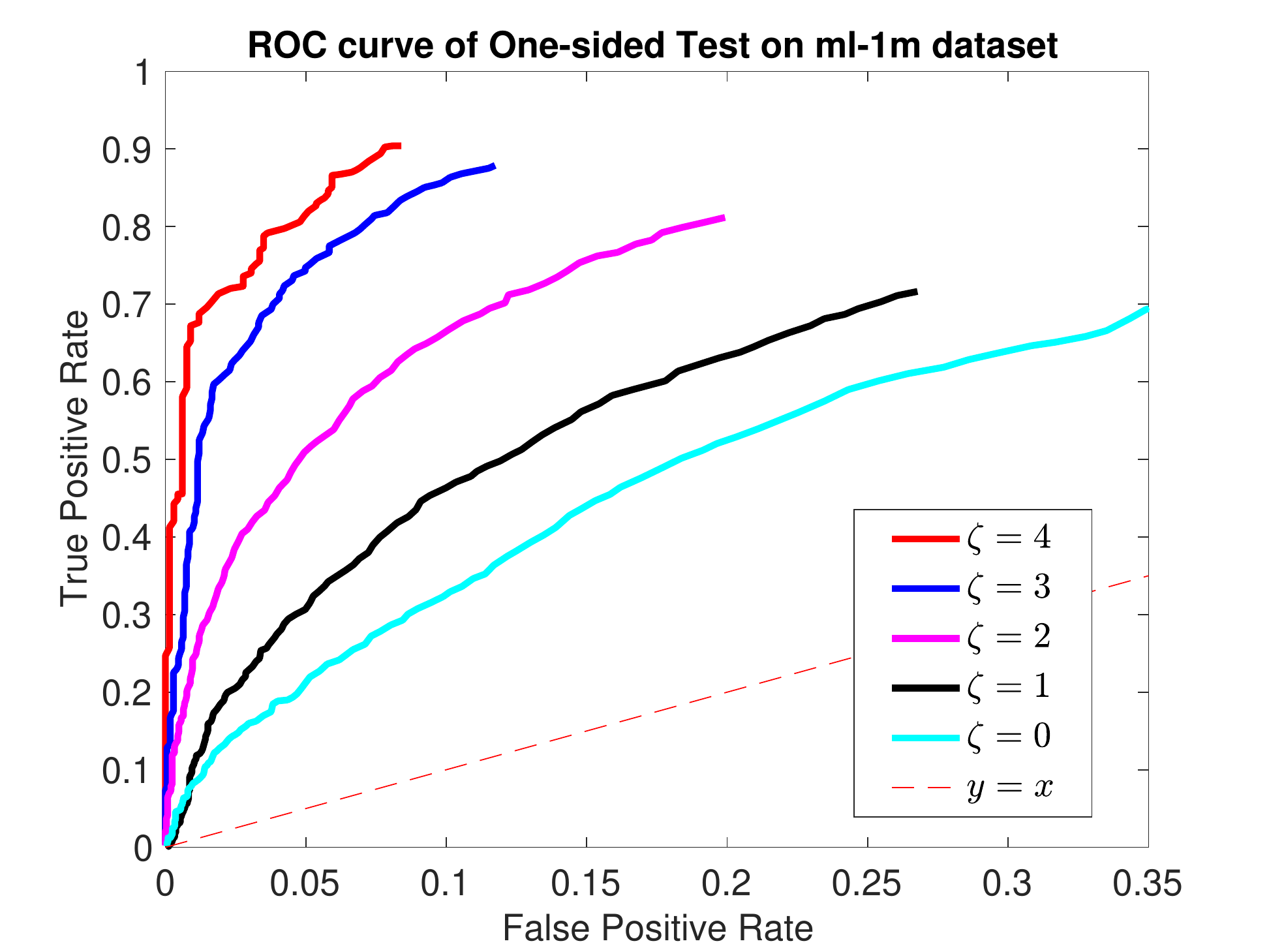}
\includegraphics[width=.32\textwidth]{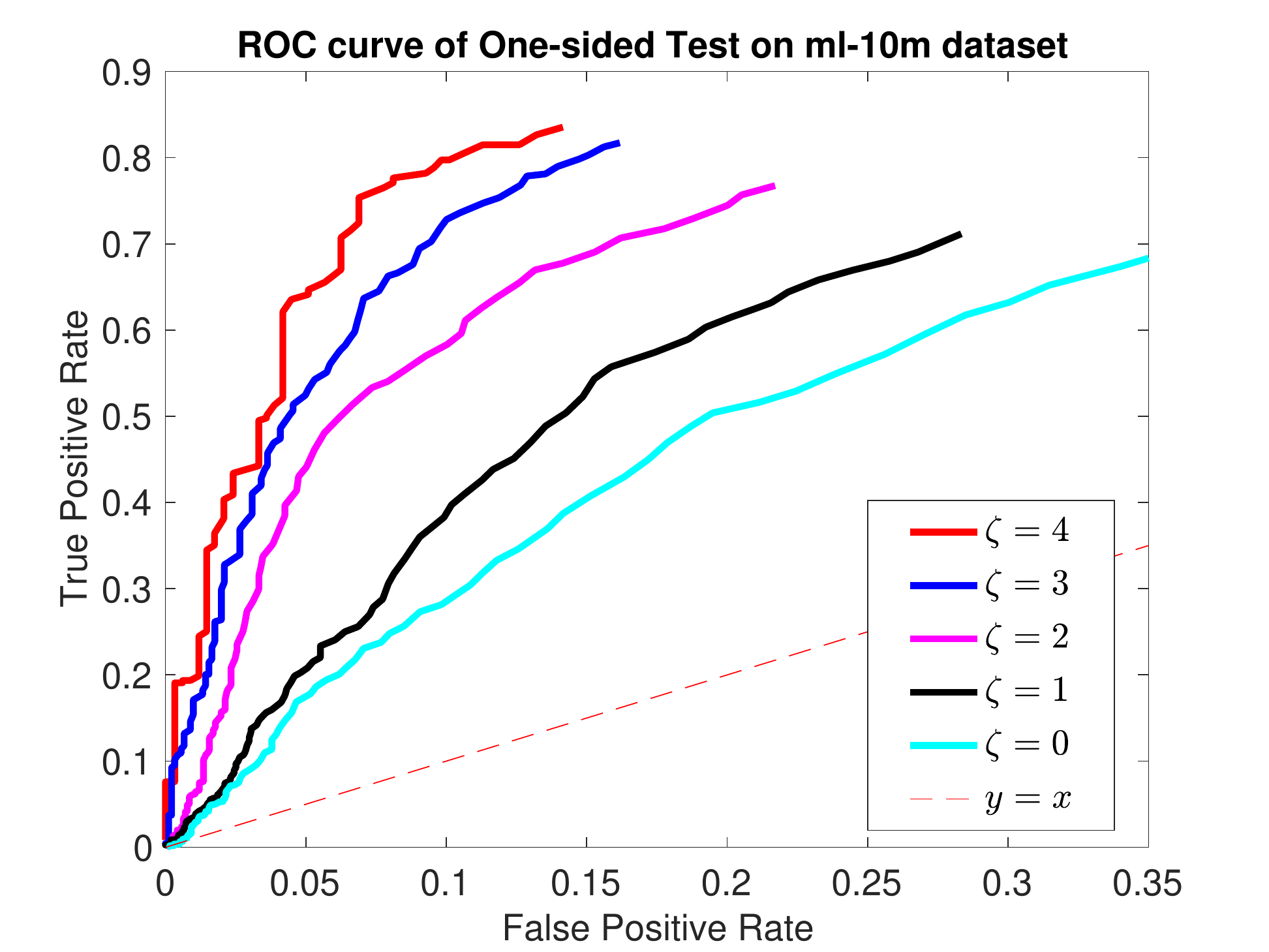}
\caption{ROC curves for one-sided tests $H_0: M(i,j_1)\leq M(i,j_2)\  {\rm v.s.}\  H_1: M(i,j_1)>M(i,j_2)$ on MovieLens datasets. The testing data are sampled such that $|Y(i,j_1)-Y(i,j_2)|\geq \zeta$. The estimated noise level $\hat\sigma_{\xi}=0.9973$ on ml-100k, $\hat\sigma_{\xi}=0.8936$ on ml-1m, and $\hat\sigma_{\xi}=0.9151$ on ml-10m. The rightmost point of each ROC curve corresponds to the significance level $\theta=0.5$ so that $z_{\theta}=0$.}
\label{fig:ml-ROC}
\end{figure}

We ran a similar experiment on the MovieLens datasets. In each simulation run, we randomly sampled $\sim800$ users and $2$ ratings each as the test data.  These ratings are sampled such that $|Y(i,j_1)-Y(i,j_2)|\geq \zeta$ for $\zeta=0,1,2,3,4$. The false positive rates and true positive rates of our proposed procedure were again recorded. The ROC curves, averaged again over $10$ runs for each value of $\zeta$, are shown in Figure~\ref{fig:ml-ROC}. This again indicates a reasonable performance of the proposed testing procedure. 
Empirically, we observe a better de-biasing approach on these datasets which is $\what M_1^{\supunbiased}=\what M_1^{\supinit}+\sum_{i\in\mathfrak{D}_2}(Y_i-\langle \what M_1^{\supinit},X_i\rangle)X_i$. The rationale is to partially replace $\what M_1^{\supinit}$'s entries with the observed training ratings. This improvement might be due to the severe heterogeneity in the numbers of observed ratings from distinct users, or due to the unknown noise distributions.

\section{Proofs}\label{sec:proofs}
Throughout the proof, we write $\gamma_n$ in short for $\gamma_{n,d_1,d_2}$. 
\subsection{De-localized perturbation of singular vectors}
Essential to our proofs is the precise characterization of the empirical singular spaces. To this end, we shall first develop bounds for the estimation error of $\what U_1,\what U_2, \what V_1$ and $\what V_2$. Recall that the matrix $(2,{\rm max})$-norm is defined as $\|A\|_{\submaxx}=\max_{j\in[d_1]}\|e_{j}^{\tran}A\|$. This can be naturally extended to a distance on Grassmannians 
$$
d_{\submaxx}(U_1,U_2):=\|U_1U_1^{\tran}-U_2U_2^{\tran}\|_{\submaxx},
$$
for $U_1,U_2\in \OO^{d\times r}$. The main goal of this subsection is to establish the following result:

\begin{theorem}\label{lem:hatUhatV_deloc}
Under Assumptions~\ref{assump:init_entry}-\ref{assump:noise}, there exists an absolute constant $C>0$ such that if $n\geq C\mumax^2rd_1\log d_1$, then with probability at least $1-5d_1^{-2}\log d_1$, 
$$
d_{\submaxx}(\what U_i, U)\leq C_2\mumax\frac{(1+\gamma_n)\sigma_{\xi}}{\lambda_r}\cdot\sqrt{\frac{rd_2d_1\log d_1}{n}}
$$
and
$$
d_{\submaxx}(\what V_i, V) \leq C_2\mumax\frac{(1+\gamma_n)\sigma_{\xi}}{\lambda_r}\cdot\sqrt{\frac{rd_1^2\log d_1}{n}}
$$
for $i=1,2$ and some absolute constant $C_2>0$.
\end{theorem}

Immediately following Theorem \ref{lem:hatUhatV_deloc} and Assumption~\ref{assump:noise}, we know that 
\begin{align*}
\big\|e_{j}^{\tran}\what U_1\big\|=\big\|e_j^{\tran}&(\what U_1\what U_1^{\tran}-UU^{\tran})\what U_1\big\|+\big\|e_j^{\tran}UU^{\tran}\what U_1\big\|\\
\leq&C_2\frac{\sigma_{\xi}}{\lambda_r}\sqrt{\frac{d_1^2d_2\log d_1}{n}}\cdot \mumax\sqrt{\frac{r}{d_1}}+\|e_j^{\tran}U\|\leq 2\mumax\sqrt{\frac{r}{d_1}}.
\end{align*}
Then, we conclude that
$$
\|\what U_i\|_{\submaxx}\leq 2\mumax\sqrt{\frac{r}{d_1}}\quad {\rm and}\quad \|\what V_i\|_{\submaxx}\leq 2\mumax\sqrt{\frac{r}{d_2}},\quad \forall i=1,2,
$$
an observation that we shall repeatedly use in the following subsections.

\subsubsection{Preliminary bounds}

Denote $\what\Delta_1=M-\what M_1^{\supinit} $ and $\what\Delta_2=M-\what M_2^{\supinit}$. We then write
\begin{equation}\label{eq:tildeM1}
\what{M}_1^{\supunbiased} =  M+\underbrace{\frac{d_1d_2}{n_0}\sum_{i=n_0+1}^n \xi_i X_i}_{\what Z_1^{(1)}}+\underbrace{\Big(\frac{d_1d_2}{n_0}\sum_{i=n_0+1}^{n}\langle\what\Delta_1,X_i \rangle X_i-\what\Delta_1\Big)}_{\what Z_2^{(1)}}
\end{equation}
and
\begin{equation}\label{eq:tildeM2}
\what{M}_2^{\supunbiased} = M+\underbrace{\frac{d_1d_2}{n_0}\sum_{i=1}^{n_0} \xi_i X_i}_{\what Z_1^{(2)}}+\underbrace{\Big(\frac{d_1d_2}{n_0}\sum_{i=1}^{n_0}\langle\what\Delta_2,X_i \rangle X_i-\what\Delta_2\Big)}_{\what Z_2^{(2)}},
\end{equation}
where $\what\Delta_1$ is independent with $\{(X_i,\xi_i)\}_{i=n_0+1}^n$, and $\what\Delta_2$ is independent with $\{(X_i,\xi_i)\}_{i=1}^{n_0}$.  Denote $\what{Z}^{(i)}=\what{Z}_1^{(i)}+\what{Z}_2^{(i)}$ and then $\what{M}_i^{\supunbiased}=M+\what{Z}^{(i)}$ for $i=1,2$. Clearly, $\EE\what{Z}^{(i)}=\EE\what{Z}_1^{(i)}+\EE\what{Z}_2^{(i)}=0$.

Observe that eq. (\ref{eq:tildeM1}, \ref{eq:tildeM2}) admit explicit representation formulas for $\what{M}_1^{\supunbiased}, \what{M}_2^{\supunbiased}$. Meanwhile, because $\|\what\Delta_1\|_{\submax}, \|\what\Delta_2\|_{\submax}=o_P(\sigma_{\xi})$, the strength of $\what{Z}^{(1)}_2$ and $\what{Z}^{(2)}_2$ are dominated by that of $\what{Z}^{(1)}_1$ and $\what{Z}^{(2)}_1$, respectively. Observe that the perturbation by $\what{Z}^{(1)}_1$ is analogous (or close) to a random perturbation with i.id. entry-wise noise.
Put it differently, the debiasing treatment by (\ref{eq:tildeM1},\ref{eq:tildeM2}) is essentially to re-randomize $\what M_1^{\supinit}$ and $\what M_2^{\supinit}$. It plays the key role in characterizing the distributions of $\what U_1, \what U_2$ and $\what V_1, \what V_2$.

We begin with several preliminary properties of $\{\what U_i\}_{i=1}^2$ and $\{\what V_i\}_{i=1}^2$. Recall that $\what U_1$ and $\what V_1$ are top-$r$ left and right singular vectors of $\what{M}_1^{\supunbiased}=M+\what{Z}_1^{(1)}+\what{Z}_2^{(1)}$. 
The following bounds for $\what{Z}^{(i)}_j$s are useful for our derivation.
\begin{lemma}\label{lem:Zbound}
There exist absolute constants $C_1, C_2>0$ such that if $n\geq C_1d_1\log d_1$, with probability at least $1-2d_1^{-2}$, the following bounds hold for $i=1,2$
$$
\big\|\what{Z}_1^{(i)} \big\|\leq C_2\sigma_{\xi}\sqrt{\frac{d_1^2d_2\log d_1}{n}}\quad{\rm and}\quad \big\|\what{Z}_2^{(i)} \big\|\leq C_2\|\what\Delta_i\|_{\submax}\sqrt{\frac{d_1^2d_2\log d_1}{n}},
$$
where the probability of the second inequality is conditioned on $\what\Delta_i$.
\end{lemma}

We shall defer the proof of Lemma \ref{lem:Zbound} to the Appendix. These bounds can be readily used to derive bounds for the empirical singular vectors under Frobenius-norm distance and operator-norm distance. Recall that for $U_1,U_2\in \OO^{d\times r}$, the Frobenius-norm distance and operator-norm distance are defined by 
$$
d_{\subf}(U_1,U_2)=\|U_1U_1^{\tran}-U_2U_2^{\tran}\|_{\subf}\quad{\rm and}\quad d_{\subo}(U_1,U_2)=\|U_1U_1^{\tran}-U_2U_2^{\tran}\|.
$$
It is well known that
$$
\min_{O\in\OO^{r\times r}}\|U_1-U_2O\|_{\subf}\leq \sqrt{2}d_{\subf}(U_1,U_2)\leq \sqrt{2}\cdot \min_{O\in\OO^{r\times r}}\|U_1-U_2O\|_{\subf}
$$
and
$$
\min_{O\in\OO^{r\times r}}\|U_1-U_2O\|\leq \sqrt{2}d_{\subo}(U_1,U_2)\leq \sqrt{2}\cdot \min_{O\in\OO^{r\times r}}\|U_1-U_2O\|.
$$
See, e.g., \cite{edelman1998geometry}.

\begin{lemma}\label{lemma:hatUhatV}
Under Assumptions~\ref{assump:init_entry}-\ref{assump:noise}, there exists an absolute constant $C>0$ such that if $n\geq Cd_1\log d_1$, then with probability at least $1-2d_1^{-2}$, 
$$
\max\{d_{\subf}(\what U_i, U), d_{\subf}(\what V_i,V)\}\leq C_2 \frac{(1+\gamma_n)\sigma_{\xi}}{\lambda_r}\cdot \sqrt{\frac{rd_1^2d_2\log{d_1}}{n}}\quad \textrm{ for } i=1,2,
$$
and
$$
\max\{d_{\subo}(\what U_i, U), d_{\subo}(\what V_i,V)\}\leq C_2 \frac{(1+\gamma_n)\sigma_{\xi}}{\lambda_r}\cdot \sqrt{\frac{d_1^2d_2\log d_1}{n}}\quad \textrm{ for } i=1,2,
$$
where $C_2>0$ is an absolute constant and $\gamma_n$ is defined by Assumption~\ref{assump:init_entry}. 
\end{lemma}

\subsection{Proof of Theorem~\ref{thm:normal_mt}}
We are now in position to prove Theorem~\ref{thm:normal_mt}. Recall that 
\begin{align*}
\what m_T-m_T=&\frac{1}{2}\big<\what U_1\what U_1^{\tran}\what{Z}^{(1)}\what V_1\what V_1^{\tran},T\big>+\frac{1}{2}\big<\what U_2\what U_2^{\tran}\what{Z}^{(2)}\what V_2\what V_2^{\tran},T\big>\\
&+\frac{1}{2}\big<\what U_1\what U_1^{\tran}M\what V_1\what V_1^{\tran}-M,T\big>+\frac{1}{2}\big<\what U_2\what U_2^{\tran}M\what V_2\what V_2^{\tran}-M,T\big>.
\end{align*}
Our strategy is to show that $\big\{\big|\big<\what U_i\what U_i^{\tran}\what{Z}^{(i)}\what V_i\what V_i^{\tran},T\big>\big|\big\}_{i=1}^2$ are negligible. Then, we prove the normal approximation of $\big\{\big<\what U_i\what U_i^{\tran}M\what V_i\what V_i^{\tran}-M,T\big>\big\}_{i=1}^2$.
We begin with the upper bounds of $\big\{\big|\big<\what U_i\what U_i^{\tran}\what{Z}^{(i)}\what V_i\what V_i^{\tran},T\big>\big|\big\}_{i=1}^2$. 
\begin{lemma}\label{lem:hatUZhatV}
Under Assumptions~\ref{assump:init_entry}-\ref{assump:noise},  and conditioned on the event in Theorem~\ref{lem:hatUhatV_deloc}, there exist  absolute constants $C_1,C_2>0$ such that with probability at least $1-2d_1^{-2}$,
\begin{align*}
\big|\big<\what U_i\what U_i^{\tran}\what{Z}^{(i)}\what V_i\what V_i^{\tran},T\big>\big|
\leq  C_1&\|T\|_{\ell_1}\mumax^2\sigma_{\xi}\sqrt{\frac{r\log d_1}{n}}\\
&+C_2\|T\|_{\ell_1}\mumax^2\frac{\sigma_{\xi}}{\lambda_r}\sqrt{\frac{rd_1^2d_2\log d_1}{n}}\cdot \sigma_{\xi}\sqrt{\frac{rd_1\log d_1}{n}}.
\end{align*}
\end{lemma}
We now prove the normal approximation of 
$$
\frac{1}{2}\big<\what U_1\what U_1^{\tran}M\what V_1\what V_1^{\tran}-M,T\big>+\frac{1}{2}\big<\what U_2\what U_2^{\tran}M\what V_2\what V_2^{\tran}-M,T\big>.
$$ 
Let $\Theta$ and $A$ be defined as in the proof of Theorem~\ref{lem:hatUhatV_deloc}. Moreover, we define 
$$
\what\Theta_1=\left(\begin{array}{cc}\what U_1&0 \\ 0&\what V_1\end{array}\right)\quad {\rm and}\quad \what\Theta_2=\left(\begin{array}{cc}\what U_2&0 \\ 0&\what V_2\end{array}\right).
$$
Then, we write
$$
\what\Theta_1\what\Theta_1^{\tran}A\what\Theta_1\what\Theta_1^{\tran}-\Theta\Theta^{\tran}A\Theta\Theta^{\tran}=\left(\begin{array}{cc}0&\what U_1\what U_1^{\tran}M\what V_1\what V_1^{\tran}-M\\ (\what U_1\what U_1^{\tran}M\what V_1\what V_1^{\tran}-M)^{\tran}&0\end{array}
\right)
$$
and
$$
\what\Theta_2\what\Theta_2^{\tran}A\what\Theta_2\what\Theta_2^{\tran}-\Theta\Theta^{\tran}A\Theta\Theta^{\tran}=\left(\begin{array}{cc}0&\what U_2\what U_2^{\tran}M\what V_2\what V_2^{\tran}-M\\ (\what U_2\what U_2^{\tran}M\what V_2\what V_2^{\tran}-M)^{\tran}&0\end{array}
\right).
$$
Denote
$$
\wtilde{T}=\left(\begin{array}{cc}0&T\\0&0\end{array}\right)\quad {\rm and}\quad \what{E}^{(i)}=\left(\begin{array}{cc}0&\what{Z}^{(i)}\\ \what{Z}^{(i)\tran}&0\end{array}\right)\quad \forall i=1,2.
$$
Therefore, we have 
\begin{align*}
\frac{1}{2}\big<\what U_1\what U_1^{\tran}M\what V_1\what V_1^{\tran}&-M,T\big>+\frac{1}{2}\big<\what U_2\what U_2^{\tran}M\what V_2\what V_2^{\tran}-M,T\big>\\
=&\frac{1}{2}\big<\what\Theta_1\what\Theta_1^{\tran}A\what\Theta_1\what\Theta_1^{\tran}-\Theta\Theta^{\tran}A\Theta\Theta^{\tran}, \wtilde{T}\big>+\frac{1}{2}\big<\what\Theta_2\what\Theta_2^{\tran}A\what\Theta_2\what\Theta_2^{\tran}-\Theta\Theta^{\tran}A\Theta\Theta^{\tran}, \wtilde{T}\big>.
\end{align*}
 By \eqref{lem:hatU-U}, we write
$$
\what\Theta_i\what\Theta_i^{\tran}-\Theta\Theta^{\tran}=\sum_{k=1}^{\infty}\calS_{A,k}(\what{E}^{(i)})\quad \forall i=1,2,
$$
and as a result, for $i=1,2$, 
\begin{align*}
\what\Theta_i\what\Theta_i^{\tran}A&\what\Theta\what\Theta^{\tran}-\Theta\Theta^{\tran}A\Theta\Theta^{\tran}\\
=&\big(\calS_{A,1}(\what{E}^{(i)})A\Theta\Theta^{\tran}+\Theta\Theta^{\tran}A\calS_{A,1}(\what{E}^{(i)})\big)+\sum_{k=2}^{\infty}\big(\calS_{A,k}(\what{E}^{(i)})A\Theta\Theta^{\tran}+\Theta\Theta^{\tran}A\calS_{A,k}(\what{E}^{(i)})\big)\\
&\hspace{2cm}+(\what\Theta_i\what\Theta_i^{\tran}-\Theta\Theta^{\tran})A(\what\Theta_i\what\Theta_i^{\tran}-\Theta\Theta^{\tran}).
\end{align*}
Then, we write
\begin{align*}
\frac{1}{2}\big<\what U_1\what U_1^{\tran}M\what V_1\what V_1^{\tran}&-M,T\big>+\frac{1}{2}\big<\what U_2\what U_2^{\tran}M\what V_2\what V_2^{\tran}-M,T\big>\\
=&\frac{1}{2}\sum_{i=1}^2\big<\big(\calS_{A,1}(\what{E}^{(i)})A\Theta\Theta^{\tran}+\Theta\Theta^{\tran}A\calS_{A,1}(\what{E}^{(i)})\big),\wtilde{T}\big>\\
+&\frac{1}{2}\sum_{i=1}^2\sum_{k=2}^{\infty}\big<\big(\calS_{A,k}(\what{E}^{(i)})A\Theta\Theta^{\tran}+\Theta\Theta^{\tran}A\calS_{A,k}(\what{E}^{(i)})\big),\wtilde{T}\big>\\
+&\frac{1}{2}\sum_{i=1}^2\big<(\what\Theta_i\what\Theta_i^{\tran}-\Theta\Theta^{\tran})A(\what\Theta_i\what\Theta_i^{\tran}-\Theta\Theta^{\tran}),\wtilde{T}\big>.
\end{align*}
By the definition of $\calS_{A,1}(\what{E}^{(i)})$, we write
\begin{align*}
\frac{1}{2}\sum_{i=1}^2\big<\calS_{A,1}&(\what{E}^{(i)})A\Theta\Theta^{\tran}+\Theta\Theta^{\tran}A\calS_{A,1}(\what{E}^{(i)}), \wtilde{T}\big>\\
=&\langle U_{\perp}U_{\perp}^{\tran}(\what{Z}^{(1)}/2+\what{Z}^{(2)}/2)VV^{\tran},T\rangle+\langle UU^{\tran}(\what{Z}^{(1)}/2+\what{Z}^{(2)}/2)V_{\perp}V_{\perp}^{\tran},T\rangle.
\end{align*}
We begin with the normal approximation of $\frac{1}{2}\sum_{i=1}^2\big<\calS_{A,1}(\what{E}^{(i)})A\Theta\Theta^{\tran}+\Theta\Theta^{\tran}A\calS_{A,1}(\what{E}^{(i)}), \wtilde{T}\big>$.
\begin{lemma}\label{lem:normal_SA1}
Under Assumptions~\ref{assump:init_entry}-\ref{assump:noise} and suppose that $n\geq C_1\mumax^2 rd_1\log d_1$, we have 
\begin{align*}
\sup_{x\in\RR}\bigg|\PP\bigg(&\frac{\frac{1}{2}\sum_{i=1}^2\big<\calS_{A,1}(\what{E}^{(i)})A\Theta\Theta^{\tran}+\Theta\Theta^{\tran}A\calS_{A,1}(\what{E}^{(i)}), \wtilde{T}\big>}{\sigma_{\xi}(\|V^{\tran}T^{\tran}\|_{\subf}^2+\|U^{\tran}T\|_{\subf}^2)^{1/2}\cdot\sqrt{d_1d_2/n}}\leq x\bigg)-\Phi(x)\bigg|\\
&\hspace{2cm}\leq C_2\frac{\mumax^4\|T\|_{\ell_1}^2}{\alpha_T^2\|T\|_{\subf}^2}\cdot\frac{r\sqrt{\log d_1}}{d_2}+\frac{3}{d_1^2}+C_3\gamma_n\sqrt{\log d_1}+C_4\mumax\sqrt{\frac{rd_1}{n}}.
\end{align*}
where $C_1, C_2, C_3,C_4>0$ are absolute constants and $\gamma_n$ is defined by Assumption~\ref{assump:init_entry}. 
\end{lemma}
Lemma~\ref{lem:sumSAk} and Lemma~\ref{lem:hatAhat_err} characterize sharp bounds for the remainder terms. 
\begin{lemma}\label{lem:sumSAk}
Under Assumptions~\ref{assump:init_entry}-\ref{assump:noise}, under the event of Theorem~\ref{lem:hatUhatV_deloc}, 
\begin{align*}
\big|\sum_{i=1}^2\sum_{k=2}^{\infty}\big<\big(\calS_{A,k}(\what{E}^{(i)})&A\Theta\Theta^{\tran}+\Theta\Theta^{\tran}A\calS_{A,k}(\what{E}^{(i)})\big),\wtilde{T}\big>\big|\\
&\leq C_2\|T\|_{\ell_1}\mumax^2\sigma_{\xi}\sqrt{\frac{rd_1\log d_1}{n}}\cdot \Big(\frac{\sigma_{\xi}}{\lambda_r}\cdot\sqrt{\frac{rd_1^2d_2\log d_1}{n}}\Big),
\end{align*}
where $C_2>0$ is some absolute constant. 
\end{lemma}

\begin{lemma}\label{lem:hatAhat_err}
Under Assumptions~\ref{assump:init_entry}-\ref{assump:noise}, on the event of Theorem~\ref{lem:hatUhatV_deloc}, 
\begin{align*}
\sum_{i=1}^2\big|\big<(\what\Theta_i\what\Theta_i^{\tran}-\Theta&\Theta^{\tran})A(\what\Theta_i\what\Theta_i^{\tran}-\Theta\Theta^{\tran}),\wtilde{T}\big>\big|\\
\leq&C_2\kappa_0\mumax^2\|T\|_{\ell_1}\sigma_{\xi}\sqrt{\frac{rd_1\log d_1}{n}}\cdot \frac{\sigma_{\xi}}{\lambda_r}\sqrt{\frac{rd_1^2d_2\log d_1}{n}}
\end{align*}
for some absolute constant $C_2>0$ and $\kappa(M)\leq\kappa_0$ denotes $M$'s condition number.  
\end{lemma}
We write 
\begin{align*}
&\frac{\what m_T-m_T}{\sigma_{\xi}(\|V^{\tran}T^{\tran}\|_{\subf}^2+\|U^{\tran}T\|_{\subf}^2)^{1/2}\cdot\sqrt{d_1d_2/n}}\\
&\hspace{2cm}=\frac{\sum_{i=1}^2\langle\what U_i\what U_i^{\tran}\what{Z}^{(i)}\what V_i\what V_i^{\tran},T \rangle/2}{\sigma_{\xi}(\|V^{\tran}T^{\tran}\|_{\subf}^2+\|U^{\tran}T\|_{\subf}^2)^{1/2}\cdot\sqrt{d_1d_2/n}}\\
&\hspace{2cm}+\frac{\sum_{i=1}^2\big<\calS_{A,1}(\what{E}^{(i)})A\Theta\Theta^{\tran}+\Theta\Theta^{\tran}A\calS_{A,1}(\what{E}^{(i)}), \wtilde{T}\big>/2}{\sigma_{\xi}(\|V^{\tran}T^{\tran}\|_{\subf}^2+\|U^{\tran}T\|_{\subf}^2)^{1/2}\cdot\sqrt{d_1d_2/n}}\\
&\hspace{2cm}+\frac{\sum_{i=1}^2\sum_{k=2}^{\infty}\big<\big(\calS_{A,k}(\what{E}^{(i)})A\Theta\Theta^{\tran}+\Theta\Theta^{\tran}A\calS_{A,k}(\what{E}^{(i)})\big),\wtilde{T}\big>/2}{\sigma_{\xi}(\|V^{\tran}T^{\tran}\|_{\subf}^2+\|U^{\tran}T\|_{\subf}^2)^{1/2}\cdot\sqrt{d_1d_2/n}}\\
&\hspace{2cm}+\frac{\sum_{i=1}^2\big<(\what\Theta_i\what\Theta_i^{\tran}-\Theta\Theta^{\tran})A(\what\Theta_i\what\Theta_i^{\tran}-\Theta\Theta^{\tran}),\wtilde{T}\big>/2}{\sigma_{\xi}(\|V^{\tran}T^{\tran}\|_{\subf}^2+\|U^{\tran}T\|_{\subf}^2)^{1/2}\cdot\sqrt{d_1d_2/n}}.
\end{align*}
Under Assumption~\ref{assump:T}, it holds that $\|V^{\tran}T^{\tran}\|_{\subf}^2+\|U^{\tran}T\|_{\subf}^2\geq \alpha_T^2 \|T\|_{\subf}^2 r/d_1$. As a result, 
$$
\sigma_{\xi}(\|V^{\tran}T^{\tran}\|_{\subf}^2+\|U^{\tran}T\|_{\subf}^2)^{1/2}\cdot\sqrt{d_1d_2/n}\geq \alpha_T\|T\|_{\subf}\sigma_{\xi}\sqrt{\frac{rd_2}{n}}.
$$
Together with Lemma~\ref{lem:hatUZhatV},\ref{lem:sumSAk} and Lemma~\ref{lem:hatAhat_err}, we get, with probability at least $1-6d_1^{-2}\log d_1$, that
\begin{align*}
&\frac{\big|(\what m_T-m_T)-\sum_{i=1}^2\big<\calS_{A,1}(\what{E}^{(i)})A\Theta\Theta^{\tran}+\Theta\Theta^{\tran}A\calS_{A,1}(\what{E}^{(i)}), \wtilde{T}\big>/2\big|}{\sigma_{\xi}(\|V^{\tran}T^{\tran}\|_{\subf}^2+\|U^{\tran}T\|_{\subf}^2)^{1/2}\cdot\sqrt{d_1d_2/n}}\\
&\hspace{1cm}\leq C_1\frac{\mumax^2\|T\|_{\ell_1}}{\alpha_T\|T\|_{\subf}}\cdot\sqrt{\frac{\log d_1}{d_2}}+C_2\kappa_0\frac{\mumax^2\|T\|_{\ell_1}}{\alpha_T\|T\|_{\subf}}\cdot \frac{\sigma_{\xi}}{\lambda_r}\sqrt{\frac{\alpha_drd_1^2d_2\log^2d_1}{n}}
\end{align*}
for some absolute constants $C_1,C_2>0$. By the normal approximation of $\sum_{i=1}^2\big<\calS_{A,1}(\what{E}^{(i)})A\Theta\Theta^{\tran}+\Theta\Theta^{\tran}A\calS_{A,1}(\what{E}^{(i)}), \wtilde{T}\big>/2$ in Lemma~\ref{lem:normal_SA1} and the Lipschitz property of $\Phi(x)$, we  get 
\begin{align*}
\sup_{x\in\RR}\Big|\PP&\Big(\frac{\what m_T-m_T}{\sigma_{\xi}(\|V^{\tran}T^{\tran}\|_{\subf}^2+\|U^{\tran}T\|_{\subf}^2)^{1/2}\cdot\sqrt{d_1d_2/n}}\leq x\Big)-\Phi(x)\Big|\\
\leq&C_1\frac{\mumax^2\|T\|_{\ell_1}}{\alpha_T\|T\|_{\subf}}\sqrt{\frac{\log d_1}{d_2}}+C_2\kappa_0\frac{\mumax^2\|T\|_{\ell_1}}{\alpha_T\|T\|_{\subf}}\cdot \frac{\sigma_{\xi}}{\lambda_r}\sqrt{\frac{\alpha_d rd_1^2d_2\log^2d_1}{n}}\\
+&C_4\frac{\mumax^4\|T\|_{\ell_1}^2}{\alpha_T^2\|T\|_{\subf}^2}\cdot\frac{r\sqrt{\log d_1}}{d_2}+\frac{6\log d_1}{d_1^2}+C_5\gamma_n\sqrt{\log d_1}+C_6\mumax\sqrt{\frac{rd_1}{n}},
\end{align*}
which concludes the proof of Theorem~\ref{thm:normal_mt}.

\subsection{Proof of Theorem~\ref{thm:conf_int}}
It suffices to prove the normal approximation of 
$$
\frac{\what m_T-m_T}{\what{\sigma}_{\xi}\what s_T\cdot\sqrt{d_1d_2/n}}
$$
with data-driven estimators $\what \sigma_{\xi}$ and $\what s_T$. Write 
\begin{align*}
&\frac{\what m_T-m_T}{\what{\sigma}_{\xi}\what s_T\cdot\sqrt{d_1d_2/n}}\\
&\hspace{2cm}=\frac{\what m_T-m_T}{\sigma_{\xi}(\| V^{\tran}T^{\tran}\|_{\subf}^2+\| U^{\tran}T\|_{\subf}^2)^{1/2}\cdot\sqrt{d_1d_2/n}}\\
&\hspace{2.5cm}+\frac{\what m_T-m_T}{\what{\sigma}_{\xi}\what s_T\cdot\sqrt{d_1d_2/n}}\cdot\Big(1-\frac{\what\sigma_{\xi}}{\sigma_{\xi}}\Big)\\
&\hspace{2.5cm}+\frac{\what m_T-m_T}{\sigma_{\xi}(\| V^{\tran}T^{\tran}\|_{\subf}^2+\| U^{\tran}T\|_{\subf}^2)^{1/2}\cdot\sqrt{d_1d_2/n}}\cdot\Big(\frac{(\| V^{\tran}T^{\tran}\|_{\subf}^2+\| U^{\tran}T\|_{\subf}^2)^{1/2}}{\what s_T}-1\Big).
\end{align*}
Recall that 
\begin{align*}
\what\sigma_{\xi}^2=\frac{1}{2n_0}\sum_{i=n_0+1}^{n}&\big(Y_i-\langle \what{M}_1^{\supinit},X_i\rangle\big)^2+\frac{1}{2n_0}\sum_{i=1}^{n_0}(Y_i-\langle \what M_2^{\supinit},X_i\rangle)^2\\
&=\frac{1}{n}\sum_{i=1}^n\xi_i^2+\frac{1}{2n_0}\sum_{i=n_0+1}^{n}\langle\what\Delta_1, X_i\rangle^2+\frac{1}{2n_0}\sum_{i=1}^{n_0}\langle\what \Delta_2,X_i \rangle^2\\
&\hspace{2cm}+\frac{1}{n_0}\sum_{i=n_0+1}^n\xi_i\langle\what{\Delta}_1,X_i \rangle+\frac{1}{n_0}\sum_{i=1}^{n_0}\xi_i\langle\what{\Delta}_2,X_i \rangle.
\end{align*}
Note that $\{(X_i,\xi_i)\}_{i=n_0+1}^n$ are independent with $\what\Delta_1$. By Bernstein inequality and under Assumption~\ref{assump:init_entry}, it is easy to show that, with probability at least $1-2d_1^{-2}$, 
\begin{align*}
\big|\what\sigma_{\xi}^2-\sigma_{\xi}^2 \big|\leq \frac{2(\|\what\Delta_1\|_{\subf}^2+\|\what\Delta_2\|_{\subf}^2)}{d_1d_2}+\frac{C_1\sigma_{\xi}^2\log d_1}{\sqrt{n}}\leq \frac{C_1\sigma_{\xi}^2\log d_1}{\sqrt{n}}+C_2\gamma_n^2\cdot \sigma_{\xi}^2. 
\end{align*}
Then, if $C_2\gamma_n^2\leq 1/3$ so that $|\what\sigma_{\xi}^2-\sigma_{\xi}^2|\leq \sigma_{\xi}^2/2$, we get 
\begin{align*}
\Big|1-\frac{\what\sigma_{\xi}}{\sigma_{\xi}}\Big|\leq \Big|1-\frac{\what\sigma_{\xi}^2}{\sigma_{\xi}^2}\Big|\leq \frac{C_1\log d_1}{\sqrt{n}}+C_2\gamma_n^2.
\end{align*}
We now bound $\big| \|TV\|_{\subf}^2-\|T\what V_1\|_{\subf}^2\big|$. Observe that $V$ and $\what V_1$ both have orthonormal columns. Then, 
\begin{align*}
\big| \|TV\|_{\subf}^2-\|T&\what V_1\|_{\subf}^2\big|=\big| \|TVV^{\tran}\|_{\subf}^2-\|T\what V_1\what V_1^{\tran}\|_{\subf}^2\big|\\
&\leq \big\|T(VV^{\tran}-\what V_1\what V_1^{\tran})\big\|_{\subf}^2+2\big|\big<T(VV^{\tran}-\what V_1\what V_1^{\tran}), TV V^{\tran}\big>\big|.
\end{align*}
Clearly, 
\begin{align*}
\big\|T(VV^{\tran}-\what V_1\what V_1^{\tran})\big\|_{\subf}^2
\leq \Big(\sum_{(j_1,j_2)\in{\rm supp}(T)} |T_{j_1,j_2}| \big\|e_{j_2}^{\tran}(VV^{\tran}-\what V_1\what V_1^{\tran})\big\|\Big)^2\\
\leq \|T\|_{\ell_1}^2\cdot \|VV^{\tran}-\what V_1\what V_1^{\tran}\|_{\submaxx}^2\leq C_1\mumax^2 \frac{\|T\|_{\ell_1}^2}{d_2}\cdot \Big(\frac{\sigma_{\xi}^2}{\lambda_r^2}\Big)\frac{rd_1^2d_2\log d_1}{n}.
\end{align*}
Similarly, 
\begin{align*}
\big|\big<T(VV^{\tran}-\what V_1\what V_1^{\tran}), TV V^{\tran}\big>\big|\leq \|TV\|_{\subf}\|T(VV^{\tran}-\what V_1\what V_1^{\tran})V\|_{\subf}\\
\leq \|TV\|_{\subf}\|T\|_{\ell_1}\|(VV^{\tran}-\what V_1\what V_1^{\tran})V\|_{\submaxx}\\
\leq \|TV\|_{\subf}\|T\|_{\ell_1}\mumax\cdot \frac{\sigma_{\xi}}{\lambda_r}\sqrt{\frac{rd_1^2\log d_1}{n}}.
\end{align*}
Therefore,
\begin{align*}
\big| \|TV\|_{\subf}^2-\|T&\what V_1\|_{\subf}^2\big|=\big| \|TVV^{\tran}\|_{\subf}^2-\|T\what V_1\what V_1^{\tran}\|_{\subf}^2\big|\\
\leq&C_1\mumax^2 \frac{\|T\|_{\ell_1}^2}{d_2}\cdot \Big(\frac{\sigma_{\xi}^2}{\lambda_r^2}\Big)\frac{rd_1^2d_2\log d_1}{n}+C_2 \|TV\|_{\subf}\|T\|_{\ell_1}\mumax\cdot \frac{\sigma_{\xi}}{\lambda_r}\sqrt{\frac{rd_1^2\log d_1}{n}}.
\end{align*}
Similar bounds can be shown for
$
\big| \|U^{\tran}T\|_{\subf}^2-\|\what U_1^{\tran}T\|_{\subf}^2\big|.
$
The same bounds also hold for $\|\what U_2^{\tran}T\|_{\subf}^2$ and $\|T\what V_2\|_{\subf }^2$. Under the event of Theorem~\ref{lem:hatUhatV_deloc}, 
\begin{align*}
\big|\what s_T^2 -&\big(\| V^{\tran}T^{\tran}\|_{\subf}^2+\| U^{\tran}T\|_{\subf}^2\big)\big|\\
\leq& C_1\mumax^2\frac{\|T\|_{\ell_1}^2}{d_2}\cdot \Big(\frac{\sigma_{\xi}^2}{\lambda_r^2}\Big)\frac{rd_1^2d_2\log d_1}{n}+C_2 \|TV\|_{\subf}\|T\|_{\ell_1}\mumax\cdot \frac{\sigma_{\xi}}{\lambda_r}\sqrt{\frac{rd_1^2\log d_1}{n}},
\end{align*}
and as a result
\begin{align*}
\bigg|\frac{\what s_T^2}{\| V^{\tran}T^{\tran}\|_{\subf}^2+\| U^{\tran}T\|_{\subf}^2}-1\bigg|\leq C_1\mumax\frac{\|T\|_{\ell_1}}{\|T\|_{\subf}\alpha_T}\cdot \frac{\sigma_{\xi}}{\lambda_r}\sqrt{\frac{\alpha_dd_1^2d_2\log d_1}{n}},
\end{align*}
where we used the fact $\| V^{\tran}T^{\tran}\|_{\subf}^2+\| U^{\tran}T\|_{\subf}^2\geq \alpha_T^2\|T\|_{\subf}^2(r/d_1)$ and also the fact 
$$
\| V^{\tran}T^{\tran}\|_{\subf}^2+\| U^{\tran}T\|_{\subf}^2\geq \max\{\|TV\|_{\subf}, \|U^{\tran}T\|_{\subf}\}\cdot \alpha_T\|T\|_{\subf}\sqrt{r/d_1},
$$
due to Assumption~\ref{assump:T}.  It also implies, under condition (\ref{eq:asymp_cond}), that 
$$
\what s_T^2\geq \big(\| V^{\tran}T^{\tran}\|_{\subf}^2+\| U^{\tran}T\|_{\subf}^2\big)/2\geq \alpha_T^2\|T\|_{\subf}^2\cdot \frac{r}{2d_1}.
$$
Then, 
\begin{align*}
\bigg| \frac{(\| V^{\tran}T^{\tran}\|_{\subf}^2+\| U^{\tran}T\|_{\subf}^2)^{1/2}}{\what s_T}-1\bigg|\leq& \bigg| \frac{\| V^{\tran}T^{\tran}\|_{\subf}^2+\| U^{\tran}T\|_{\subf}^2}{\what s_T^2}-1\bigg|\\
\leq&C_2\mumax\frac{\|T\|_{\ell_1}}{\|T\|_{\subf}\alpha_T}\cdot \frac{\sigma_{\xi}}{\lambda_r}\sqrt{\frac{\alpha_dd_1^2d_2\log d_1}{n}}.
\end{align*}
By the normal approximation in Theorem~\ref{thm:normal_mt}, there is an event $\calE_2$ with 
\begin{align*}
\PP(\calE_2)\geq &1-C_1\frac{\|T\|_{\ell_1}\mumax^2}{\|T\|_{\subf}\alpha_T}\sqrt{\frac{\log d_1}{d_2}}-C_2\kappa_0\frac{\|T\|_{\ell_1}\mumax^2}{\|T\|_{\subf}\alpha_T}\cdot \frac{\sigma_{\xi}}{\lambda_r}\sqrt{\frac{\alpha_drd_1^2d_2\log^2d_1}{n}}\\
-&C_3\frac{\|T\|_{\ell_1}^2\mumax^2}{\|T\|_{\subf}^2\alpha_T^2}\cdot\frac{r\sqrt{\log d_1}}{d_1}-\frac{6\log d_1}{d_1^2}-C_4\gamma_n\sqrt{\log d_1}-C_5\mumax\sqrt{\frac{rd_1}{n}},
\end{align*}
so that on event $\calE_2$, 
$$
\frac{\what m_T-m_T}{\what{\sigma}_{\xi}\what s_T\cdot\sqrt{d_1d_2/n}}\leq C_6\sqrt{\log d_1}
$$
and
$$
\frac{\what m_T-m_T}{\sigma_{\xi}(\| V^{\tran}T^{\tran}\|_{\subf}^2+\| U^{\tran}T\|_{\subf}^2)^{1/2}\cdot\sqrt{d_1d_2/n}}\leq C_6\sqrt{\log d_1}.
$$
Therefore, under event $\calE_2$, with probability at least $1-2d_1^{-2}$, 
\begin{align}
\bigg|\frac{\what m_T-m_T}{\what{\sigma}_{\xi}\what s_T\cdot\sqrt{d_1d_2/n}}\cdot\Big(1-\frac{\what\sigma_{\xi}}{\sigma_{\xi}}\Big)\bigg|
\leq \frac{C_1\log^{3/2}d_1}{\sqrt{n}}+C_2\gamma_n^2\sqrt{\log d_1}\label{eq:ci_err1}
\end{align}
and
\begin{align}
&\bigg|\frac{\what m_T-m_T}{\sigma_{\xi}(\| V^{\tran}T^{\tran}\|_{\subf}^2+\| U^{\tran}T\|_{\subf}^2)^{1/2}\cdot\sqrt{d_1d_2/n}}\cdot\Big(\frac{(\| V^{\tran}T^{\tran}\|_{\subf}^2+\| U^{\tran}T\|_{\subf}^2)^{1/2}}{\what s_T}-1\Big)\bigg|\nonumber\\
&\hspace{2cm}\leq C_2\mumax\frac{\|T\|_{\ell_1}}{\|T\|_{\subf}\alpha_T}\cdot \frac{\sigma_{\xi}}{\lambda_r}\sqrt{\frac{\alpha_drd_1^2d_2\log d_1^2}{n}}\label{eq:ci_err2}.
\end{align}
As a result, if 
$$
\lim_{d_1,d_2\to\infty}\max\bigg\{\frac{\mumax^2\|T\|_{\ell_1}}{\alpha_T\|T\|_{\subf}}\sqrt{\frac{r\log d_1}{d_2}}, \frac{\kappa_0\mumax^2\|T\|_{\ell_1}}{\alpha_T\|T\|_{\subf}}\cdot \frac{\sigma_{\xi}}{\lambda_r}\sqrt{\frac{\alpha_drd_1^2d_2\log d_1^2}{n}}, \gamma_n\sqrt{\log d_1}\bigg\}= 0,
$$
then
$$
\frac{\what m_T-m_T}{\what{\sigma}_{\xi}\what s_T\cdot\sqrt{d_1d_2/n}}\overset{{\rm d}}{\longrightarrow} \calN(0,1)
$$
as $d_1,d_2\to\infty$.

\subsection{Proof of Theorem~\ref{thm:gd}}

We begin with the accuracy of $\what G^{(t)}$. 
By the definition of $\what G^{(t)}$, we have 
\begin{equation}\label{eq:hatG_t}
\frac{d_1d_2}{N_0}\sum_{j\in\mathfrak{D}_{2t}}\big(\langle \what U^{(t)}\what G^{(t)}(\what V^{(t)})^{\tran}, X_j\rangle -Y_j\big)\what U^{(t)\tran}X_j\what V^{(t)}=0.
\end{equation}
To this end, let $\what O_U^{(t)}$ and $\what O_V^{(t)}$ be any orthogonal matrices so that 
\begin{equation}
\max\big(\|\what U^{(t)}-U\what O_U^{(t)}\|, \|\what V^{(t)}-V\what O_V^{(t)}\|\big)\leq \frac{1}{C_1\mumax \kappa_0^2\sqrt{r}}
\end{equation}
for some large enough constant $C_1>0$. 

\begin{lemma}\label{lem:hatG}
Suppose that $\|\what U^{(t)}\|_{\submaxx}\leq 2\mumax\sqrt{r/d_1}$ and $\|\what V^{(t)}\|_{\submaxx}\leq 2\mumax\sqrt{r/d_2}$ and if $n\geq C_2\mumax^4 r^3(r^2+\log d_1)\log d_1$, then with probability at least $1-3d_1^{-2}$,
\begin{align*}
\big\|\what G^{(t)}-&\what O_U^{(t)\tran}\Lambda\what O_V^{(t)} \big\|\\
\leq& C_3\sigma_{\xi}\sqrt{\frac{rd_1d_2\log d_1}{N_0}}+2\|\Lambda\|\cdot \big(\|\what U^{(t)}-U\what O_U^{(t)}\|^2+\|\what V^{(t)}-V\what O_V^{(t)}\|^2\big)\\
+&C_4\|\Lambda\|\Big(\sqrt{\frac{r}{d_2}}\|\what U^{(t)}-U\what O_U^{(t)}\|_{\submaxx}+\sqrt{\frac{r}{d_1}}\|\what V^{(t)}-V\what O_V^{(t)}\|_{\submaxx}\big)\cdot\mumax\sqrt{\frac{rd_1d_2\log d_1}{N_0}}
\end{align*}
for some absolute constants $C_2,C_3,C_4>0$.
\end{lemma}
Let $\what G^{(t)}=\what L_G^{(t)}\what\Lambda^{(t)}\what R_G^{(t)\tran}$ denote $\what G^{(t)}$'s SVD where $\what L_G^{(t)}$, $\what R_G^{(t)}$ are both $r\times r$ orthogonal matrices and $\what \Lambda^{(t)}$ is a diagonal matrix. 
Recall the gradient descent step of Algorithm~\ref{algo:GD},
\begin{align*}
\what U^{(t+0.5)}=\what U^{(t)}\what L_G^{(t)}-\eta \cdot \frac{d_1d_2}{N_0}\sum_{j\in\mathfrak{D}_{2t+1}}\big<\what U^{(t)}\what G^{(t)}\what V^{(t)\tran}-U\Lambda V^{\tran}, X_j\big>X_j\what V^{(t)}\what R^{(t)}_G(\what \Lambda^{(t)})^{-1}\\
-\eta\cdot \frac{d_1d_2}{N_0}\sum_{j\in\mathfrak{D}_{2t+1}}\xi_j X_j \what V^{(t)}\what R^{(t)}_G(\what \Lambda^{(t)})^{-1}.
\end{align*}
Observe that $(\what U^{(t)}, \what V^{(t)},\what G^{(t)})$ are independent with $\mathfrak{D}_{2t+1}$. Then, we write
\begin{align*}
\what U^{(t+0.5)}=\what U^{(t)}\what L_G^{(t)}-\eta\cdot (\what U^{(t)}\what G^{(t)}\what V^{(t)\tran}-U\Lambda V^{\tran})\what V^{(t)}\what R_G^{(t)}(\what \Lambda^{(t)})^{-1}+\what{E}_{V}^{(t)}+\what{E}_{\xi,V}^{(t)},
\end{align*}
where 
\begin{align*}
\what{E}_{V}^{(t)}=\eta\cdot (\what U^{(t)}&\what G^{(t)}\what V^{(t)\tran}-U\Lambda V^{\tran})\what V^{(t)}\what R_G^{(t)}(\what \Lambda^{(t)})^{-1}\\
-&\eta\cdot \frac{d_1d_2}{N_0}\sum_{j\in\mathfrak{D}_{2t+1}}\big<\what U^{(t)}\what G^{(t)}\what V^{(t)\tran}-U\Lambda V^{\tran}, X_j\big>X_j\what V^{(t)}\what R_G^{(t)}(\what \Lambda^{(t)})^{-1}
\end{align*}
and 
$$
\what{E}_{\xi,V}^{(t)}=-\eta\cdot \frac{d_1d_2}{N_0}\sum_{j\in\mathfrak{D}_{2t+1}}\xi_j X_j \what V^{(t)}\what R_G^{(t)}(\what \Lambda^{(t)})^{-1}.
$$
Note that 
\begin{align*}
 (\what U^{(t)}\what G^{(t)}\what V^{(t)\tran}-U\Lambda& V^{\tran})\what V^{(t)}\what R_G^{(t)}(\what \Lambda^{(t)})^{-1}= \what U^{(t)}\big(\what G^{(t)}-\what O_U^{(t)\tran}\Lambda\what O_V^{(t)}\big)\what R^{(t)}_G(\what \Lambda^{(t)})^{-1}\\
 +&\big(\what U^{(t)}\what L_G^{(t)}-U\what O_U^{(t)}\what L_G^{(t)}\big)\what L_G^{(t)\tran}\what O_U^{(t)\tran}\Lambda \what O_V^{(t)}\what R_G^{(t)}(\what \Lambda^{(t)})^{-1}\\
 &+U\Lambda(\what V^{(t)}\what O_V^{(t)\tran}-V)^{\tran}\what V^{(t)}\what R_G^{(t)}(\what \Lambda^{(t)})^{-1}.
\end{align*}
Therefore,
\begin{align}
\what U^{(t+0.5)}=&U\what O_U^{(t)}\what L_G^{(t)}+(\what U^{(t)}\what L_G^{(t)}-U\what O^{(t)}_U\what L_G^{(t)})\big(I-\eta\cdot \what L_G^{(t)\tran}\what O_U^{(t)\tran}\Lambda \what O_V^{(t)}\what R_G^{(t)}(\what \Lambda^{(t)})^{-1}\big)\nonumber\\
-&\eta\cdot \what U^{(t)}\big(\what G^{(t)}-\what O_U^{(t)\tran}\Lambda\what O_V^{(t)}\big)\what R^{(t)}_G(\what \Lambda^{(t)})^{-1}\nonumber\\
&-\eta\cdot U\Lambda(\what V^{(t)}\what O_V^{(t)\tran}-V)^{\tran}\what V^{(t)}\what R_G^{(t)}(\what \Lambda^{(t)})^{-1}\nonumber\\
&\hspace{2cm}+\what E_V^{(t)}+\what E_{\xi,V}^{(t)}. \label{eq:hatUt0.5}
\end{align}

\begin{lemma}\label{lem:gd_infty}
Under Assumptions~\ref{assump:incoh} and \ref{assump:noise} and the assumptions of Lemma~\ref{lem:hatG}, suppose that $\|\what U^{(t)}\|_{\submaxx}\leq 2\mumax\sqrt{r/d_1}$, $\|\what V^{(t)}\|_{\submaxx}\leq 2\mumax\sqrt{r/d_2}$,
$$
\max\Big\{\|\what U^{(t)}-U\what O_U^{(t)}\|, \|\what V^{(t)}-V\what O_V^{(t)}\| \Big\}\leq 1/(C_1\mumax\kappa_0\sqrt{r\alpha_d}),
$$
and $n\geq C_2\alpha_d\kappa_0^2\mumax^4r^2d_1(r+\log d_1)\log d_1$ for some large enough constant $C_1,C_2>0$, if $\eta\in [0.25, 0.75]$, then the following bound holds with probability at least $1-2d_1^{-2}$,
\begin{align*}
\big\|\what U^{(t+0.5)}-U\what O_U^{(t)}\what L_G^{(t)} \big\|_{\submaxx}\leq \Big(1-&\frac{9\eta}{10}\Big)\|\what U^{(t)}-U\what O_U^{(t)}\|_{\submaxx}+C_3\eta\frac{\sigma_{\xi}}{\lambda_r}\sqrt{\frac{rd_1d_2\log d_1}{N_0}}\\
&+\frac{\eta}{8}\cdot \big(\|\what U^{(t)}-U\what O_U^{(t)}\|_{\submaxx}+\|\what V^{(t)}-V\what O_V^{(t)}\|_{\submaxx}\big);
\end{align*}
and with probability at least $1-2d_1^{-2}$,
\begin{align*}
\big\{|1-\lambda_r&(\what U^{(t+0.5)})|, |1-\lambda_1(\what U^{(t+0.5)})| \big\}
\leq C_3\eta\frac{\sigma_{\xi}}{\lambda_r}\cdot\sqrt{\frac{d_1^2d_2\log d_1}{N_0}}\\
+C_4&(\kappa_0\eta+\kappa_0^2\eta^2)\cdot\big(\|\what U^{(t)}-U\what O_U^{(t)}\|^2+\|\what V^{(t)}-V\what O_V^{(t)}\|^2\big)\\
+C_5&\eta\kappa_0\cdot\Big(\sqrt{\frac{r}{d_2}}\|\what U^{(t)}-U\what O_U^{(t)}\|_{\submaxx}+\sqrt{\frac{r}{d_1}}\|\what V^{(t)}-V\what O_V^{(t)}\|_{\submaxx}\Big)\cdot\mumax\sqrt{\frac{rd_1^2d_2}{N_0}},
\end{align*}
for some absolute constants $C_3,C_4,C_5>0$. 
\end{lemma}

By Lemma~\ref{lem:gd_infty}, we denote the SVD of $\what U^{(t+0.5)}$ by $\what U^{(t+0.5)}=\what U^{(t+1)}\what \Sigma_U^{(t+1)}\what K_U^{(t+1)\tran}$ where $\what\Sigma_U^{(t+1)}$ is diagonal and
\begin{align*}
&\|\what\Sigma_U^{(t+1)}-I\|\\
\leq& C_3\eta\frac{\sigma_{\xi}}{\lambda_r}\cdot\sqrt{\frac{d_1^2d_2\log d_1}{N_0}}+C_4(\kappa_0\eta+\kappa_0^2\eta^2)\cdot\big(\|\what U^{(t)}-U\what O_U^{(t)}\|^2+\|\what V^{(t)}-V\what O_V^{(t)}\|^2\big)\\
+C_5&\eta\kappa_0\cdot\Big(\sqrt{\frac{r}{d_2}}\|\what U^{(t)}-U\what O_U^{(t)}\|_{\submaxx}+\sqrt{\frac{r}{d_1}}\|\what V^{(t)}-V\what O_V^{(t)}\|_{\submaxx}\Big)\cdot\mumax\sqrt{\frac{rd_1^2d_2}{N_0}}.
\end{align*}
By $\what U^{(t+1)}\what \Sigma_U^{(t+1)}\what K_U^{(t+1)\tran}=U\what O_U^{(t)}\what L_G^{(t)}+\big(\what U^{(t+0.5)}-U\what O_U^{(t)}\what L_G^{(t)}\big)$, we write
\begin{align*}
\what U^{(t+1)}=U\what O_U^{(t)}\what L_G^{(t)}\what K_U^{(t+1)}(\what \Sigma_U^{(t+1)})^{-1}+\big(\what U^{(t+0.5)}-U\what O_U^{(t)}\what L_G^{(t)}\big)\what K_U^{(t+1)}(\what \Sigma_U^{(t+1)})^{-1}
\end{align*}
and obtain
\begin{align}
\what U^{(t+1)}-U\what O_U^{(t)}\what L_G^{(t)}\what K_U^{(t+1)}=U\what O_U^{(t)}\what L_G^{(t)}\what K_U^{(t+1)}\big((\what \Sigma_U^{(t+1)})^{-1}-I\big)\nonumber\\
+\big(\what U^{(t+0.5)}-U\what O_U^{(t)}\what L_G^{(t)}\big)\what K_U^{(t+1)}(\what \Sigma_U^{(t+1)})^{-1}.\label{eq:hatUt+1}
\end{align}
Note that $\what O_U^{(t)}\what L_G^{(t)}\what K_U^{(t+1)}$ is an $r\times r$ orthogonal matrix. The Assumptions of Lemma~\ref{lem:gd_infty} can guarantee $\lambda_r(\what \Sigma_U^{(t+1)})\geq 1-\eta/20$ so that $\|(\what \Sigma_U^{(t+1)})^{-1}\|\leq 1+\eta/10$. 

Therefore,
\begin{align*}
\big\|\what U^{(t+1)}-U&\what O_U^{(t)}\what L_G^{(t)}\what K_U^{(t+1)}\big\|_{\submaxx}
\leq& \|U\|_{\submaxx}\cdot \big\|(\what \Sigma_U^{(t+1)})^{-1}-I\big\|+(1+\eta/10)\|\what U^{(t+0.5)}-U\what O_U^{(t)}\what L_G^{(t)}\|_{\submaxx}.
\end{align*}
Then, by Lemma~\ref{lem:gd_infty}, 
\begin{align*}
\big\|&\what U^{(t+1)}-U\what O_U^{(t)}\what L_G^{(t)}\what K_U^{(t+1)}\big\|_{\submaxx}
\leq \mumax\sqrt{\frac{r}{d_1}} \big\|(\what \Sigma_U^{(t+1)})^{-1}-I\big\|+(1+\eta/10)\|\what U^{(t+0.5)}-U\what O_U^{(t)}\what L_G^{(t)}\|_{\submaxx}\\
\leq& C_3\eta\mumax\cdot \frac{\sigma_{\xi}}{\lambda_r}\sqrt{\frac{rd_1d_2\log d_1}{N_0}}
+C_4(\kappa_0\eta+\kappa_0^2\eta^2)\mumax\cdot \sqrt{\frac{r}{d_1}}\big(\|\what U^{(t)}-U\what O_U^{(t)}\|^2+\|\what V^{(t)}-V\what O_V^{(t)}\|^2\big)\\
&\hspace{0.5cm}+C_5 \eta\cdot\Big(\sqrt{\frac{r}{d_2}}\|\what U^{(t)}-U\what O_U^{(t)}\|_{\submaxx}+\sqrt{\frac{r}{d_1}}\|\what V^{(t)}-V\what O_V^{(t)}\|_{\submaxx}\Big)\cdot\mumax^2\kappa_0\sqrt{\frac{r^2d_1 d_2}{N_0}}\\
&\hspace{0.5cm}+\Big(1-\frac{4\eta}{5}\Big)\|\what U^{(t)}-U\what O_U^{(t)}\|_{\submaxx}+C_3\eta\frac{\sigma_{\xi}}{\lambda_r}\sqrt{\frac{rd_1d_2\log d_1}{N_0}}\\
&\hspace{1cm}+\frac{\eta}{7}\cdot \big(\|\what U^{(t)}-U\what O_U^{(t)}\|_{\submaxx}+\|\what V^{(t)}-V\what O_V^{(t)}\|_{\submaxx}\big)\\
&\hspace{0.5cm}\leq \Big(1-\frac{4\eta}{5}\Big)\|\what U^{(t)}-U\what O_U^{(t)}\|_{\submaxx}+\frac{\eta}{6}\cdot \big(\|\what U^{(t)}-U\what O_U^{(t)}\|_{\submaxx}+\|\what V^{(t)}-V\what O_V^{(t)}\|_{\submaxx}\big)\\
&\hspace{1cm}+C_3\eta\frac{\sigma_{\xi}}{\lambda_r}\sqrt{\frac{rd_1d_2\log d_1}{N_0}},
\end{align*}
where the last inequality holds as long as $n\geq C_5\alpha_d\kappa_0^2\mumax^4 d_1r^3 \log d_1$, $\|\what U^{(t)}-U\what O_U^{(t)}\|+\|\what V^{(t)}-V\what O_V^{(t)}\|\leq (C_6\kappa_0^2\mumax\sqrt{r\alpha_d})^{-1}$  for some large enough constants $C_5,C_6>0$. Then, 
\begin{align*}
\big\|\what U^{(t+1)}-U\what O_U^{(t)}\what L_G^{(t)}\what K_U^{(t+1)}\big\|_{\submaxx}\leq C_3\eta\frac{\sigma_{\xi}}{\lambda_r}\sqrt{\frac{rd_1d_2\log d_1}{N_0}}\\
+ \Big(1-\frac{4\eta}{5}\Big)\|\what U^{(t)}-U\what O_U^{(t)}\|_{\submaxx}+\frac{\eta}{6}\cdot \big(\|\what U^{(t)}-U\what O_U^{(t)}\|_{\submaxx}+\|\what V^{(t)}-V\what O_V^{(t)}\|_{\submaxx}\big).
\end{align*}
Similarly, the gradient descent step for $\what V^{(t)}$ reads
\begin{align*}
\what V^{(t+0.5)}=\what V^{(t)}\what R_G^{(t)}-\eta \cdot \frac{d_1d_2}{N_0}\sum_{j\in\mathfrak{D}_{2t+1}}\big<\what U^{(t)}\what G^{(t)}\what V^{(t)\tran}-U\Lambda V^{\tran}, X_j\big>X_j^{\tran}\what U^{(t)}\what L^{(t)}_G(\what \Lambda^{(t)})^{-1}\\
-\eta\cdot \frac{d_1d_2}{N_0}\sum_{j\in\mathfrak{D}_{2t+1}}\xi_j X_j^{\tran} \what U^{(t)}\what L^{(t)}_G(\what \Lambda^{(t)})^{-1}.
\end{align*}
Let $\what V^{(t+0.5)}=\what V^{(t+1)}\what\Sigma_V^{(t+1)} \what K_V^{(t+1)\tran}$ denote $\what V^{(t+0.5)}$'s SVD where $\what K_V^{(t+1)}$ is an orthogonal matrix. In the same fashion,  with probability at least $1-4d_1^{-2}$, 
\begin{align*}
\big\|\what V^{(t+1)}-V\what O_V^{(t)}\what R_G^{(t)}\what K_V^{(t+1)}\big\|_{\submaxx}\leq C_3\eta\frac{\sigma_{\xi}}{\lambda_r}\sqrt{\frac{rd_1d_2\log d_1}{N_0}}\\
+ \Big(1-\frac{4\eta}{5}\Big)\|\what V^{(t)}-V\what O_V^{(t)}\|_{\submaxx}+\frac{\eta}{6}\cdot \big(\|\what U^{(t)}-U\what O_U^{(t)}\|_{\submaxx}+\|\what V^{(t)}-V\what O_V^{(t)}\|_{\submaxx}\big).
\end{align*}
Then, we conclude with 
\begin{align}
\big\|\what U^{(t+1)}-U\what O_U^{(t)}\what L_G^{(t)}\what K_U^{(t+1)}\big\|_{\submaxx}+\big\|\what V^{(t+1)}-V\what O_V^{(t)}\what R_G^{(t)}\what K_V^{(t+1)}\big\|_{\submaxx}\nonumber\\
\leq C_3\eta\frac{\sigma_{\xi}}{\lambda_r}\sqrt{\frac{rd_1d_2\log d_1}{N_0}}+\Big(1-\frac{2\eta}{3}\Big)\cdot \big(\|\what U^{(t)}-U\what O_U^{(t)}\|_{\submaxx}+\|\what V^{(t)}-V\what O_V^{(t)}\|_{\submaxx}\big),\label{eq:hatUVt+1}
\end{align}
where both $\what O_U^{(t)}\what L_G^{(t)}\what K_U^{(t+1)}$ and $\what O_V^{(t)}\what R_G^{(t)}\what K_V^{(t+1)} $ are orthogonal matrices. 

The contraction property of the iterations is then proved after replacing $\what O_U^{(t)}$ and $\what O_V^{(t)}$ with the orthogonal matrices defined in Theorem~\ref{thm:gd}. It suffices to show that 
\begin{equation}\label{eq:UVt_opt}
\max\Big\{\|\what U^{(t)}-U\what O_U^{(t)}\|, \|\what V^{(t)}-V\what O_V^{(t)}\|\Big\}\leq \frac{1}{C_6\mumax \kappa_0^2\sqrt{r\alpha_d}},
\end{equation}
and $\|\what U^{(t)}\|_{\submaxx}\leq 2\mumax\sqrt{r/d_1}, \|\what V^{(t)}\|_{\submaxx}\leq 2\mumax\sqrt{r/d_2}$
for all $1\leq t\leq m$ and  some large constant $C_6>0$.

We first show $\|\what U^{(t)}\|_{\submaxx}\leq 2\mumax\sqrt{r/d_1}, \|\what V^{(t)}\|_{\submaxx}\leq 2\mumax\sqrt{r/d_2}$ for all $1\leq t\leq m$. By the contraction property (\ref{eq:hatUVt+1}), it suffices to show $\|\what U^{(1)}-U\what O_U^{(1)}\|_{\submaxx}\leq \mumax\sqrt{r/d_1}, \|\what V^{(1)}-V\what O_V^{(1)}\|_{\submaxx}\leq \mumax\sqrt{r/d_1}$ and $C_3(\sigma_{\xi}/\lambda_r)\sqrt{rd_1d_2\log d_1/N_0}\leq \mumax\sqrt{r/d_1}$ where the last inequality holds automatically under Assumption~\ref{assump:noise}. Similarly as the proof of Theorem~\ref{lem:hatUhatV_deloc}, with probability at least $1-5d_1^{-2}\log d_1$, 
\begin{align*}
d_{\submaxx}(\what U^{(1)},U)\leq& C_2\mumax\frac{\sigma_{\xi}+\|M\|_{\submax}}{\lambda_r}\cdot \sqrt{\frac{rd_2d_1\log d_1}{N_0}},\\
 d_{\submaxx}(\what V^{(1)},V)\leq& C_2\mumax\frac{\sigma_{\xi}+\|M\|_{\submax}}{\lambda_r}\cdot \sqrt{\frac{rd_1d_1\log d_1}{N_0}}.
\end{align*}
Since $\|M\|_{\submax}\leq \|\Lambda\| \|U\|_{\submaxx}\|V\|_{\submaxx}\leq \mumax^2\|\Lambda\|\sqrt{r^2/d_1d_2}$, it implies $\|\what U^{(1)}-U\what O_U^{(1)}\|_{\submaxx}\leq \mumax\sqrt{r/d_1}$ and  $\|\what V^{(1)}-V\what O_V^{(1)}\|_{\submaxx}\leq \mumax\sqrt{r/d_1}$ as long as 
$$
n\geq C_2\alpha_d\mumax^4\kappa_0^2 r^2d_1\log^2d_1\quad {\rm and}\quad C_2\frac{\sigma_{\xi}}{\lambda_r}\cdot \sqrt{\frac{\alpha_d d_1^2d_2\log^2d_1}{n}}\leq 1
$$
for some large enough constant $C_2>0$. 

We then show (\ref{eq:UVt_opt}) for all $t=1,\cdots,m$. By eq. (\ref{eq:hatUt0.5_3}), we write
\begin{align*}
\what U^{(t+1)}&\what \Sigma_U^{(t+1)}-U\what O_U^{(t)}\what L_G^{(t)}\what K_U^{(t+1)}\\
&=(\what U^{(t)}\what L_G^{(t)}-U\what O^{(t)}_U\what L_G^{(t)})\big(I-\eta\cdot \what L_G^{(t)\tran}\what O_U^{(t)\tran}\Lambda \what O_V^{(t)}\what R_G^{(t)}(\what \Lambda^{(t)})^{-1}\big)\what K_U^{(t+1)}\nonumber\\
-&\eta\cdot \what U^{(t)}\big(\what G^{(t)}-\what O_U^{(t)\tran}\Lambda\what O_V^{(t)}\big)\what R^{(t)}_G(\what \Lambda^{(t)})^{-1}\what K_U^{(t+1)}\\
-&\eta\cdot U\Lambda(\what V^{(t)}\what O_V^{(t)\tran}-V)^{\tran}\what V^{(t)}\what R_G^{(t)}(\what \Lambda^{(t)})^{-1}\what K_U^{(t+1)}+\what E_V^{(t)}\what K_U^{(t+1)}+\what E_{\xi,V}^{(t)}\what K_U^{(t+1)}.
\end{align*} 
Similar as the proof of Lemma~\ref{lem:gd_infty} and (\ref{eq:hatUt0.5_3}), we can write 
\begin{align*}
\big\| \what U^{(t+1)}\what \Sigma_U^{(t+1)}&-U\what O_U^{(t)}\what L_G^{(t)}\what K_U^{(t+1)}\big\|
\leq (1-0.9\eta)\|\what U^{(t)}-U\what O_U^{(t)}\|+2\eta\frac{\|\what G^{(t)}-\what O_U^{(t)\tran}\Lambda\what O_V^{(t)}\|}{\lambda_r}\\
&\hspace{2cm}+2\eta\kappa_0\|\what V^{(t)}-V\what O_V^{(t)}\|^2+\big\|\what E_V^{(t)}\what K_U^{(t+1)}+\what E_{\xi,V}^{(t)}\what K_U^{(t+1)} \big\|,
\end{align*}
and as a result
\begin{align*}
\big\| \what U^{(t+1)}&-U\what O_U^{(t)}\what L_G^{(t)}\what K_U^{(t+1)}\big\|\\
\leq& (1-0.9\eta)\|\what U^{(t)}-U\what O_U^{(t)}\|+2\eta\frac{\|\what G^{(t)}-\what O_U^{(t)\tran}\Lambda\what O_V^{(t)}\|}{\lambda_r}+2\eta\kappa_0\|\what V^{(t)}-V\what O_V^{(t)}\|^2\\
&\hspace{2cm}+\big\|\what E_V^{(t)}\what K_U^{(t+1)}+\what E_{\xi,V}^{(t)}\what K_U^{(t+1)} \big\|+\|\what \Sigma_U^{(t+1)}-I\|.
\end{align*}
Then, by Lemma~\ref{lem:hatG}-\ref{lem:gd_infty} and the upper bound of $\|\what E_V^{(t)}+\what E^{(t)}_{\xi,V}\|$ in the proof of Lemma~\ref{lem:gd_infty},
\begin{align*}
\big\|&\what U^{(t+1)}-U\what O_U^{(t)}\what L_G^{(t)}\what K_U^{(t+1)}\big\|\leq  (1-0.8\eta)\|\what U^{(t)}-U\what O_U^{(t)}\|+ C_3\eta\frac{\sigma_{\xi}}{\lambda_r}\cdot\sqrt{\frac{d_1^2d_2\log d_1}{N_0}}\\
&+C_4(\kappa_0\eta+\kappa_0^2\eta^2)\cdot\big(\|\what U^{(t)}-U\what O_U^{(t)}\|^2+\|\what V^{(t)}-V\what O_V^{(t)}\|^2\big)\\
&+C_5\eta\kappa_0\cdot\Big(\sqrt{\frac{r}{d_2}}\|\what U^{(t)}-U\what O_U^{(t)}\|_{\submaxx}+\sqrt{\frac{r}{d_1}}\|\what V^{(t)}-V\what O_V^{(t)}\|_{\submaxx}\Big)\cdot\mumax\sqrt{\frac{rd_1^2d_2}{N_0}}.
\end{align*}
Similarly, we can get the bound for $\|\what V^{(t+1)}-V\what O_V^{(t)}\what R_G^{(t)}\what K_V^{(t+1)} \|$ and as a result
\begin{align*}
\big\|&\what U^{(t+1)}-U\what O_U^{(t)}\what L_G^{(t)}\what K_U^{(t+1)}\big\|+\|\what V^{(t+1)}-V\what O_V^{(t)}\what R_G^{(t)}\what K_V^{(t+1)} \|\\
&\leq  (1-0.8\eta)\big(\|\what U^{(t)}-U\what O_U^{(t)}\|+\|\what V^{(t)}-U\what O_V^{(t)}\|\big)+ C_3\eta\frac{\sigma_{\xi}}{\lambda_r}\cdot\sqrt{\frac{d_1^2d_2\log d_1}{N_0}}\\
&+C_4(\kappa_0\eta+\kappa_0^2\eta^2)\cdot\big(\|\what U^{(t)}-U\what O_U^{(t)}\|^2+\|\what V^{(t)}-V\what O_V^{(t)}\|^2\big)\\
&+C_5\eta\kappa_0\cdot\Big(\sqrt{\frac{r}{d_2}}\|\what U^{(t)}-U\what O_U^{(t)}\|_{\submaxx}+\sqrt{\frac{r}{d_1}}\|\what V^{(t)}-V\what O_V^{(t)}\|_{\submaxx}\Big)\cdot\mumax\sqrt{\frac{rd_1^2d_2}{N_0}}.
\end{align*}
By the previous proof, with probability at least $1-2d_1^{-2}$,
$$
\|\what U^{(t)}-U\what O_U^{(t)}\|_{\submaxx}\leq C_1\mumax\sqrt{\frac{r}{d_1}}\quad
{\rm and}
\quad 
\|\what V^{(t)}-V\what O_V^{(t)}\|_{\submaxx}\leq C_1\mumax\sqrt{\frac{r}{d_2}}.
$$
If $\|\what U^{(t)}-U\what O_U^{(t)}\|+\|\what V^{(t)}-V\what O_V^{(t)}\|\leq 1/(3C_4\mumax\kappa_0^2\sqrt{r\alpha_d})$, we get 
\begin{align*}
\big\|&\what U^{(t+1)}-U\what O_U^{(t)}\what L_G^{(t)}\what K_U^{(t+1)}\big\|+\|\what V^{(t+1)}-V\what O_V^{(t)}\what R_G^{(t)}\what K_V^{(t+1)} \|\\
\leq&C_3\eta\frac{\sigma_{\xi}}{\lambda_r}\cdot\sqrt{\frac{d_1^2d_2\log d_1}{N_0}}+\big(1-\frac{\eta}{2}\big)\cdot \big(\|\what U^{(t)}-U\what O_U^{(t)}\|+\|\what V^{(t)}-V\what O_V^{(t)}\|\big)
+2C_5\eta \mumax^2\kappa_0\sqrt{\frac{ r^3d_1}{N_0}}\\
&\hspace{2cm}\leq \frac{1}{3C_4\mumax\kappa_0^2\sqrt{r\alpha_d}},
\end{align*}
where the last inequality holds as long as $\eta\leq 0.75$ and $n\geq C_6\alpha_d\mumax^6\kappa_0^6r^3d_1\log^2d_1$ and $\mumax\kappa_0^2(\sigma_{\xi}/\lambda_r)\cdot\sqrt{\alpha_d r d_1^2d_2\log^2d_1/n}\leq C_7^{-1}$ for some large enough constants $C_6,C_7>0$. Then, it suffices to prove $\|\what U^{(1)}-U\what O_U^{(1)}\|+\|\what V^{(1)}-V\what O_V^{(1)}\|\leq 1/(3C_4\mumax\kappa_0^2\sqrt{r\alpha_d})$ where, by Davis-Kahan theorem, with probability at least $1-2d_1^{-2}$,
\begin{align*}
\|\what U^{(1)}\what U^{(1)\tran}-&UU^{\tran}\|+\|\what V^{(1)}\what V^{(1)\tran}-VV^{\tran}\|
\leq&C_4\frac{\sigma_{\xi}+\|M\|_{\submax}}{\lambda_r}\cdot \sqrt{\frac{d_1^2d_2\log d_1}{N_0}}\leq \frac{1}{3C_4\mumax\kappa_0^2\sqrt{r\alpha_d}},
\end{align*}
as long as $n\geq C_5\alpha_d\kappa_0^6\mumax^6r^3d_1\log^2d_1$ and $C_6\mumax\kappa_0^2(\sigma_{\xi}/\lambda_r)\cdot\sqrt{\alpha_d rd_1^2d_2\log^2d_1/n}\leq 1$. We then conclude the proof of the first statement of Theorem~\ref{thm:gd}. 

We now prove the second statement. Recall that $N_0\asymp n/\log d_1$, by the first statement with $\eta=0.75$, we get with probability at least $1-4md_1^{-2}$,
\begin{align*}
\big\|&\what U^{(m)}-U\what O_U^{(m)}\big\|_{\submaxx}+\big\|\what V^{(m)}-V\what O_V^{(m)}\big\|_{\submaxx}-2C_3\frac{\sigma_{\xi}}{\lambda_r}\sqrt{\frac{rd_1d_2\log^2d_1}{n}}\\
\leq&\Big(\frac{1}{2}\Big)^m\cdot\bigg(\big\|\what U^{(1)}-U\what O_U^{(1)}\big\|_{\submaxx}+\big\|\what V^{(1)}-V\what O_V^{(1)}\big\|_{\submaxx}-2C_3\frac{\sigma_{\xi}}{\lambda_r}\sqrt{\frac{rd_1d_2\log^2d_1}{n}}\bigg).
\end{align*}
Similar as the proof of Theorem~\ref{lem:hatUhatV_deloc}, with probability at least $1-d_1^{-2}$,
\begin{align*}
\|\what U^{(1)}-U\what O_U^{(1)}&\|_{\submaxx}+ \|\what V^{(1)}-V\what O_V^{(1)}\|_{\submaxx}
\leq& C_4\mumax\frac{\sigma_{\xi}}{\lambda_r}\sqrt{\frac{rd_1^2\log^2d_1}{n}}+C_5\mumax\frac{\|M\|_{\submax}}{\lambda_r}\sqrt{\frac{rd_1^2\log^2d_1}{n}}.
\end{align*}
Therefore, if $m=2\ceil{\log(\alpha_d\|M\|_{\submax}/\sigma_{\xi})}\leq 2C_1\ceil{\log d_1}$, we get 
\begin{align*}
\big\|\what U^{(m)}-U&\what O_U^{(m)}\big\|_{\submaxx}+\big\|\what V^{(m)}-V\what O_V^{(m)}\big\|_{\submaxx}\leq C_4\frac{\sigma_{\xi}}{\lambda_r}\sqrt{\frac{rd_1d_2\log^2d_1}{n}},
\end{align*}
which holds with probability at least $1-4C_1d_1^{-2}\log d_1$. Then, by Lemma~\ref{lem:hatG},
\begin{align*}
\big\|\what M^{(m)}-M \big\|_{\submax}\leq 2\|\Lambda\|\mumax\cdot\Big(\sqrt{\frac{r}{d_2}}\|\what U^{m}-U\what O_U^{(m)}\|_{\submaxx}+\sqrt{\frac{r}{d_1}}\|\what V^{(m)}-V\what O_V^{(m)}\|_{\submaxx}\Big)\\
+\mumax^2\sqrt{\frac{r^2}{d_1d_2}}\|\what G^{(m)}-\what O_U^{(m)\tran}\Lambda\what O_V^{(m)}\|
\leq C_3\mumax\kappa_0 \sigma_{\xi}\sqrt{\frac{r^2d_1\log^2d_1}{n}}.
\end{align*}

\subsection{Proof of Theorem~\ref{lem:hatUhatV_deloc}}
W.L.O.G., we only prove the bounds for $d_{\submaxx}(\what U_1, U)$ and $d_{\submaxx}(\what V_1,V)$.
To this end, define the $(d_1+d_2)\times (2r)$ matrices
$$
\Theta=\left(\begin{array}{cc} U& 0\\ 0&V\end{array}\right)\quad {\rm and}\quad \what\Theta_1=\left(\begin{array}{cc} \what U_1& 0\\ 0&\what V_1\end{array}\right).
$$
We also define the $(d_1+d_2)\times (d_1+d_2)$ matrices
$$
A=\left(\begin{array}{cc}   0&M\\ M^{\tran}&0\end{array}\right)\quad {\rm and}\quad \what{E}^{(1)}=\left(\begin{array}{cc}  0&\what{Z}^{(1)}\\ \what{Z}^{(1)\tran}&0\end{array}\right).
$$
Let $U_{\perp}\in\RR^{d_1\times (d_1-r)}$ and $V_{\perp}\in\RR^{d_2\times (d_2-r)}$ so that $(U_{\perp},U)$ and $(V_{\perp},V)$ are orthogonal matrices. For any positive integer $s\geq1$, we define 
$$
\mathfrak{P}^{-s}=
\begin{cases}
\left(\begin{array}{cc}U\Lambda^{-s}U^{\tran}&0\\0&V\Lambda^{-s}V^{\tran}\end{array}\right),&\textrm{ if } s \textrm{ is even};\\
\left(\begin{array}{cc}0&U\Lambda^{-s}V^{\tran}\\V\Lambda^{-s}U^{\tran}&0\end{array}\right),&\textrm{ if } s \textrm{ is odd}.
\end{cases}
$$
Define also
$$
\mathfrak{P}^0=\mathfrak{P}^{\perp}=\left(\begin{array}{cc}U_{\perp}U_{\perp}^{\tran}&0\\ 0&V_{\perp}V_{\perp}^{\tran}\end{array}\right).
$$
As shown by \cite{xia2019normal}, 
if $\lambda_r\geq 2\|\what Z^{(1)}\|$, then 
\begin{align}
\what{\Theta}_1\what\Theta^{\tran}_1-\Theta\Theta^{\tran}=\sum_{k=1}^{\infty}\underbrace{\sum_{\bs:s_1+\cdots+s_{k+1}=k}(-1)^{1+\tau(\bs)}\cdot \mathfrak{P}^{-s_1}\what E^{(1)}\mathfrak{P}^{-s_2}\cdots\mathfrak{P}^{-s_k}\what{E}^{(1)}\mathfrak{P}^{-s_{k+1}}}_{\calS_{A,k}(\what{E}^{(1)})},\label{lem:hatU-U}
\end{align}
where $s_1,\cdots,s_{k+1}\geq 0$ are integers and  $\tau(\bs)=\sum_{i=1}^{k+1}{\bf 1}(s_i>0)$. 
We aim to prove sharp upper bounds for  $\|\what U_1\what U_1^{\tran}-UU^{\tran}\|_{\submaxx}$ and $\|\what V_1\what V_1^{\tran}-VV^{\tran}\|_{\submaxx}$. Note that 
$$
\what\Theta_1\what\Theta^{\tran}_1-\Theta\Theta^{\tran}=\left(\begin{array}{cc}\what U_1\what U_1^{\tran}-UU^{\tran}&0\\0& \what V_1\what V_1^{\tran}-VV^{\tran} \end{array}\right).
$$
Therefore, it suffices to investigate $\|\what\Theta_1\what\Theta^{\tran}_1-\Theta\Theta^{\tran}\|_{\submaxx}$. By \eqref{lem:hatU-U}, we obtain
\begin{align*}
\|\what\Theta_1\what\Theta^{\tran}_1-\Theta\Theta^{\tran}\|_{\submaxx}\leq \sum_{k=1}^{\infty}\sum_{\bs:s_1+\cdots+s_{k+1}=k}\big\|\mathfrak{P}^{-s_1}\what E^{(1)}\mathfrak{P}^{-s_2}\cdots\mathfrak{P}^{-s_k}\what{E}^{(1)}\mathfrak{P}^{-s_{k+1}}\big\|_{\submaxx}.
\end{align*}
Denote $e_{j}$ the $j$-th canonical basis vector in $\RR^{d_1+d_2}$ for any $j\in[d_1+d_2]$.  
Recall Assumption~\ref{assump:incoh} and the definition of $\mathfrak{P}^{-s}$, it is obvious that for all $s\geq 1$,
$$
\max_{j\in[d_1]}\big\|e_{j}^{\tran}\mathfrak{P}^{-s} \big\|\leq \mumax\sqrt{\frac{r}{d_1}}\cdot \|\Lambda^{-s}\|\quad{\rm and}\quad \max_{j\in[d_2]}\big\|e_{j+d_1}^{\tran}\mathfrak{P}^{-s} \big\|\leq \mumax\sqrt{\frac{r}{d_2}}\cdot \|\Lambda^{-s}\|. 
$$
For any $(s_1,\cdots,s_{k+1})$ such that $\sum_{j=1}^{k+1}s_j=k$ and $s_1\geq 1$, we have 
\begin{align*}
\big\|e_j^{\tran}\mathfrak{P}^{-s_1}\what E^{(1)}\mathfrak{P}^{-s_2}\cdots&\mathfrak{P}^{-s_k}\what{E}^{(1)}\mathfrak{P}^{-s_{k+1}}\big\|\leq \|e_j^{\tran}\mathfrak{P}^{-s_1}\|\cdot \big\|\what E^{(1)}\mathfrak{P}^{-s_2}\cdots\mathfrak{P}^{-s_k}\what{E}^{(1)}\mathfrak{P}^{-s_{k+1}}\big\|\\
\leq& \|e_j^{\tran}\mathfrak{P}^{-s_1}\|\cdot \|\what{E}^{(1)}\|^k\|\Lambda^{-1}\|^{k-s_1}.
\end{align*}
By Lemma~\ref{lem:Zbound}, there exists an event $\calE_0$ with $\PP(\calE_0)\geq 1-2d_1^{-2}$ so that on $\calE_0$, 
\begin{equation}\label{eq:delta_def}
\big\|\what{E}^{(1)} \big\|\leq \underbrace{C_2(1+\gamma_n)\sigma_{\xi}\sqrt{\frac{d_1^2d_2\log d_1}{n}}}_{\delta}.
\end{equation}
Therefore, on event $\calE_0$, if $s_1\geq 1$, then
$$
\max_{j\in[d_1]}\big\|e_{j}^{\tran}\mathfrak{P}^{-s_1}\what E^{(1)}\mathfrak{P}^{-s_2}\cdots\mathfrak{P}^{-s_k}\what{E}^{(1)}\mathfrak{P}^{-s_{k+1}}\big\|\leq\Big(\frac{\delta}{\lambda_r}\Big)^k\cdot\mumax\sqrt{\frac{r}{d_1}}
$$
and
$$
\max_{j\in[d_2]}\big\|e_{j+d_1}^{\tran}\mathfrak{P}^{-s_1}\what E^{(1)}\mathfrak{P}^{-s_2}\cdots\mathfrak{P}^{-s_k}\what{E}^{(1)}\mathfrak{P}^{-s_{k+1}}\big\|\leq\Big(\frac{\delta}{\lambda_r}\Big)^k\cdot\mumax\sqrt{\frac{r}{d_2}},
$$
where $\delta$ is defined in (\ref{eq:delta_def}). 

As a result, it suffices to prove the upper bounds for $\|\mathfrak{P}^{-s_1}\what E^{(1)}\mathfrak{P}^{-s_2}\cdots\mathfrak{P}^{-s_k}\what{E}^{(1)}\mathfrak{P}^{-s_{k+1}}\|_{\submaxx}$ for $s_1=0$. Because $s_1+\cdots+s_{k+1}=k$, there must exists $s_j\geq 1$ for some $j\geq 2$. It then suffices to prove the upper bounds for $\|\mathfrak{P}^{\perp}(\mathfrak{P}^{\perp}\what{E}^{(1)}\mathfrak{P}^{\perp})^k\what{E}^{(1)}\Theta\|_{\submaxx}$ with $k\geq 0$.  Note that we used the fact $\Theta\Theta^{\tran}\mathfrak{P}^{-s}\Theta\Theta^{\tran}=\mathfrak{P}^{-s}$ for any integer $s\geq 1$. 
\begin{lemma}\label{lem:PE_2max}
Under the event $\calE_0$ where (\ref{eq:delta_def}) holds, there exist absolute constants $C_1,C_2>0$ so that, for all $k\geq 0$,  the following bounds hold with probability at least $1-2(k+1)d_1^{-2}$, 
$$
\max_{j\in[d_1]}\big\|e_{j}^{\tran}\mathfrak{P}^{\perp}(\mathfrak{P}^{\perp}\what{E}^{(1)}\mathfrak{P}^{\perp})^k\what{E}^{(1)}\Theta \big\|\leq C_1(C_2\delta)^{k+1}\cdot\mumax\sqrt{\frac{r}{d_1}},
$$
$$
\max_{j\in[d_2]}\big\|e_{j+d_1}^{\tran}\mathfrak{P}^{\perp}(\mathfrak{P}^{\perp}\what{E}^{(1)}\mathfrak{P}^{\perp})^k\what{E}^{(1)}\Theta \big\|\leq C_1(C_2\delta)^{k+1}\cdot\mumax\sqrt{\frac{r}{d_2}}
$$
where $\delta$ is defined in (\ref{eq:delta_def}) and $\mumax$ is the incoherence constant in Assumption~\ref{assump:incoh}.
\end{lemma}
We shall defer the proof of Lemma \ref{lem:PE_2max} to Appendix.

By Lemma~\ref{lem:PE_2max} and (\ref{eq:delta_def}), choosing $\kmax=\ceil{2\log d_1}$ yields that, for all $\bs=(s_1,\cdots,s_{k+1})$ with $\sum_{j=1}^{k+1}s_j=k$,
\begin{align*}
\max_{j\in[d_1]}\big\|e_{j}^{\tran}\mathfrak{P}^{-s_1}\what E^{(1)}\mathfrak{P}^{-s_2}\cdots\mathfrak{P}^{-s_k}\what{E}^{(1)}\mathfrak{P}^{-s_{k+1}} \big\|\leq C_1\Big(\frac{C_2\delta}{\lambda_r}\Big)^k\cdot \mumax\sqrt{\frac{r}{d_1}},
\end{align*}
which holds for all $k\leq \kmax$ with probability at least $1-4d_1^{-2}\log d_1$, under event $\calE_0$. Then,
\begin{align*}
\max_{j\in[d_1]}\big\|e_{j}^{\tran}(\what U_1\what U_1^{\tran}-UU^{\tran})\big\|=\max_{j\in[d_1]}\big\|e_j^{\tran}(\what\Theta\what\Theta^{\tran}-\Theta\Theta^{\tran})\big\|,
\end{align*}
where we abuse the notations that $e_{j}\in\RR^{d_1}$ on the left hand side and $e_j\in\RR^{d_1+d_2}$ on the right hand side. Then, by the representation formula \eqref{lem:hatU-U},
\begin{align*}
\max_{j\in[d_1]}\big\|e_{j}^{\tran}(\what U_1\what U_1^{\tran}-UU^{\tran})\big\|\leq& \max_{j\in[d_1]}\sum_{k=1}^{\kmax}\sum_{\bs:s_1+\cdots+s_{k+1}=k}\big\|e_j^{\tran} \mathfrak{P}^{-s_1}\what E^{(1)}\mathfrak{P}^{-s_2}\cdots\mathfrak{P}^{-s_k}\what{E}^{(1)}\mathfrak{P}^{-s_{k+1}} \big\|\\
&+\sum_{k=\kmax+1}^{\infty}\sum_{\bs:s_1+\cdots+s_{k+1}=k}\big\| \mathfrak{P}^{-s_1}\what E^{(1)}\mathfrak{P}^{-s_2}\cdots\mathfrak{P}^{-s_k}\what{E}^{(1)}\mathfrak{P}^{-s_{k+1}} \big\|_{\submaxx}.
\end{align*}
Obviously, 
$${\rm Card}\Big(\big\{(s_1,\cdots,s_{k+1}): \sum_{j=1}^{k+1}s_j=k, s_j\in\ZZ, s_j\geq 0\big\}\Big)\leq 4^k.
$$
Therefore, under event $\calE_0$, 
\begin{align*}
\max_{j\in[d_1]}\big\|e_{j}^{\tran}(\what U_1\what U_1^{\tran}-UU^{\tran})\big\|\leq& \sum_{k=1}^{\kmax} C_1\Big(\frac{4C_2\delta}{\lambda_r}\Big)^k\cdot\mumax\sqrt{\frac{r}{d_1}}
+\sum_{k=\kmax+1}^{\infty}4^k\cdot \Big(\frac{\delta}{\lambda_r}\Big)^k\\
\leq& C_1\sum_{k=1}^{\kmax}\Big(\frac{4C_2\delta}{\lambda_r}\Big)^k\cdot \mumax\sqrt{\frac{r}{d_1}}+\sum_{k=\kmax+1}^{\infty}\Big(\frac{4\delta}{\lambda_r}\Big)^k.
\end{align*}
If $8C_2\delta/ \lambda_r\leq 1$ and $C_2>4$, then 
\begin{align*}
\max_{j\in[d_1]}\big\|e_{j}^{\tran}(\what U_1\what U_1^{\tran}-UU^{\tran})\big\|\leq& C_1\frac{\delta}{\lambda_r}\cdot \mumax\sqrt{\frac{r}{d_1}}+2\Big(\frac{4\delta}{\lambda_r}\Big)^{\kmax+1}\\
\leq&C_1\frac{\delta}{\lambda_r}\cdot \mumax\sqrt{\frac{r}{d_1}}+\frac{8\delta}{\lambda_r}\cdot \Big(\frac{1}{2C_2}\Big)^{\ceil{2\log d_1}}\\
\leq& C_2\frac{(1+\gamma_n)\sigma_{\xi}}{\lambda_r}\sqrt{\frac{d_1^2d_2\log d_1}{n}}\cdot \mumax\sqrt{\frac{r}{d_1}}. 
\end{align*}
Therefore, 
$$
\PP\Big(\big\|\what U_1\what U_1^{\tran}-UU^{\tran} \big\|_{\submaxx}\leq C_2\frac{(1+\gamma_n)\sigma_{\xi}}{\lambda_r}\sqrt{\frac{d_1^2d_2\log d_1}{n}}\cdot \mumax\sqrt{\frac{r}{d_1}}\Big)\geq 1-5d_1^{-2}\log d_1.
$$
Similarly, on the same event, 
$$
\big\|\what V_1\what V_1^{\tran}-VV^{\tran} \big\|_{\submaxx}\leq C_2\frac{(1+\gamma_n)\sigma_{\xi}}{\lambda_r}\sqrt{\frac{d_1^2d_2\log d_1}{n}}\cdot \mumax\sqrt{\frac{r}{d_2}},
$$
which proves the claimed bound.

\bibliographystyle{plainnat}
\bibliography{refer}

\newpage
\appendix


\section{Proof of Lemma~\ref{lem:Zbound}}
W.L.O.G., we only prove the bounds for $\|\what{Z}_1^{(1)}\|$ and $\|\what{Z}_2^{(1)}\|$. 
Recall that $\what{Z}_1^{(1)}$ is defined by
$$
\what{Z}_1^{(1)}=\frac{d_1d_2}{n_0}\sum_{i=n_0+1}^n \xi_i X_i
$$
where $\{(\xi_i,X_i)\}_{i=n_0+1}^n$ are i.i.d. The $\psi_\alpha$-norm of a random variable $Y$ is defined by $\|Y\|_{\psi_\alpha}=\inf\{t>0: \EE\exp^{|Y/t|^{\alpha}}\leq 2\}$ for $\alpha\in[1,2]$. Since $\xi$ is sub-Gaussian, we obtain $\|\xi\|_{\psi_2}\lesssim \sigma_{\xi}$. Clearly,
$$
\big\|\|\xi_i X_i\| \big\|_{\psi_2}\leq \|\xi_i\|_{\psi_2}\lesssim \sigma_{\xi}
$$
where we used the fact $X_i\in\mathfrak{E}=\{e_{j_1}e_{j_2}^{\tran}: j_1\in[d_1], j_2\in[d_2]\}$. Meanwhile, 
\begin{align*}
\big\|\EE \xi_i^2 X_iX_i^{\tran}\big\|=\bigg\|\sigma_{\xi}^2\cdot \frac{1}{d_1d_2}\sum_{j_1=1}^{d_1}\sum_{j_2=1}^{d_2} e_{j_1}e_{j_2}^{\tran}e_{j_2}e_{j_1}^{\tran} \bigg\|=\big\|\frac{\sigma_{\xi}^2}{d_1}\cdot I_{d_1} \big\|\leq \frac{\sigma_{\xi}^2}{d_1}.
\end{align*}
Similar bounds also hold for $\|\EE\xi_i^2 X_i^{\tran}X_i\|$ and we conclude with 
$$
\max\big\{\big\|\EE \xi_i^2 X_iX_i^{\tran}\big\|,  \|\EE\xi_i^2 X_i^{\tran}X_i\|\big\}\leq \frac{\sigma_{\xi}^2}{d_2}. 
$$
By matrix Bernstein inequality (\cite{koltchinskii2011neumann, minsker2017some, tropp2012user}), for all $t>0$, the following bound holds with probability at least $1-e^{-t}$,
$$
\big\|\what{Z}_1^{(1)}\big\|\leq C_1\sigma_{\xi}\sqrt{\frac{d_1^2d_2(t+\log d_1)}{n}}+C_2\sigma_{\xi}\frac{d_1d_2(t+\log d_1)}{n}.
$$
By setting $t=2\log d_1$ and the fact $n\geq C_3d_1\log d_1$, we conclude with 
$$
\PP\bigg(\|\what{Z}_1^{(1)}\|\geq C_1\sigma_{\xi}\sqrt{\frac{d_1^2d_2(t+\log d_1)}{n}}\bigg)\leq \frac{1}{d_1^2}.
$$
The upper bound for $\|\what{Z}_2^{(1)}\|$ can be derived in the same fashion by observing that 
$$
\big\|d_1d_2\langle\what\Delta_1,X_i \rangle X_i -\what\Delta_1\big\|\leq d_1d_2\|\what\Delta_1\|_{\submax}+\|\what\Delta_1\|\leq 2d_1d_2\|\what\Delta_1\|_{\submax}
$$
and
\begin{align*}
\big\|\EE\big(d_1d_2\langle\what\Delta_1,&X_i \rangle X_i -\what\Delta_1\big) \big(d_1d_2\langle\what\Delta_1,X_i \rangle X_i -\what\Delta_1\big)^{\tran}\big\|\\
\leq&\big\|d_1^2d_2^2\EE\langle\what\Delta_1, X_i \rangle^2X_iX_i^{\tran} \big\|+\|\what\Delta_1\|^2\leq d_1^2d_2\|\what\Delta_1\|_{\submax}^2+\|\what\Delta_1\|^2\leq 2d_1^2d_2\|\what\Delta_1\|_{\submax}^2.
\end{align*}

\section{Proof of Lemma~\ref{lemma:hatUhatV}}
W.L.O.G., we only prove the upper bounds for $(\what U_1, \what V_1)$ since the proof for $(\what U_2, \what V_2)$  is identical. 
Recall from Assumption~\ref{assump:init_entry} that $\|\what\Delta_1\|_{\submax}\leq C_1\gamma_n\cdot \sigma_{\xi}$ with probability at least $1-d_1^{-2}$. 
To this end, we conclude with 
$$
\PP\bigg(\|\what{Z}^{(1)}\|\geq  C_2(1+\gamma_n)\sigma_{\xi}\sqrt{\frac{d_1^2d_2\log d_1}{n}}\bigg)\leq \frac{2}{d_1^2}. 
$$
Recall that $\what{M}_1^{\supunbiased}=M+\what{Z}^{(1)}$. By Davis-Kahan Theorem (\cite{davis1970rotation}) or Wedin's $\sin\Theta$ Theorem (\cite{wedin1972perturbation}), we get 
$$
\max\{d_{\subo}(\what U_1, U), d_{\subo}(\what V_1, V)\}\leq \frac{\sqrt{2}\|\what{Z}^{(1)}\|}{\lambda_r}\leq C_2\frac{\sqrt{2}(1+\gamma_n)\sigma_{\xi}}{\lambda_r}\cdot \sqrt{\frac{d_1^2d_2\log d_1}{n}}
$$
where the last inequality holds with probability at least $1-2d_1^{-2}$. Similarly, with the same probability, 
$$
\max\{d_{\subf}(\what U_1, U), d_{\subf}(\what V_1, V)\}\leq \frac{\sqrt{2r}\|\what{Z}^{(1)}\|}{\lambda_r}\leq C_2\frac{\sqrt{2}(1+\gamma_n)\sigma_{\xi}}{\lambda_r}\cdot \sqrt{\frac{rd_1^2d_2\log d_1}{n}}
$$
which concludes the proof of Lemma~\ref{lemma:hatUhatV}.

\section{Proof of Lemma~\ref{lem:PE_2max}}
For notational simplicity, we write $\what E=\what E^{(1)}$ in this section.

\subsection{Case 0: $k=0$}
W.L.O.G., we bound $\|e_{j}^{\tran}\mathfrak{P}^{\perp}\what{E}\Theta\|$ for $j\in[d_1]$. Clearly, 
\begin{align*}
\|e_{j}^{\tran}\mathfrak{P}^{\perp}\what{E}\Theta\|\leq&\|e_j^{\tran}\Theta\Theta^{\tran}\what{E}\Theta\|+\|e_{j}^{\tran}\what{E}\Theta\|
\leq \delta\mumax\cdot\sqrt{\frac{r}{d_1}}+\|e_{j}^{\tran}\what{E}\Theta\|
\end{align*}
where $\delta$ denotes the upper bound of $\|\what{E}\|$ defined in (\ref{eq:delta_def}) and the last inequality is due to  $\|U\|_{\submaxx}\leq \mumax\sqrt{r/d_1}$.  By the definitions of $\what{E}$ and $\Theta$,  $\|e_{j}^{\tran}\what{E}\Theta\|=\|e_j^{\tran}\what{Z}^{(1)}V\|$ where we abuse the notations and denote $e_j$ the canonical basis vectors in $\RR^{d_1}$. 

Recall that $\what{Z}^{(1)}=\what{Z}_1^{(1)}+\what{Z}_2^{(1)}$.  We write 
\begin{align*}
e_{j}^{\tran}\what{Z}^{(1)}V=\frac{d_1d_2}{n_0}\sum_{i=n_0+1}^n \xi_i e_{j}^{\tran}X_i V+\Big(\frac{d_1d_2}{n_0}\sum_{i=n_0+1}^n\langle\what\Delta_1,X_i \rangle e_{j}^{\tran}X_i V-e_{j}^{\tran}\what\Delta_1 V\Big).
\end{align*}
Clearly, 
$$
\big\|\|\xi_i e_j^{\tran}X_iV\|\big\|_{\psi_2}\leq \sigma_{\xi} \|V\|_{\submaxx}\leq \sigma_{\xi}\mumax\cdot\sqrt{\frac{r}{d_2}}
$$
and 
$$
\EE \xi_i^2 e_j^{\tran}X_i VV^{\tran}X_i^{\tran}e_j\leq \frac{\sigma_{\xi}^2}{d_1d_2}\cdot \tr(VV^{\tran})\leq \frac{r\sigma_{\xi}^2}{d_1d_2}.
$$
Then, by Bernstein inequality, we get
\begin{align*}
\PP\bigg(\big\|e_{j}^{\tran}\what{Z}_1^{(1)}V \big\|\geq  C_1\sigma_{\xi}\sqrt{\frac{rd_1d_2(t+\log d_1)}{n}}+C_2\mumax\sigma_{\xi}\frac{d_1\sqrt{rd_2}(t+\log d_1)}{n}\bigg)\leq e^{-t}
\end{align*}
for all $t>0$ and some absolute constants $C_1,C_2>0$.  Similarly, 
\begin{align*}
\PP\bigg(\big\|e_{j}^{\tran}\what{Z}_2^{(1)}V \big\|\geq  C_1\|\what\Delta_1\|_{\submax}\sqrt{\frac{rd_1d_2(t+\log d_1)}{n}}+C_2\mumax\|\what\Delta_1\|_{\submax}\frac{d_1\sqrt{rd_2}(t+\log d_1)}{n}\bigg)\\
\leq e^{-t}
\end{align*}
By setting $t=3\log d_1$ and observing $\|\what\Delta_1\|_{\submax}\leq C_1\gamma_n\cdot\sigma_{\xi}$, we conclude that 
$$
\|e_j^{\tran}\what{Z}^{(1)}V\|\leq  C_1(1+\gamma_n)\sigma_{\xi}\sqrt{\frac{rd_1d_2\log d_1}{n}}=\delta\cdot \sqrt{\frac{r}{d_1}}
$$
which holds with probability at least $1-2d_1^{-3}$ and we used the assumption $n\geq C\mumax^2 rd_1\log d_1$ for some large enough constant $C>0$. As a result, 
$$
\PP\Big(\max_{j\in[d_1]}\|e_j^{\tran}\mathfrak{P}^{\perp}\what{E}\Theta\|\leq 2\delta\mumax\cdot \sqrt{\frac{r}{d_1}}\Big)\geq 1-2d_1d_1^{-3}.
$$
Following the same arguments, we can prove the bound for $\max_{j\in[d_2]}\|e_{d_1+j}^{\tran}\mathfrak{P}^{\perp}\what{E}\Theta\|$. Therefore, with probability at least $1-2d_1^{-2}$,
$$
\max_{j\in[d_1]}\|e_j^{\tran}\mathfrak{P}^{\perp}\what{E}\Theta\|\leq 2\delta\mumax\cdot \sqrt{\frac{r}{d_1}}\quad {\rm and}\quad \max_{j\in[d_2]}\|e_{d_1+j}^{\tran}\mathfrak{P}^{\perp}\what{E}\Theta\|\leq 2\delta\mumax\cdot\sqrt{\frac{r}{d_2}}
$$
where $\delta$ is defined by (\ref{eq:delta_def}).

\subsection{Case 1: $k=1$}
W.L.O.G., we bound $\max_{j\in[d_1]}\big\|e_{j}^{\tran}\mathfrak{P}^{\perp}\what{E}\mathfrak{P}^{\perp}\what{E}\Theta\big\|$.
Observe that 
\begin{align*}
\big\|e_{j}^{\tran}\mathfrak{P}^{\perp}\what{E}\mathfrak{P}^{\perp}\what{E}\Theta\big\|\leq& \big\|e_{j}^{\tran}\Theta\Theta^{\tran}\what{E}\mathfrak{P}^{\perp}\what{E}\Theta\big\|+\big\|e_{j}^{\tran}\what{E}\mathfrak{P}^{\perp}\what{E}\Theta\big\|
\leq& \delta^2\mumax\cdot\sqrt{\frac{r}{d_1}}+\big\|e_{j}^{\tran}\what{E}\mathfrak{P}^{\perp}\what{E}\Theta\big\|.
\end{align*}
By the definition of $\what{E}$ and $\mathfrak{P}^{\perp}$, we have
$$
\what{E}\mathfrak{P}^{\perp}\what{E}\Theta=\left(\begin{array}{cc}\what{Z}^{(1)}V_{\perp}V_{\perp}^{\tran}\what{Z}^{(1)\tran}U&0\\0&\what{Z}^{(1)\tran}U_{\perp}U_{\perp}^{\tran}\what{Z}^{(1)}V \end{array}\right).
$$
It suffices to prove the upper bound for $\big\|e_{j}^{\tran}\what{Z}^{(1)}V_{\perp}V_{\perp}^{\tran}\what{Z}^{(1)\tran}U\big\|$. 
Define $\mathfrak{I}_{j}=e_{j}e_{j}^{\tran}\in\RR^{d_1\times d_1}$ and $\mathfrak{I}_{j}^{\perp}=\calI-\mathfrak{I}_{j}$. Then, write $\what{Z}^{(1)}=\mathfrak{I}_{j}\what{Z}^{(1)}+\mathfrak{I}_{j}^{\perp}\what{Z}^{(1)}$ and 
\begin{align*}
e_{j}^{\tran}\what{Z}^{(1)}V_{\perp}V_{\perp}^{\tran}\what{Z}^{(1)\tran}U=e_{j}^{\tran}\what{Z}^{(1)}V_{\perp}V_{\perp}^{\tran}\big(\mathfrak{I}_{j}\what{Z}^{(1)}\big)^{\tran}U+e_{j}^{\tran}\what{Z}^{(1)}V_{\perp}V_{\perp}^{\tran}\big(\mathfrak{I}_{j}^{\perp}\what{Z}^{(1)}\big)^{\tran}U.
\end{align*}
As a result, 
\begin{align*}
\big\|e_{j}^{\tran}\what{Z}^{(1)}V_{\perp}V_{\perp}^{\tran}\what{Z}^{(1)\tran}U\big\|\leq& \big\|e_{j}^{\tran}\what{Z}^{(1)}V_{\perp}V_{\perp}^{\tran}\big(\mathfrak{I}_{j}\what{Z}^{(1)}\big)^{\tran}U\big\|+\big\|e_{j}^{\tran}\what{Z}^{(1)}V_{\perp}V_{\perp}^{\tran}\big(\mathfrak{I}_{j}^{\perp}\what{Z}^{(1)}\big)^{\tran}U\big\|\\
\leq& \|e_{j}^{\tran}U\|\|e_{j}^{\tran}\what{Z}^{(1)}V_{\perp}V_{\perp}^{\tran}\what{Z}^{(1)\tran}e_j\|+\big\|e_{j}^{\tran}\what{Z}^{(1)}V_{\perp}V_{\perp}^{\tran}\big(\mathfrak{I}_{j}^{\perp}\what{Z}^{(1)}\big)^{\tran}U\big\|\\
\leq& \delta^2\mumax\cdot\sqrt{\frac{r}{d_1}}+\big\|e_{j}^{\tran}\what{Z}^{(1)}V_{\perp}V_{\perp}^{\tran}\big(\mathfrak{I}_{j}^{\perp}\what{Z}^{(1)}\big)^{\tran}U\big\|.
\end{align*}
Recall
$$
\what{Z}^{(1)}=\frac{d_1d_2}{n_0}\sum_{i=n_0+1}^{n}\xi_i X_i+\Big(\frac{d_1d_2}{n_0}\sum_{i=n_0+1}^n \langle \what\Delta_1, X_i\rangle X_i-\what\Delta_1\Big).
$$
Define 
$$
\calN_{j}=\big\{n_0+1\leq i\leq n: e_{j}^{\tran}X_i\neq 0\big\}\quad {\rm and}\quad \calN_{j}^{\scriptscriptstyle c}=\big\{n_0+1\leq i\leq n: i\notin \calN_{j}\big\}.
$$
By Chernoff bound, we get that if $n\geq C_1d_1\log d_1$ for a large enough absolute constant $C_1>0$, then 
\begin{equation}\label{eq:event1}
\PP\Big(\bigcap_{j=1}^{d_1}\Big\{\frac{n_0}{2d_1}\leq \big|\calN_{j}\big|\leq \frac{2n_0}{d_1}\Big\}\Big)\geq 1-e^{-c_1n/d_1}
\end{equation}
for some absolute constant $c_1>0$. Denote the above event by $\calE_{1}$ with $\PP\big(\calE_{1}\big)\geq 1-e^{-c_1n/d_1}$.

We now prove the upper bound $\big\|e_{j}^{\tran}\what{Z}^{(1)}V_{\perp}V_{\perp}^{\tran}\big(\mathfrak{I}_{j}^{\perp}\what{Z}^{(1)}\big)^{\tran}U\big\|$ conditioned on $\calN_{j}$.  To this end, by the definitions of $\mathfrak{I}_j$ and $\mathfrak{I}_j^{\perp}$, we write
\begin{align*}
e_{j}^{\tran}\what{Z}^{(1)}V_{\perp}V_{\perp}^{\tran}\big(\mathfrak{I}_{j}^{\perp}\what{Z}^{(1)}\big)^{\tran}U=&\frac{d_1d_2}{n_0}\sum_{i\in\calN_{j}} \xi_i e_{j}^{\tran}X_iV_{\perp}V_{\perp}^{\tran}\big(\mathfrak{I}_{j}^{\perp}\what{Z}^{(1)}\big)^{\tran}U\\
+&\Big[\frac{d_1d_2}{n_0}\sum_{i\in\calN_{j}}\langle\what\Delta_1,X_i \rangle e_{j}^{\tran}X_iV_{\perp}V_{\perp}^{\tran}\big(\mathfrak{I}_{j}^{\perp}\what{Z}^{(1)}\big)^{\tran}U-e_{j}^{\tran}\what\Delta_1 V_{\perp}V_{\perp}^{\tran}\big(\mathfrak{I}_{j}^{\perp}\what{Z}^{(1)}\big)^{\tran}U\Big].
\end{align*}
Note that, conditioned on $\calN_{j}$, $\{\xi_i, X_i\}_{i\in\calN_{j}}$ are independent with $\mathfrak{I}_{j}^{\perp}\what{Z}^{(1)}$. Conditioned on $\calN_{j}$ and $\mathfrak{I}_{j}^{\perp}\what{Z}^{(1)}$, the following facts are obvious. 
\begin{align*}
\Big\|\xi_i \big\|e_{j}^{\tran}X_iV_{\perp}V_{\perp}^{\tran}\big(\mathfrak{I}_{j}^{\perp}\what{Z}^{(1)}\big)^{\tran}U\big\| \Big\|_{\psi_2}\leq& \sigma_{\xi}\big\|V_{\perp}V_{\perp}^{\tran}\big(\mathfrak{I}_{j}^{\perp}\what{Z}^{(1)}\big)^{\tran}U\big\|_{\submaxx}.
\end{align*}
By the results of {\it Case 0} when $k=0$, we have 
\begin{align*}
\big\|V_{\perp}V_{\perp}^{\tran}\big(\mathfrak{I}_{j}^{\perp}\what{Z}^{(1)}\big)^{\tran}U\big\|_{\submaxx}\leq& \big\|V_{\perp}V_{\perp}^{\tran}\big(\mathfrak{I}_{j}\what{Z}^{(1)}\big)^{\tran}U\big\|_{\submaxx}+ \big\|V_{\perp}V_{\perp}^{\tran}\what{Z}^{(1)\tran}U\big\|_{\submaxx}\\
\leq&\delta\mumax\cdot\sqrt{\frac{r}{d_1}}+\delta\mumax\cdot\sqrt{\frac{r}{d_2}}\leq 2\delta\mumax\cdot\sqrt{\frac{r}{d_2}}.
\end{align*}
Meanwhile, (note that conditioned on $i\in\calN_j$, $X_i\overset{{\rm d}}{=}e_{j}e_{k}^{\tran}$ with $k$ being uniformly distributed over $[d_2]$)
\begin{align*}
\EE \xi_i^2 e_{j}^{\tran}X_iV_{\perp}V_{\perp}^{\tran}&\big(\mathfrak{I}_{j}^{\perp}\what{Z}^{(1)}\big)^{\tran}UU^{\tran}\big(\mathfrak{I}_{j}^{\perp}\what{Z}^{(1)}\big)V_{\perp}V_{\perp}^{\tran}X_i^{\tran}e_{j}\Big| \big(\mathfrak{I}_{j}^{\perp}\what{Z}^{(1)}\big), i\in\calN_j\\
&=\frac{\sigma_{\xi}^2}{d_2}\|V_{\perp}V_{\perp}^{\tran}(\mathfrak{I}_j^{\perp}\what{Z}^{(1)})^{\tran}U\|_{\subf}^2\leq \frac{r\sigma_{\xi}^2}{d_2}\cdot\big\|\big(\mathfrak{I}_{j}^{\perp}\what{Z}^{(1)}\big)^{\tran}U\big\|^2
\end{align*}
By Bernstein inequality, for all $t>0$, we get 
\begin{align*}
\PP\bigg(\Big\|\sum_{i\in\calN_{j}}\xi_i e_{j}^{\tran}X_iV_{\perp}&V_{\perp}^{\tran}\big(\mathfrak{I}_{j}^{\perp}\what{Z}^{(1)}\big)^{\tran}U\Big\|\geq C_1\sigma_{\xi}|\calN_{j}|^{1/2}\sqrt{\frac{r(t+\log d_1)}{d_2}}\big\|\big(\mathfrak{I}_{j}^{\perp}\what{Z}^{(1)}\big)^{\tran}U\big\|\\
+&C_2(t+\log d_1)\sigma_{\xi}\cdot\big\|V_{\perp}V_{\perp}^{\tran}\big(\mathfrak{I}_{j}^{\perp}\what{Z}^{(1)}\big)^{\tran}U\big\|_{\scriptscriptstyle \sf 2,max}
\bigg|\mathfrak{I}_{j}^{\perp}\what{Z}^{(1)}, \calN_{j}\bigg)\geq 1-e^{-t}
\end{align*}
for some absolute constants $C_1,C_2>0$. 

On event $\calE_0$, 
$$
\big\|\big(\mathfrak{I}_{j}^{\perp}\what{Z}^{(1)}\big)^{\tran}U\big\|\leq \|\what{Z}^{(1)}\|\leq \delta.
$$
By setting $t=3\log d_1$,  then with probability at least $1-d_1^{-3}$,
\begin{align*}
\Big\|\frac{d_1d_2}{n_0}\sum_{i=n_0+1}^n\xi_i e_{j}^{\tran}X_iV_{\perp}V_{\perp}^{\tran}&\big(\mathfrak{I}_{j}^{\perp}\what{Z}^{(1)}\big)^{\tran}U\Big\|\\
\leq& C_1\delta\cdot \sigma_{\xi}\sqrt{\frac{rd_1d_2\log d_1}{n}}+C_2\delta\cdot \mumax\sigma_{\xi}\frac{\sqrt{rd_2d_1^2\log^2d_1}}{n}\\
\leq& 2C_1\delta\cdot \sigma_{\xi}\sqrt{\frac{rd_1d_2\log d_1}{n}}\leq C_2\delta^2\cdot \sqrt{\frac{r}{d_1}}
\end{align*}
conditioned on $\calE_0\cap \calE_1$. The second inequality holds as long as $n\geq C\mumax^2d_1\log d_1$ for some large enough constant $C>0$. 

Similarly, since $\what\Delta_1$ is independent with $\{(X_i,\xi_i)\}_{i=n_0+1}^n$, with probability at least $1-d_1^{-3}$, 
$$
\Big\|\frac{d_1d_2}{n_0}\sum_{i=n_0+1}^n\langle\what\Delta_1,X_i \rangle e_{j}^{\tran}X_iV_{\perp}V_{\perp}^{\tran}\big(\mathfrak{I}_{j}^{\perp}\what{Z}^{(1)}\big)^{\tran}U-e_{j}^{\tran}\what\Delta_1 V_{\perp}V_{\perp}^{\tran}\big(\mathfrak{I}_{j}^{\perp}\what{Z}^{(1)}\big)^{\tran}U\Big\|\leq C_2\delta^2\cdot \sqrt{\frac{r}{d_1}}
$$
as long as $\|\what\Delta_1\|_{\submax}\leq \sigma_{\xi}$.  
Therefore, conditioned on $\calE_0\cap\calE_1$, with probability at least $1-2d_1^{-3}$,  
$$
\big\|e_{j}^{\tran}\what{Z}^{(1)}V_{\perp}V_{\perp}^{\tran}\big(\mathfrak{I}_{j}^{\perp}\what{Z}^{(1)}\big)^{\tran}U\big\|\leq C_2\delta^2\cdot \sqrt{\frac{r}{d_1}}.
$$
Therefore,  conditioned on $\calE_0\cap\calE_1$, 
\begin{align*}
\PP\bigg(\max_{j\in[d_1]}\|e_{j}^{\tran}\mathfrak{P}^{\perp}\what{E}\mathfrak{P}^{\perp}\what{E}\Theta\|\geq C_2\delta^2\cdot \mumax\sqrt{\frac{r}{d_1}}\bigg)\leq 2d_1^{-2}.
\end{align*}
Finally, conditioned on $\calE_0\cap\calE_1$, with probability at least $1-4d_1^{-2}$, 
$$
\max_{j\in[d_1]}\|e_{j}^{\tran}\mathfrak{P}^{\perp}\what{E}\mathfrak{P}^{\perp}\what{E}\Theta\|\leq C_2\delta^2\cdot \mumax\sqrt{\frac{r}{d_1}}\quad {\rm and}\quad \max_{j\in[d_2]}\|e_{d_1+j}^{\tran}\mathfrak{P}^{\perp}\what{E}\mathfrak{P}^{\perp}\what{E}\Theta\|\leq C_2\delta^2\cdot \mumax\sqrt{\frac{r}{d_2}}.
$$

\subsection{General $k\geq 2$}
{\it (Induction Assumption)} Suppose that for all $1\leq k_0\leq k$ with $k\geq 2$,  the following bounds hold, under events $\calE_0\cap \calE_1$, with probability at least $1-2kd_1^{-2}$
\begin{equation}\label{eq:induct1}
\max_{j\in[d_1]}\big\|e_{j}^{\tran}\mathfrak{P}^{\perp}\big(\mathfrak{P}^{\perp}\what{E}\mathfrak{P}^{\perp}\big)^{k_0-1}\what{E}\Theta\big\|\leq C_1(C_2\delta)^{k_0}\cdot \mumax\sqrt{\frac{r}{d_1}}
\end{equation}
and
\begin{equation}\label{eq:induct2}
\max_{j\in[d_2]}\big\|e_{j+d_1}^{\tran}\mathfrak{P}^{\perp}\big(\mathfrak{P}^{\perp}\what{E}\mathfrak{P}^{\perp}\big)^{k_0-1}\what{E}\Theta\big\|\leq C_1(C_2\delta)^{k_0}\cdot \mumax\sqrt{\frac{r}{d_2}}
\end{equation}
where $C_1, C_2>0$ are some absolute constants.

Based on the {\it Induction Assumption}, 
we  prove the upper bound for $\big\|\big(\mathfrak{P}^{\perp}\what{E}\mathfrak{P}^{\perp}\big)^{k}\what{E}\Theta\big\|_{\scriptscriptstyle\sf 2,max}$. 
W.O.L.G, we consider $\|e_{j}^{\tran}\big(\mathfrak{P}^{\perp}\what{E}\mathfrak{P}^{\perp}\big)^{k}\what{E}\Theta\big\|$ for any $j\in[d_1]$. To this end,  define the dilation operator $\mathfrak{D}$ so that 
$$
\mathfrak{D}(M)=\left(\begin{array}{cc}
0&M\\
M^{\tran}&0
\end{array}
\right).
$$
Then, $\what{E}=\mathfrak{D}(\what{Z}^{(1)})$. Similarly, define the following projectors on $\what{E}$, 
$$
\calP_{j}(\what{E})=\mathfrak{D}\big(e_{j}e_{j}^{\tran}\what{Z}^{(1)}\big)\quad {\rm and}\quad \calP_{j}^{\perp}(\what{E})=\mathfrak{D}\big(\mathfrak{I}_{j}^{\perp}\what{Z}^{(1)}\big).
$$
On event $\calE_0$,
\begin{align*}
\big\|e_{j}^{\tran}\big(\mathfrak{P}^{\perp}\what{E}&\mathfrak{P}^{\perp}\big)^{k}\what{E}\Theta \big\|=\big\|e_{j}^{\tran}\mathfrak{P}^{\perp}\what{E}\big(\mathfrak{P}^{\perp}\what{E}\mathfrak{P}^{\perp}\big)^{k-1}\what{E}\Theta \big\|\\
\leq&\big\|e_{j}^{\tran}\Theta\Theta^{\tran}\what{E}\big(\mathfrak{P}^{\perp}\what{E}\mathfrak{P}^{\perp}\big)^{k-1}\what{E}\Theta \big\|+\big\|e_{j}^{\tran}\what{E}\big(\mathfrak{P}^{\perp}\what{E}\mathfrak{P}^{\perp}\big)^{k-1}\what{E}\Theta \big\|\\
\leq& \delta^{k+1}\cdot \mumax\sqrt{\frac{r}{d_1}}+\big\|e_{j}^{\tran}\what{E}\big(\mathfrak{P}^{\perp}\what{E}\mathfrak{P}^{\perp}\big)^{k-1}\what{E}\Theta \big\|.
\end{align*}
We then write
\begin{align*}
e_{j}^{\tran}\what{E}\big(\mathfrak{P}^{\perp}\what{E}\mathfrak{P}^{\perp}\big)^{k-1}\what{E}\Theta=e_{j}^{\tran}\what{E}&\mathfrak{P}^{\perp}\calP_{j}(\what{E})\mathfrak{P}^{\perp}\big(\mathfrak{P}^{\perp}\what{E}\mathfrak{P}^{\perp}\big)^{k-2}\what{E}\Theta\\
+&e_{j}^{\tran}\what{E}\mathfrak{P}^{\perp}\calP_{j}^{\perp}(\what{E})\mathfrak{P}^{\perp}\big(\mathfrak{P}^{\perp}\what{E}\mathfrak{P}^{\perp}\big)^{k-2}\what{E}\Theta.
\end{align*}
By the {\it Induction Assumption},  under event $\calE_0$, 
\begin{align*}
\big\| e_{j}^{\tran}\what{E}\mathfrak{P}^{\perp}\calP_{j}(\what{E})\mathfrak{P}^{\perp}\big(\mathfrak{P}^{\perp}\what{E}\mathfrak{P}^{\perp}\big)^{k-2}\what{E}\Theta\big\|\leq& \|\what{E}\|\cdot \big\|\calP_{j}(\what{E})\mathfrak{P}^{\perp}\big(\mathfrak{P}^{\perp}\what{E}\mathfrak{P}^{\perp}\big)^{k-2}\what{E}\Theta\big\|_{\submaxx}\\
\leq& C_1(C_2\delta)^{k}\delta\cdot \mumax\sqrt{\frac{r}{d_1}}.
\end{align*}
Similarly, 
\begin{align*}
e_{j}^{\tran}\what{E}\mathfrak{P}^{\perp}&\calP_{j}^{\perp}(\what{E})\mathfrak{P}^{\perp}\big(\mathfrak{P}^{\perp}\what{E}\mathfrak{P}^{\perp}\big)^{k-2}\what{E}\Theta\\
=&e_{j}^{\tran}\what{E}\mathfrak{P}^{\perp}\calP_{j}^{\perp}(\what{E})\mathfrak{P}^{\perp}\big(\mathfrak{P}^{\perp}\what{E}\big)^{k-1}\Theta\\
=&e_{j}^{\tran}\what{E}\big(\mathfrak{P}^{\perp}\calP_{j}^{\perp}(\what{E})\big)^{k}\Theta+\sum_{t=1}^{k-1}e_{j}^{\tran}\what{E}\big(\mathfrak{P}^{\perp}\calP_{j}^{\perp}(\what{E})\big)^t\mathfrak{P}^{\perp}\calP_{j}(\what{E})\big(\mathfrak{P}^{\perp}\what{E}\big)^{k-1-t}\Theta.
\end{align*}
By the {\it Induction Assumption} and under event $\calE_0$,  
\begin{align*}
\big\|e_{j}^{\tran}\what{E}\big(\mathfrak{P}^{\perp}\calP_{j}^{\perp}(&\what{E})\big)^t\mathfrak{P}^{\perp}\calP_{j}(\what{E})\big(\mathfrak{P}^{\perp}\what{E}\big)^{k-1-t}\Theta \big\|\\
\leq & \|\what{E}\|^{t+1}\big\|\calP_{j}(\what{E})\big(\mathfrak{P}^{\perp}\what{E}\big)^{k-1-t}\Theta \big\|_{\scriptscriptstyle\sf 2,max}\leq C_1 (C_2\delta)^{k-t}\delta^{t+1}\cdot\mumax\sqrt{\frac{r}{d_1}}
\end{align*}
which holds for all $1\leq t\leq k-1$. Therefore, we conclude that on event $\calE_0\cap \calE_1$, with probability at least $1-2kd_1^{-2}$,  
\begin{align*}
\big\| e_{j}^{\tran}\what{E}\big(\mathfrak{P}^{\perp}\what{E}\mathfrak{P}^{\perp}\big)^{k-1}\what{E}\Theta\big\|\leq& \big\| e_{j_1}^{\tran}\what{E}\big(\mathfrak{P}^{\perp}\calP_{j_1}^{\perp}(\what{E})\big)^{k}\Theta\big\|+C_1\mumax\sqrt{\frac{r}{d_1}}\cdot\sum_{t=0}^{k}(C_2\delta)^{k-t}\delta^{t+1}\\
\leq&\big\| e_{j_1}^{\tran}\what{E}\big(\mathfrak{P}^{\perp}\calP_{j_1}^{\perp}(\what{E})\big)^{k}\Theta\big\|+C_1\mumax\sqrt{\frac{r}{d_1}}(C_2\delta)^{k+1}\cdot\sum_{t=0}^{k}C_2^{-(t+1)}\\
\leq&\big\| e_{j_1}^{\tran}\what{E}\big(\mathfrak{P}^{\perp}\calP_{j_1}^{\perp}(\what{E})\big)^{k}\Theta\big\|+\frac{C_1}{2}\mumax\sqrt{\frac{r}{d_1}}(C_2\delta)^{k+1}
\end{align*}
as long as $C_2>4$. 
We now bound $\big\| e_{j}^{\tran}\what{E}\big(\mathfrak{P}^{\perp}\calP_{j}^{\perp}(\what{E})\big)^{k}\Theta\big\|$. The idea is the same to {\it Case 1} and we shall utilize the independence between $e_{j}^{\tran}\what{E}$ and $\calP_{j}^{\perp}(\what{E})$, conditioned on $\calN_{j}$. Indeed, conditioned on $\calN_{j}$ and $\calP_{j}^{\perp}(\what{E})$, by Bernstein inequality, for all $t\geq 0$, 
\begin{align*}
\PP\bigg(\big\| e_{j}^{\tran}\what{E}\big(\mathfrak{P}^{\perp}\calP_{j}^{\perp}&(\what{E})\big)^{k}\Theta\big\|\geq C_1\sigma_{\xi}\sqrt{\frac{rd_1d_2 (t+\log d_1)}{n}}\|\calP_j^{\perp}\what E\|^k\\
&+C_2\sigma_{\xi}\frac{d_1d_2(t+\log d_1)}{n}\big\|\big(\mathfrak{P}^{\perp}\calP_{j}^{\perp}(\what{E})\big)^{k}\Theta\big\|_{\scriptscriptstyle \sf 2,max}
\bigg| \calN_{j}, \calP_{j}^{\perp}(\what{E})\bigg)\leq e^{-t}.
\end{align*}
Again, by the {\it Induction Assumption} and under event $\calE_0$, 
\begin{align*}
\big\|\big(\mathfrak{P}^{\perp}&\calP_{j}^{\perp}(\what{E})\big)^{k}\Theta \big\|_{\scriptscriptstyle \sf 2,max}\leq \big\|\big(\mathfrak{P}^{\perp}\calP_{j}^{\perp}(\what{E})\big)^{k-1}\mathfrak{P}^{\perp}\calP_{j}(\what{E})\Theta \big\|_{\scriptscriptstyle \sf 2,max}+\big\|\big(\mathfrak{P}^{\perp}\calP_{j}^{\perp}(\what{E})\big)^{k-1}\mathfrak{P}^{\perp}\what{E}\Theta \big\|_{\scriptscriptstyle \sf 2,max}\\
\leq& C_1\delta^{k}\cdot\mumax\sqrt{\frac{r}{d_2}}+\big\|\big(\mathfrak{P}^{\perp}\calP_{j}^{\perp}(\what{E})\big)^{k-1}\mathfrak{P}^{\perp}\what{E}\Theta \big\|_{\scriptscriptstyle \sf 2,max}\\
\leq &C_1\delta^{k}\cdot\mumax\sqrt{\frac{r}{d_2}}+\|(\mathfrak{P}^{\perp}\what{E})^{k}\Theta\|_{\scriptscriptstyle \sf 2,max}+\sum_{t=1}^{k-1}\Big\|\big(\mathfrak{P}^{\perp}\calP_{j}^{\perp}(\what{E})\big)^{k-t-1}\mathfrak{P}^{\perp}\calP_{j}(\what{E})\big(\mathfrak{P}^{\perp}\what{E}\big)^t\Theta\Big\|_{\scriptscriptstyle \sf 2,max}\\
&\hspace{2cm}\leq C_1\delta^k\mumax\sqrt{\frac{r}{d_2}}+C_1\mumax\sqrt{\frac{r}{d_2}}\cdot \sum_{t=0}^{k-1}(C_2\delta)^{t+1}\delta^{k-t-1}\leq 2C_1\mumax\sqrt{\frac{r}{d_2}}\cdot (C_2\delta)^{k}.
\end{align*}
By setting $t=3\log d_1$, conditioned on {\it Induction Assumption}, with probability at least $1-d_1^{-3}$, 
\begin{align*}
\big\|e_{j}^{\tran}\what{E}\big(\mathfrak{P}^{\perp}\calP_j^{\perp}(\what{E})\big)^k\Theta \big\|\leq& C_1\delta^k \cdot \sigma_{\xi}\sqrt{\frac{rd_1d_2\log d_1}{n}}+2C_1(C_2\delta)^k\cdot \mumax\sigma_{\xi}\frac{\sqrt{d_2d_1^2\log d_1^2}}{n}\\
\leq& 2C_1(C_2\delta)^{k}\cdot\sigma_{\xi}\sqrt{\frac{rd_1d_2\log d_1}{n}}\leq \frac{C_1}{2}\sqrt{\frac{r}{d_1}}(C_2\delta)^{k+1}
\end{align*}
where the last inequality holds as long as $n\geq C\mumax^2rd_1\log d_1$ for a large enough $C>0$.

Therefore, conditioned on {\it Induction Assumption}, with probability at least $1-d_1^{-2}$
$$
\max_{j\in[d_1]} \big\|e_j^{\tran}(\mathfrak{P}^{\perp}\what{E}\mathfrak{P}^{\perp})^{k}\what{E}\Theta \big\|\leq C_1(C_2\delta)^{k+1}\cdot \mumax\sqrt{\frac{r}{d_1}}.
$$
Finally, we conclude that, under event $\calE_0\cap \calE_1$, with probability at least $1-2(k+1)d_{\max}^{-2}$ so that for all $1\leq k_0\leq k+1$, 
\begin{equation}\label{eq:induct1}
\max_{j\in[d_1]}\big\|e_{j}^{\tran}\mathfrak{P}^{\perp}\big(\mathfrak{P}^{\perp}\what{E}\mathfrak{P}^{\perp}\big)^{k_0-1}\what{E}\Theta\big\|\leq C_1(C_2\delta)^{k_0}\cdot \mumax\sqrt{\frac{r}{d_1}}
\end{equation}
and
\begin{equation}
\max_{j\in[d_2]}\big\|e_{j+d_1}^{\tran}\mathfrak{P}^{\perp}\big(\mathfrak{P}^{\perp}\what{E}\mathfrak{P}^{\perp}\big)^{k_0-1}\what{E}\Theta\big\|\leq C_1(C_2\delta)^{k_0}\cdot \mumax\sqrt{\frac{r}{d_2}}
\end{equation}
where $C_1, C_2>0$ are some absolute constants. We conclude the proof of Lemma~\ref{lem:PE_2max}

\section{Proof of Lemma~\ref{lem:hatUZhatV}}
W.O.L.G., we only prove the upper bound for $\big|\big<\what U_1\what U_1^{\tran}\what{Z}^{(1)}\what V_1\what V_1^{\tran},T\big>\big|$. 
Clearly, 
\begin{align*}
\big|\big<\what U_1\what U_1^{\tran}\what{Z}^{(1)}\what V_1\what V_1^{\tran},T\big> \big|\leq \|T\|_{\ell_1}\cdot \big\|\what U_1\what U_1^{\tran}\what{Z}^{(1)}\what V_1\what V_1^{\tran} \big\|_{\submax}.
\end{align*}
It suffices to prove the upper bound for $\big\|\what U_1\what U_1^{\tran}\what{Z}^{(1)}\what V_1\what V_1^{\tran} \big\|_{\submax}$. By Theorem~\ref{lem:hatUhatV_deloc}, 
\begin{align*}
\big\|\what U_1\what U_1^{\tran}\what{Z}^{(1)}\what V_1\what V_1^{\tran} \big\|_{\submax}\leq& \big\| U U^{\tran}\what{Z}^{(1)} VV^{\tran} \big\|_{\submax}+\big\|(\what U_1\what U_1^{\tran}- U U^{\tran})\what{Z}^{(1)} VV^{\tran} \big\|_{\submax}\\
+&\big\| U U^{\tran}\what{Z}^{(1)} (\what V_1\what V_1^{\tran}-VV^{\tran}) \big\|_{\submax}+\big\|(\what U_1\what U_1^{\tran}- U U^{\tran})\what{Z}^{(1)} (\what V_1\what V_1^{\tran}-VV^{\tran}) \big\|_{\submax}\\
\leq&\big\| U U^{\tran}\what{Z}^{(1)} VV^{\tran} \big\|_{\submax}+\|\what{Z}^{(1)}\|\|\what U_1\what U_1^{\tran}- U U^{\tran}\|_{\submaxx}\|V\|_{\submaxx}\\
+\|\what{Z}^{(1)}\|\|\what V_1&\what V_1^{\tran}- V V^{\tran}\|_{\submaxx}\|U\|_{\submaxx}+\|\what{Z}^{(1)}\|\|\what V_1\what V_1^{\tran}- V V^{\tran}\|_{\submaxx}\|\what U_1\what U_1^{\tran}- U U^{\tran}\|_{\submaxx}\\
\leq&\big\| U U^{\tran}\what{Z}^{(1)} VV^{\tran} \big\|_{\submax}+C_2\mumax^2\frac{\sigma_{\xi}}{\lambda_r}\sqrt{\frac{rd_1^2d_2\log d_1}{n}}\cdot \sigma_{\xi}\sqrt{\frac{rd_1\log d_1}{n}}
\end{align*}
which holds under the event in Theorem~\ref{lem:hatUhatV_deloc}. Now, we prove the bound for $\|UU^{\tran}\what{Z}^{(1)}VV^{\tran}\|_{\submax}$. For any $j_1\in[d_1], j_2\in[d_2]$, we write
\begin{align*}
e_{j_1}^{\tran}UU^{\tran}\what{Z}^{(1)}VV^{\tran}e_{j_2}=&\frac{d_1d_2}{n_0}\sum_{i=n_0+1}^n \xi_i e_{j_1}^{\tran}UU^{\tran}X_iVV^{\tran}e_{j_2}\\
+&e_{j_1}^{\tran}UU^{\tran}\Big(\frac{d_1d_2}{n_0}\sum_{i=n_0+1}^n\langle\what \Delta_1,X_i \rangle X_i-\what\Delta_1\Big)VV^{\tran}e_{j_2}.
\end{align*}
Clearly, 
$$
\big\|  \xi_i e_{j_1}^{\tran}UU^{\tran}X_iVV^{\tran}e_{j_2}\big\|_{\psi_2}\leq \sigma_{\xi}\mumax^4\cdot \frac{r^2}{d_1d_2}
$$
and 
\begin{align*}
\EE \big(\xi_i e_{j_1}^{\tran}UU^{\tran}X_iVV^{\tran}e_{j_2}\big)^2=&\frac{\sigma_{\xi}^2}{d_1d_2}\sum_{i_1=1}^{d_1}\sum_{i_2=1}^{d_2}(e_{j_1}^{\tran}UU^{\tran}e_{i_1})^2(e_{i_2}^{\tran}VV^{\tran}e_{j_2})^2\\
=&\frac{\sigma_{\xi}^2}{d_1d_2}\|e_{j_1}^{\tran}U\|^2\|e_{j_2}^{\tran}V\|^2\leq \frac{\mumax^4\sigma_{\xi}^2r}{d_1^2d_2^2}.
\end{align*}
By Bernstein inequality, for all $t\geq 0$, with probability at least $1-e^{-t}$,
\begin{align*}
\Big|\frac{d_1d_2}{n_0}\sum_{i=n_0+1}^n \xi_i e_{j_1}^{\tran}UU^{\tran}X_iVV^{\tran}e_{j_2}\Big|\leq C_1\mumax^2\sigma_{\xi}\cdot\sqrt{\frac{r(t+\log d_1)}{n}}+C_2\mumax^4\sigma_{\xi}\cdot\frac{r^2(t+\log d_1)}{n}.
\end{align*}
By setting $t=3\log d_1$ and the union bound for all $j_1\in[d_1], j_2\in[d_2]$, we conclude that 
\begin{align*}
\PP\Big(\|UU^{\tran}\what{Z}_1^{(1)}VV^{\tran}\|_{\submax}\geq C_1\mumax^2\sigma_{\xi}\sqrt{\frac{r\log d_1}{n}}\Big)\leq d_1^{-2}
\end{align*}
as long as $n\geq C_3\mumax^4 r\log d_1$. Similar bounds also hold for $\|UU^{\tran}\what{Z}_2^{(1)}VV^{\tran}\|_{\submax}$. Therefore, conditioned on the event of Theorem~\ref{lem:hatUhatV_deloc}, with probability at least $1-2d_1^{-2}$, 
$$
\|\what U_1\what U_1^{\tran}\what{Z}^{(1)}\what V_1\what V_1^{\tran}\|_{\submax}\leq C_1\mumax^2\sigma_{\xi}\sqrt{\frac{r\log d_1}{n}}+C_2\mumax^2\frac{\sigma_{\xi}}{\lambda_r}\sqrt{\frac{rd_1^2d_2\log d_1}{n}}\cdot \sigma_{\xi}\sqrt{\frac{rd_1\log d_1}{n}}
$$
which concludes the proof of Lemma~\ref{lem:hatUZhatV}.

\section{Proof of Lemma~\ref{lem:normal_SA1}}
We aim to show the normal approximation of 
\begin{align*}
\frac{1}{2}\sum_{i=1}^2\big(\big<U_{\perp}U_{\perp}^{\tran}&\what{Z}^{(i)}VV^{\tran},T\big>+\big<UU\what{Z}^{(i)}V_{\perp}V_{\perp}^{\tran},T\big>\big)\\
=&\big<U_{\perp}U_{\perp}^{\tran}(\what{Z}^{(1)}/2+\what{Z}^{(2)}/2)VV^{\tran},T\big>+\big<UU(\what{Z}^{(1)}/2+\what{Z}^{(2)}/2)V_{\perp}V_{\perp}^{\tran},T\big>.
\end{align*}
 Recall that $\what{Z}^{(i)}=\what{Z}_1^{(i)}+\what{Z}_2^{(i)}$ where 
$$
\what{Z}_1^{(1)}=\frac{d_1d_2}{n_0}\sum_{i=n_0+1}^{n}\xi_i X_i\quad {\rm and}\quad \what{Z}_1^{(2)}=\frac{d_1d_2}{n_0}\sum_{i=1}^{n_0}\xi_i X_i
$$
so that (recall that $n=2n_0$)
$$
\frac{\what{Z}^{(1)}_1+\what{Z}^{(2)}_1}{2}=\frac{d_1d_2}{n}\sum_{i=1}^n \xi_i X_i
$$
and
$$
\frac{\what{Z}^{(1)}_2+\what{Z}^{(2)}_2}{2}=\frac{d_1d_2}{n}\Big(\sum_{i=1}^{n_0} \langle \what\Delta_2, X_i\rangle X_i+\sum_{i=n_0+1}^{n} \langle \what\Delta_1, X_i\rangle X_i\Big)-\frac{\what\Delta_1+\what\Delta_2}{2}.
$$
Therefore, write 
\begin{align*}
\big<U_{\perp}U_{\perp}^{\tran}(\what{Z}_1^{(1)}/2+\what{Z}^{(2)}_1/2)&VV^{\tran},T\big>+\big<UU(\what{Z}_1^{(1)}/2+\what{Z}^{(2)}_1/2)V_{\perp}V_{\perp}^{\tran},T\big>\\
&=\frac{d_1d_2}{n}\sum_{i=1}^n \xi_i\big(\big<U_{\perp}U_{\perp}^{\tran}X_iVV^{\tran},T\big>+\big<UU^{\tran}X_iV_{\perp}V_{\perp}^{\tran},T\big>\big)
\end{align*}
which is a sum of i.i.d. random variables: $\xi\big(\big<U_{\perp}U_{\perp}^{\tran}XVV^{\tran},T\big>+\big<UU^{\tran}XV_{\perp}V_{\perp}^{\tran},T\big>\big)$. To apply Berry-Essen theorem, we calculate its second and third moments. Clearly,
\begin{align*}
\EE\xi^2\big(\big<U_{\perp}U_{\perp}^{\tran}&XVV^{\tran},T\big>+\big<UUXV_{\perp}V_{\perp}^{\tran},T\big>\big)^2\\
=&\sigma_{\xi}^2\EE\big(\big<U_{\perp}U_{\perp}^{\tran}XVV^{\tran},T\big>+\big<UUXV_{\perp}V_{\perp}^{\tran},T\big>\big)^2\\
=&\sigma_{\xi}^2\big(\EE\big<U_{\perp}U_{\perp}^{\tran}XVV^{\tran},T\big>^2+\EE\big<UUXV_{\perp}V_{\perp}^{\tran},T\big>^2\\
&\hspace{2cm}+2\EE\big<U_{\perp}U_{\perp}^{\tran}XVV^{\tran},T\big>\big<UU^{\tran}XV_{\perp}V_{\perp}^{\tran},T\big>\big).
\end{align*}
Recall that $X$ is uniformly distributed over $\mathfrak{E}=\{e_{j_1}e_{j_2}^{\tran}: j_1\in[d_1], j_2\in[d_2]\}$. Therefore, 
\begin{align*}
\EE\big<U_{\perp}U_{\perp}^{\tran}XVV^{\tran},T\big>^2=&\frac{1}{d_1d_2}\sum_{i_1=1}^{d_1}\sum_{i_2=1}^{d_2}(e_{i_2}^{\tran}VV^{\tran}T^{\tran}U_{\perp}U_{\perp}^{\tran}e_{i_1})^2\\
=&\frac{1}{d_1d_2}\|VV^{\tran}T^{\tran}U_{\perp}U_{\perp}^{\tran}\|_{\subf}^2=\frac{1}{d_1d_2}\|V^{\tran}T^{\tran}U_{\perp}\|_{\subf}^2.
\end{align*}
Similarly,  $\EE\big<UUXV_{\perp}V_{\perp}^{\tran},T\big>^2=\|U^{\tran}TV_{\perp}\|_{\subf}^2/(d_1d_2)$. Meanwhile,
\begin{align*}
\EE\big<U_{\perp}U_{\perp}^{\tran}&XVV^{\tran},T\big>\big<UU^{\tran}XV_{\perp}V_{\perp}^{\tran},T\big>\\
=&\frac{1}{d_1d_2}\sum_{i_1=1}^{d_1}\sum_{i_2=1}^{d_2} (e_{i_2}^{\tran}VV^{\tran}T^{\tran}U_{\perp}U_{\perp}^{\tran}e_{i_1})\cdot (e_{i_1}^{\tran}UU^{\tran}TV_{\perp}V_{\perp}^{\tran}e_{i_2})\\
=&\frac{1}{d_1d_2}\sum_{i_2=1}^{d_2} e_{i_2}^{\tran}VV^{\tran}T^{\tran}U_{\perp}U_{\perp}^{\tran}UU^{\tran}TV_{\perp}V_{\perp}^{\tran}e_{i_2}\\
=&0
\end{align*}
where we used the fact $U^{\tran}U_{\perp}=0$. As a result, the second moment is
$$
\EE\xi^2\big(\big<U_{\perp}U_{\perp}^{\tran}XVV^{\tran},T\big>+\big<UUXV_{\perp}V_{\perp}^{\tran},T\big>\big)^2=\frac{\sigma_{\xi}^2}{d_1d_2}\big(\|V^{\tran}T^{\tran}U_{\perp}\|_{\subf}^2+\|U^{\tran}TV_{\perp}\|_{\subf}^2\big).
$$
Next, we bound the third moment of $\xi\big(\big<U_{\perp}U_{\perp}^{\tran}XVV^{\tran},T\big>+\big<UUXV_{\perp}V_{\perp}^{\tran},T\big>\big)$. By the sub-Gaussian Assumption~\ref{assump:noise}, we have
\begin{align*}
\EE|\xi|^3\big|\big<U_{\perp}&U_{\perp}^{\tran}XVV^{\tran},T\big>+\big<UU^{\tran}XV_{\perp}V_{\perp}^{\tran},T\big>\big|^3\\
\leq& C_2\sigma_{\xi}^3\cdot \EE\big|\big<U_{\perp}U_{\perp}^{\tran}XVV^{\tran},T\big>+\big<UU^{\tran}XV_{\perp}V_{\perp}^{\tran},T\big>\big|^3\\
=&C_2\sigma_{\xi}^3\cdot \frac{1}{d_1d_2}\sum_{i_1=1}^{d_1}\sum_{i_2=1}^{d_2}\big|\big<U_{\perp}U_{\perp}^{\tran}e_{i_1}e_{i_2}^{\tran}VV^{\tran},T\big>+\big<UU^{\tran}e_{i_1}e_{i_2}^{\tran}V_{\perp}V_{\perp}^{\tran},T\big>\big|^3\\
\leq&\frac{C_3\sigma_{\xi}^3}{d_1d_2}\cdot\sum_{i_1=1}^{d_1}\sum_{i_2=1}^{d_2}\big(\big|\big<U_{\perp}U_{\perp}^{\tran}e_{i_1}e_{i_2}^{\tran}VV^{\tran},T\big>\big|^3+\big|\big<UU^{\tran}e_{i_1}e_{i_2}^{\tran}V_{\perp}V_{\perp}^{\tran},T\big>\big|^3\big).
\end{align*}
Clearly, 
\begin{align*}
\big|\big<U_{\perp}U_{\perp}^{\tran}e_{i_1}e_{i_2}^{\tran}VV^{\tran},T\big>\big|=\big|\big<U_{\perp}^{\tran}e_{i_1}e_{i_2}^{\tran}V,U_{\perp}^{\tran}TV\big>\big|\leq \|U_{\perp}^{\tran}TV\|_{\subf}\big\| U_{\perp}^{\tran}e_{i_1}e_{i_2}^{\tran}V\big\|\\
\leq \|U_{\perp}^{\tran}TV\|_{\subf} \mumax\sqrt{\frac{r}{d_2}}
\leq \|U_{\perp}^{\tran}TV\|_{\subf} \cdot\mumax\sqrt{\frac{r}{d_2}}.
\end{align*}
Similar bound also holds for $|\langle UU^{\tran}e_{i_1}e_{i_2}^{\tran}V_{\perp}V_{\perp}^{\tran},T\rangle |$. Then, 
\begin{align*}
\EE|\xi|^3\big|&\big<U_{\perp}U_{\perp}^{\tran}XVV^{\tran},T\big>+\big<UU^{\tran}XV_{\perp}V_{\perp}^{\tran},T\big>\big|^3\\
\leq& C_3\frac{\sqrt{r}\mumax\sigma_{\xi}^3}{d_1d_2\sqrt{d_2}}\\
\times&\sum_{i_1=1}^{d_1}\sum_{i_2=1}^{d_2}\big(\big|\big<U_{\perp}U_{\perp}^{\tran}e_{i_1}e_{i_2}^{\tran}VV^{\tran},T\big>\big|^2\|U_{\perp}^{\tran}TV\|_{\subf}+\big|\big<UU^{\tran}e_{i_1}e_{i_2}^{\tran}V_{\perp}V_{\perp}^{\tran},T\big>\big|^2\|U^{\tran}TV_{\perp}\|_{\subf}\big).
\end{align*}
We write
\begin{align*}
\sum_{i_1=1}^{d_1}\sum_{i_2=1}^{d_2}\big(&\big|\big<U_{\perp}U_{\perp}^{\tran}e_{i_1}e_{i_2}^{\tran}VV^{\tran},T\big>\big|^2\\
=&\sum_{i_1=1}^{d_1}\sum_{i_2=1}^{d_2}(e_{i_1}^{\tran}U_{\perp}U_{\perp}^{\tran}TVV^{\tran}e_{i_2})^2=\|U_{\perp}U_{\perp}^{\tran}TVV^{\tran}\|_{\subf}^2=\|U_{\perp}^{\tran}TV\|_{\subf}^2
\end{align*}
Therefore,
\begin{align*}
\EE|\xi|^3\big|\big<U_{\perp}&U_{\perp}^{\tran}XVV^{\tran},T\big>+\big<UU^{\tran}XV_{\perp}V_{\perp}^{\tran},T\big>\big|^3\\
\leq&C_1\frac{\sigma_{\xi}^3\mumax \sqrt{r}}{d_1d_2\sqrt{d_2}}\big(\|U_{\perp}^{\tran}TV\|_{\subf}^3+\|U^{\tran}TV_{\perp}\|_{\subf}^3\big).
\end{align*}
By Berry-Essen theorem (\cite{berry1941accuracy}, \cite{esseen1956moment}), we get 
\begin{align}
\sup_{x\in\RR}\bigg|\PP\bigg(&\frac{\langle U_{\perp}U_{\perp}^{\tran}(\what{Z}^{(1)}_1/2+\what{Z}^{(2)}_1/2)VV^{\tran},T\rangle+\langle UU^{\tran}(\what{Z}^{(1)}_1/2+\what{Z}^{(2)}_1/2)V_{\perp}V_{\perp}^{\tran},T\rangle}{\sigma_{\xi}(\|V^{\tran}T^{\tran}U_{\perp}\|_{\subf}^2+\|U^{\tran}TV_{\perp}\|_{\subf}^2)^{1/2}\cdot\sqrt{d_1d_2/n}}\leq x\bigg)-\Phi(x)\bigg|\nonumber\\
\leq&C_4\mumax \sqrt{\frac{rd_1}{n}}\cdot \frac{\|U^{\tran}TV_{\perp}\|_{\subf}^3+\|U_{\perp}^{\tran}TV\|_{\subf}^3}{(\|V^{\tran}T^{\tran}U_{\perp}\|_{\subf}^2+\|U^{\tran}TV_{\perp}\|_{\subf}^2)^{3/2}}\leq C_4\mumax\sqrt{\frac{rd_1}{n}},\label{eq:berryessen1}
\end{align}
where $\Phi(x)$ denotes the c.d.f. of standard normal distributions. 
By  Assumption~\ref{assump:incoh}, we write
\begin{align}
\|U^{\tran}&TV\|_{\subf}\nonumber\\
&=\Big\|\sum_{j_1, j_2} T_{j_1,j_2}U^{\tran}e_{j_1}e_{j_2}^{\tran}V\Big\|_{\subf}
\leq&\sum_{j_1,j_2}|T_{j_1,j_2}|\cdot \|U^{\tran}e_{j_1}\|\|V^{\tran}e_{j_2}\|
\leq \|T\|_{\ell_1}\cdot\mumax^2\sqrt{\frac{r^2}{d_1d_2}}.\label{eq:VTUperp}
\end{align}
We then replace $\|V^{\tran}T^{\tran}U_{\perp}\|_{\subf}^2$ and $\|U^{\tran}TV_{\perp}\|_{\subf}^2$ with $\|TV\|_{\subf}^2$ and $\|U^{\tran}T\|_{\subf}^2$, respectively, to simplify the representation. We write
\begin{align*}
&\frac{\langle U_{\perp}U_{\perp}^{\tran}(\what{Z}^{(1)}_1/2+\what{Z}^{(2)}_1/2)VV^{\tran},T\rangle+\langle UU^{\tran}(\what{Z}^{(1)}_1/2+\what{Z}^{(2)}_1/2)V_{\perp}V_{\perp}^{\tran},T\rangle}{\sigma_{\xi}(\|V^{\tran}T^{\tran}U_{\perp}\|_{\subf}^2+\|U^{\tran}TV_{\perp}\|_{\subf}^2)^{1/2}\cdot\sqrt{d_1d_2/n}}\\
&=\frac{\langle U_{\perp}U_{\perp}^{\tran}(\what{Z}^{(1)}_1/2+\what{Z}^{(2)}_1/2)VV^{\tran},T\rangle+\langle UU^{\tran}(\what{Z}^{(1)}_1/2+\what{Z}^{(2)}_1/2)V_{\perp}V_{\perp}^{\tran},T\rangle}{\sigma_{\xi}(\|TV\|_{\subf}^2+\|U^{\tran}T\|_{\subf}^2)^{1/2}\cdot\sqrt{d_1d_2/n}}\\
&\hspace{2cm}+\frac{\langle U_{\perp}U_{\perp}^{\tran}(\what{Z}^{(1)}_1/2+\what{Z}^{(2)}_1/2)VV^{\tran},T\rangle+\langle UU^{\tran}(\what{Z}^{(1)}_1/2+\what{Z}^{(2)}_1/2)V_{\perp}V_{\perp}^{\tran},T\rangle}{\sigma_{\xi}(\|V^{\tran}T^{\tran}U_{\perp}\|_{\subf}^2+\|U^{\tran}TV_{\perp}\|_{\subf}^2)^{1/2}\cdot\sqrt{d_1d_2/n}}\\
&\hspace{5cm}\times\bigg(1-\frac{(\|V^{\tran}T^{\tran}U_{\perp}\|_{\subf}^2+\|U^{\tran}TV_{\perp}\|_{\subf}^2)^{1/2}}{(\|TV\|_{\subf}^2+\|U^{\tran}T\|_{\subf}^2)^{1/2}}\bigg).
\end{align*}
By Bernstein inequality, there exists an event $\calE_2$ with $\PP(\calE_2)\geq 1-d_1^{-2}$ so that under $\calE_2$, 
$$
\bigg|\frac{\langle U_{\perp}U_{\perp}^{\tran}(\what{Z}^{(1)}_1/2+\what{Z}^{(2)}_1/2)VV^{\tran},T\rangle+\langle UU^{\tran}(\what{Z}^{(1)}_1/2+\what{Z}^{(2)}_1/2)V_{\perp}V_{\perp}^{\tran},T\rangle}{\sigma_{\xi}(\|V^{\tran}T^{\tran}U_{\perp}\|_{\subf}^2+\|U^{\tran}TV_{\perp}\|_{\subf}^2)^{1/2}\cdot\sqrt{d_1d_2/n}} \bigg|\leq C_2\sqrt{\log d_1}
$$
for some large enough constant $C_2>0$. On the other hand, by Assumption~\ref{assump:T}, 
\begin{align*}
1-\frac{(\|V^{\tran}T^{\tran}U_{\perp}\|_{\subf}^2+\|U^{\tran}TV_{\perp}\|_{\subf}^2)^{1/2}}{(\|TV\|_{\subf}^2+\|U^{\tran}T\|_{\subf}^2)^{1/2}}&\leq 1-\frac{\|V^{\tran}T^{\tran}U_{\perp}\|_{\subf}^2+\|U^{\tran}TV_{\perp}\|_{\subf}^2}{\|TV\|_{\subf}^2+\|U^{\tran}T\|_{\subf}^2}\\
=&\frac{\|V^{\tran}T^{\tran}U\|_{\subf}^2+\|U^{\tran}TV\|_{\subf}^2}{\|TV\|_{\subf}^2+\|U^{\tran}T\|_{\subf}^2}\\
\leq&\frac{\mumax^4\|T\|_{\ell_1}^2}{\alpha_T^2\|T\|_{\subf}^2}\cdot\frac{r}{d_2}
\end{align*}
where the last inequality is due to (\ref{eq:VTUperp}).  Therefore, we conclude that, under event $\calE_2$, 
\begin{align*}
&\bigg|\frac{\langle U_{\perp}U_{\perp}^{\tran}(\what{Z}^{(1)}_1/2+\what{Z}^{(2)}_1/2)VV^{\tran},T\rangle+\langle UU^{\tran}(\what{Z}^{(1)}_1/2+\what{Z}^{(2)}_1/2)V_{\perp}V_{\perp}^{\tran},T\rangle}{\sigma_{\xi}(\|V^{\tran}T^{\tran}U_{\perp}\|_{\subf}^2+\|U^{\tran}TV_{\perp}\|_{\subf}^2)^{1/2}\cdot\sqrt{d_1d_2/n}}\\
&\hspace{2cm}- \frac{\langle U_{\perp}U_{\perp}^{\tran}(\what{Z}^{(1)}_1/2+\what{Z}^{(2)}_1/2)VV^{\tran},T\rangle+\langle UU^{\tran}(\what{Z}^{(1)}_1/2+\what{Z}^{(2)}_1/2)V_{\perp}V_{\perp}^{\tran},T\rangle}{\sigma_{\xi}(\|TV\|_{\subf}^2+\|U^{\tran}T\|_{\subf}^2)^{1/2}\cdot\sqrt{d_1d_2/n_0}}\bigg|\\
&\hspace{4cm}\leq C_2\frac{\mumax^4\|T\|_{\ell_1}^2}{\alpha_T^2\|T\|_{\subf}^2}\cdot\frac{r\sqrt{\log d_1}}{d_2}.
\end{align*}
By the Lipschitz property of $\Phi(x)$, it is obvious that (see, e.g., \cite{xia2018confidence, xia2019normal})
\begin{align*}
\sup_{x\in\RR}\bigg|\PP\bigg(&\frac{\langle U_{\perp}U_{\perp}^{\tran}(\what{Z}^{(1)}_1/2+\what{Z}^{(2)}_1/2)VV^{\tran},T\rangle+\langle UU^{\tran}(\what{Z}^{(1)}_1/2+\what{Z}^{(2)}_1/2)V_{\perp}V_{\perp}^{\tran},T\rangle}{\sigma_{\xi}(\|V^{\tran}T^{\tran}\|_{\subf}^2+\|U^{\tran}T\|_{\subf}^2)^{1/2}\cdot\sqrt{d_1d_2/n}}\leq x\bigg)-\Phi(x)\bigg|\\
\leq&C_4\mumax \sqrt{\frac{rd_1}{n}}+C_2\frac{\mumax^4\|T\|_{\ell_1}^2}{\alpha_T^2\|T\|_{\subf}^2}\cdot\frac{r\sqrt{\log d_1}}{d_2}+\frac{1}{d_1^2}.
\end{align*}
Next, we prove the upper bound for 
$$
\big<U_{\perp}U_{\perp}^{\tran}(\what{Z}_2^{(1)}/2+\what{Z}^{(2)}_2/2)VV^{\tran},T\big>+\big<UU^{\tran}(\what{Z}_2^{(1)}/2+\what{Z}^{(2)}_2/2)V_{\perp}V_{\perp}^{\tran},T\big>.
$$
We write
\begin{align*}
\big<U_{\perp}U_{\perp}^{\tran}&\what{Z}_2^{(1)}VV^{\tran},T\big>\\
=&\frac{d_1d_2}{n_0}\sum_{i=n_0+1}^{n}\langle\what\Delta_1,X_i \rangle\tr\big(T^{\tran}U_{\perp}U_{\perp}^{\tran}X_iVV^{\tran}\big)-\tr\big(T^{\tran}U_{\perp}U_{\perp}^{\tran}\what\Delta_1 VV^{\tran}\big).
\end{align*}
Observe that 
\begin{align*}
\big| \langle\what\Delta_1,X_i \rangle\tr\big(T^{\tran}U_{\perp}U_{\perp}^{\tran}X_iVV^{\tran}\big)&\big|\leq \|\what\Delta_1\|_{\submax}\|U_{\perp}^{\tran}TV\|_{\subf}\cdot \mumax\sqrt{\frac{r}{d_2}}.
\end{align*}
Moreover,
\begin{align*}
\EE\langle \what\Delta_1, X_i\rangle^2\Big(\tr\big(&T^{\tran}U_{\perp}U_{\perp}^{\tran}X_iVV^{\tran}\big)\Big)^2\leq \|\what\Delta_1\|_{\submax}^2\cdot \EE\Big(\tr\big(T^{\tran}U_{\perp}U_{\perp}^{\tran}X_iVV^{\tran}\big)\Big)^2\\
=&\frac{\|\what\Delta_1\|_{\submax}^2}{d_1d_2}\sum_{i_1=1}^{d_1}\sum_{i_2=1}^{d_2}e_{i_1}^{\tran}U_{\perp}U_{\perp}^{\tran}TVV^{\tran}e_{i_2}e_{i_2}^{\tran}VV^{\tran}T^{\tran}U_{\perp}U_{\perp}^{\tran}e_{i_1}\\
=&\frac{\|\what\Delta_1\|_{\submax}^2}{d_1d_2}\|U_{\perp}^{\tran}TV\|_{\subf}^2.
\end{align*}
By Bernstein inequality, with probability at least $1-d_1^{-2}$,
\begin{align*}
&\frac{\big|\big<U_{\perp}U_{\perp}^{\tran}\what{Z}_2^{(1)}VV^{\tran},T\big>+\big<UU^{\tran}\what{Z}_2^{(1)}V_{\perp}V_{\perp}^{\tran},T\big>\big|}{\|U_{\perp}^{\tran}TV\|_{\subf}+\|U^{\tran}TV_{\perp}\|_{\subf}}\\
&\hspace{3cm}\leq C_2\|\what\Delta_1\|_{\submax}\sqrt{\frac{d_1d_2\log d_1}{n}}+C_3\mumax\|\what\Delta_1\|_{\submax}\cdot\frac{\sqrt{rd_1^2d_2\log d_1}}{n}\\
&\hspace{3cm} \leq C_2\|\what\Delta_1\|_{\submax}\sqrt{\frac{d_1d_2\log d_1}{n}}
\end{align*}
where the last bound holds as long as $n\geq C\mumax^2rd_1\log d_1$ for a large enough constant $C>0$. Recall from Assumption~\ref{assump:init_entry} that 
$$
\PP\Big(\|\what\Delta_1\|_{\submax}^2\leq C_2\gamma_n^2\cdot\sigma_{\xi}^2\Big)\geq 1-d_1^{-2}.
$$
Therefore, with probability at least $1-2d_1^{-2}$, for $i=1,2$, 
$$
\frac{\big|\big<U_{\perp}U_{\perp}^{\tran}\what{Z}_2^{(i)}VV^{\tran},T\big>+\big<UU^{\tran}\what{Z}_2^{(i)}V_{\perp}V_{\perp}^{\tran},T\big>\big|}{\sigma_{\xi}(\|V^{\tran}T^{\tran}\|_{\subf}^2+\|U^{\tran}T\|_{\subf}^2)^{1/2}\cdot\sqrt{d_1d_2/n}}\leq C_3\gamma_n\sqrt{\log d_1}.
$$
By Lipschitz property of $\Phi(x)$, then 
\begin{align*}
\sup_{x\in\RR}\bigg|\PP\bigg(&\frac{\langle U_{\perp}U_{\perp}^{\tran}(\what{Z}^{(1)}/2+\what{Z}^{(2)}/2)VV^{\tran},T\rangle+\langle UU^{\tran}(\what{Z}^{(1)}/2+\what{Z}^{(2)}/2)V_{\perp}V_{\perp}^{\tran},T\rangle}{\sigma_{\xi}(\|V^{\tran}T^{\tran}\|_{\subf}^2+\|U^{\tran}T\|_{\subf}^2)^{1/2}\cdot\sqrt{d_1d_2/n}}\leq x\bigg)-\Phi(x)\bigg|\\
\leq&C_2\frac{\mumax^4\|T\|_{\ell_1}^2}{\alpha_T^2\|T\|_{\subf}^2}\cdot\frac{r\sqrt{\log d_1}}{d_2}+\frac{3}{d_1^2}+C_3\gamma_n\sqrt{\log d_1}+C_4\mumax\sqrt{\frac{rd_1}{n}}.
\end{align*}
We conclude the proof of Lemma~\ref{lem:normal_SA1}.

\section{Proof of Lemma~\ref{lem:sumSAk}}
The following fact is clear. 
\begin{align*}
\big|\sum_{i=1}^2\sum_{k=2}^{\infty}\big<&\big(\calS_{A,k}(\what{E}^{(i)})A\Theta\Theta^{\tran}+\Theta\Theta^{\tran}A\calS_{A,k}(\what{E}^{(i)})\big),\wtilde{T}\big>\big|\\
\leq&\sum_{i=1}^2\sum_{k=2}^{\infty}\big|\big<\big(\calS_{A,k}(\what{E}^{(i)})A\Theta\Theta^{\tran}+\Theta\Theta^{\tran}A\calS_{A,k}(\what{E}^{(i)})\big),\wtilde{T}\big>\big|\\
\leq&\|T\|_{\ell_1}\cdot\sum_{i=1}^2\sum_{k=2}^{\infty}\max_{\substack{j_1\in[d_1]\\j_2\in[d_2]}}\big|e_{j_1}^{\tran}\big(\calS_{A,k}(\what{E}^{(i)})A\Theta\Theta^{\tran}+\Theta\Theta^{\tran}A\calS_{A,k}(\what{E}^{(i)})\big)e_{d_1+j_2}\big|.
\end{align*}
Observe that for $i=1,2$
\begin{align*}
\big|e_{j_1}^{\tran}\big(\calS_{A,k}(\what{E}^{(i)})A\Theta\Theta^{\tran}\big)e_{j_2+d_1} \big|\leq \mumax\sqrt{\frac{r}{d_2}}\cdot \big\|e_{j_1}^{\tran}\calS_{A,k}(\what{E}^{(i)})A\Theta\big\|
\end{align*}
and
\begin{align*}
\big|e_{j_1}^{\tran}\big(\Theta\Theta^{\tran}A\calS_{A,k}(\what{E}^{(i)})\big)e_{d_1+j_2}\big|\leq \mumax\sqrt{\frac{r}{d_1}}\cdot \|e_{d_1+j_2}^{\tran}\calS_{A,k}(\what{E}^{(i)})A\Theta\|.
\end{align*}
Recall that 
$$
\calS_{A,k}(\what{E}^{(i)})=\sum_{\bs:s_1+\cdots+s_{k+1}=k} \mathfrak{P}^{-s_1}\what{E}^{(i)}\mathfrak{P}^{-s_2}\cdots\mathfrak{P}^{-s_k}\what{E}^{(i)}\mathfrak{P}^{-s_{k+1}},\quad \forall i=1,2.
$$
Then, we  write
$$
e_{j_1}^{\tran}\calS_{A,k}(\what{E}^{(i)})A\Theta=\sum_{\bs:s_1+\cdots+s_{k+1}=k} e_{j_1}^{\tran}\mathfrak{P}^{-s_1}\what{E}^{(i)}\mathfrak{P}^{-s_2}\cdots\mathfrak{P}^{-s_k}\what{E}^{(i)}\mathfrak{P}^{-s_{k+1}}A\Theta.
$$
Clearly, if $s_{k+1}=0$, then $\mathfrak{P}^{-s_{k+1}}A=\mathfrak{P}^{\perp}A=0$. Therefore, it suffices to focus on $s_{k+1}\geq 1$. Then, 
$$
e_{j_1}^{\tran}\calS_{A,k}(\what{E}^{(i)})A\Theta=\sum_{\substack{\bs:s_1+\cdots+s_{k+1}=k\\ s_{k+1}\geq 1}} e_{j_1}^{\tran}\mathfrak{P}^{-s_1}\what{E}^{(i)}\mathfrak{P}^{-s_2}\cdots\mathfrak{P}^{-s_k}\what{E}^{(i)}\Theta\Theta^{\tran}\mathfrak{P}^{-s_{k+1}}A\Theta.
$$
Let $\kmax=2\lceil\log d_1\rceil$. Then, for all $k\leq \kmax$, $i=1,2$ and by Lemma~\ref{lem:PE_2max} (and the arguments for the cases $s_1\geq 1$), 
\begin{align*}
\max_{j_1\in[d_1]}\big\|e_{j_1}^{\tran}\mathfrak{P}^{-s_1}\what{E}^{(i)}&\mathfrak{P}^{-s_2}\cdots\mathfrak{P}^{-s_k}\what{E}^{(i)}\Theta\Theta^{\tran}\mathfrak{P}^{-s_{k+1}}A\Theta\big\|\\
\leq&\max_{j_1\in[d_1]}\|e_{j_1}^{\tran}\mathfrak{P}^{-s_1}\what{E}^{(i)}\mathfrak{P}^{-s_2}\cdots\mathfrak{P}^{-s_k}\what{E}^{(i)}\Theta\|\cdot \|\mathfrak{P}^{-s_{k+1}}A\|\\
\leq&C_1\Big(\frac{C_2\delta}{\lambda_r}\Big)^{k-1}\delta\cdot\mumax\sqrt{\frac{r}{d_1}} 
\end{align*}
where $\delta$ is the upper bound of $\|\what{E}^{(i)}\|$ defined by (\ref{eq:delta_def}). Therefore, conditioned on event $\calE_0$ (see (\ref{eq:delta_def})) and the event of Lemma~\ref{lem:PE_2max} , for all $k\leq \kmax$, 
\begin{align*}
\max_{j_1\in[d_1]}\big\| e_{j_1}^{\tran}\calS_{A,k}(\what{E}^{(i)})A\Theta\big\|\leq C_1\Big(\frac{4C_2\delta}{\lambda_r}\Big)^{k-1}\delta\cdot\mumax\sqrt{\frac{r}{d_1}}\leq C_1\Big(\frac{4C_2\delta}{\lambda_r}\Big)^{k-1}\delta\cdot\mumax\sqrt{\frac{r}{d_1}}. 
\end{align*}
As a result, we get 
\begin{align*}
\max_{\substack{j_1\in[d_1]\\j_2\in[d_2]}}\sum_{i=1}^2\sum_{k=2}^{\kmax}\big|e_{j_1}^{\tran}\big(&\calS_{A,k}(\what{E}^{(i)})A\Theta\Theta^{\tran}+\Theta\Theta^{\tran}A\calS_{A,k}(\what{E}^{(i)})\big)e_{d_1+j_2}\big|\\
\leq&C_1\mumax^2\frac{r}{\sqrt{d_1d_2}}\delta \cdot\sum_{k=2}^{\kmax}\Big(\frac{4C_2\delta}{\lambda_r}\Big)^{k-1}\leq C_1\mumax^2\frac{r}{\sqrt{d_1d_2}}\delta \cdot\frac{\delta}{\lambda_r}
\end{align*}
where the last inequality holds since $8C_2\delta /\lambda_r<1$ by Assumption~\ref{assump:noise}. Moreover, on event $\calE_0$, we have
\begin{align*}
\max_{\substack{j_1\in[d_1]\\j_2\in[d_2]}}\sum_{i=1}^2\sum_{k=\kmax+1}^{\infty}\big|e_{j_1}^{\tran}\big(&\calS_{A,k}(\what{E}^{(i)})A\Theta\Theta^{\tran}+\Theta\Theta^{\tran}A\calS_{A,k}(\what{E}^{(i)})\big)e_{d_1+j_2}\big|\\
\leq&2\mumax\sqrt{\frac{r}{d_2}}\cdot \sum_{k=\kmax+1}^{\infty}\sum_{\bs: s_1+\cdots+s_{k+1}=k}\delta\cdot \Big(\frac{\delta}{\lambda_r}\Big)^{k-1}\\
\leq&2\mumax\sqrt{\frac{r}{d_2}}\cdot \sum_{k=\kmax+1}^{\infty}\delta\cdot \Big(\frac{4\delta}{\lambda_r}\Big)^{k-1}\leq 2\delta \mumax\sqrt{\frac{r}{d_2}}\cdot\Big(\frac{4\delta}{\lambda_r}\Big)^{\kmax}\\
\leq& 2\mumax\sqrt{\frac{r}{d_1^2d_2}}\delta\cdot\frac{\delta}{\lambda_r}
\end{align*}
where the last inequality is due to $(1/2)^{\log d_1}\leq d_1^{-1}$. Therefore, under the event of Theorem~\ref{lem:hatUhatV_deloc}, 
\begin{align*}
\big|\sum_{i=1}^2\sum_{k=2}^{\infty}\big<\big(\calS_{A,k}(\what{E}^{(i)})A\Theta\Theta^{\tran}+&\Theta\Theta^{\tran}A\calS_{A,k}(\what{E}^{(i)})\big),\wtilde{T}\big>\big|\\
\leq&C_2\|T\|_{\ell_1}\mumax^2\Big(\frac{\delta}{\lambda_r}\Big)\cdot \frac{r\delta}{\sqrt{d_1d_2}}
\end{align*}
which concludes the proof by replacing $\delta$ with $C\sigma_{\xi}\sqrt{d_1^2d_2\log d_1/n}$.

\section{Proof of Lemma~\ref{lem:hatAhat_err}}
By the definitions of $A$ and $\{\hat\Theta_i\}_{i=1}^2$, we have 
$$
\big<(\what\Theta_i\what\Theta_i^{\tran}-\Theta\Theta^{\tran})A(\what\Theta_i\what\Theta_i^{\tran}-\Theta\Theta^{\tran}),\wtilde{T}\big>=\big<(\what U_i\what U_i^{\tran}-UU^{\tran})M(\what V_i\what V_i^{\tran}-VV^{\tran}), T\big>
$$
for $i=1,2$. 
Then, 
\begin{align*}
\big|\big<(\what\Theta_i\what\Theta_i^{\tran}-\Theta&\Theta^{\tran})A(\what\Theta_i\what\Theta_i^{\tran}-\Theta\Theta^{\tran}),\wtilde{T}\big>\big|\\
\leq& \big\| (\what U_i\what U_i^{\tran}-UU^{\tran})M(\what V_i\what V_i^{\tran}-VV^{\tran})\big\|_{\submax}\cdot \|T\|_{\ell_1}\\
\leq& \|T\|_{\ell_1}\cdot \|\Lambda\| \|\what U_i\what U_i^{\tran}-UU^{\tran}\|_{\submaxx}\|\what V_i\what V_i^{\tran}-VV^{\tran}\|_{\submaxx}.
\end{align*}
Therefore, under the event of Theorem~\ref{lem:hatUhatV_deloc}, 
\begin{align*}
\big|\big<(\what\Theta_i\what\Theta_i^{\tran}-\Theta&\Theta^{\tran})A(\what\Theta_i\what\Theta_i^{\tran}-\Theta\Theta^{\tran}),\wtilde{T}\big>\big|\\
\leq&C_2\kappa_0\mumax^2\|T\|_{\ell_1}\sigma_{\xi}\sqrt{\frac{r^2d_1\log d_1}{n}}\cdot \frac{\sigma_{\xi}}{\lambda_r}\sqrt{\frac{d_1^2d_2\log d_1}{n}}.
\end{align*}

\section{Proof of Lemma~\ref{lem:hatG}}
By eq. (\ref{eq:hatG_t}), we write 
\begin{align*}
\frac{d_1d_2}{N_0}\sum_{j\in\mathfrak{D}_{2t}}\langle\what U^{(t)}\what G^{(t)}\what V^{(t)\tran}-U\Lambda V^{\tran}, X_j\rangle &\what U^{(t)\tran} X_j\what V^{(t)}
-\frac{d_1d_2}{N_0}\sum_{j\in\mathfrak{D}_{2t}}\xi_j\what U^{(t)\tran} X_j \what V^{(t)}=0
\end{align*}
where, due to data splitting, $(\what U^{(t)},\what V^{(t)})$ are independent with $\mathfrak{D}_{2t}$. 
Note that 
\begin{align*}
\what U^{(t)}\what G^{(t)}&\what V^{(t)\tran}-U\Lambda V^{\tran}\\
=&\what U^{(t)}(\what G^{(t)}-\what O_U^{(t)\tran}\Lambda\what{O}_V^{(t)})\what V^{(t)\tran}+\big(\what U^{(t)}\what O_U^{(t)\tran}\Lambda(\what V^{(t)}\what O_V^{(t)\tran})^{\tran}-U\Lambda V^{\tran}\big).
\end{align*}
Then, 
\begin{align*}
\what G^{(t)}-&\what O_U^{(t)\tran}\Lambda\what{O}_V^{(t)}\\
=&\big(\what G^{(t)}-\what O_U^{(t)\tran}\Lambda\what{O}_V^{(t)}\big)-\frac{d_1d_2}{N_0}\sum_{j\in\mathfrak{D}_{2t}}\langle\what G^{(t)}-\what O_U^{(t)\tran}\Lambda\what{O}_V^{(t)}, \what U^{(t)\tran} X_j\what V^{(t)}\rangle \what U^{(t)\tran} X_j\what V^{(t)}\\
&-\frac{d_1d_2}{N_0}\sum_{j\in\mathfrak{D}_{2t}}\langle \big(\what U^{(t)}\what O_U^{(t)\tran}\Lambda(\what V^{(t)}\what O_V^{(t)\tran})^{\tran}-U\Lambda V^{\tran}\big), X_j\rangle \what U^{(t)\tran} X_j\what V^{(t)}\\
&\hspace{2cm}+\frac{d_1d_2}{N_0}\sum_{j\in\mathfrak{D}_{2t}}\xi_j\what U^{(t)\tran} X_j \what V^{(t)}.
\end{align*}
Since $\|\what U^{(t)}\|\leq 2\mumax\sqrt{r/d_1},\|\what V^{(t)}\|\leq 2\mumax\sqrt{r/d_2}$, then 
$$
\Big\|\big\| \xi_i  \what U^{(t)\tran} X_i \what V^{(t)}\big\| \Big\|_{\psi_2}\lesssim \sigma_{\xi}\cdot \mumax^2\sqrt{\frac{r^2}{d_1d_2}}
$$
where the $\psi_2$-norm of a random variable $Z$ is defined by $\|Z\|_{\psi_2}=\min\{C>0: \exp(|Z|^2/C^2)\leq 2\}$. 
Meanwhile, 
\begin{align*}
\big\|\EE (\xi^2 \what U^{(t)\tran} X\what V^{(t)}\what V^{(t)\tran}X^{\tran}\what U^{(t)})\big\|=\sigma_{\xi}^2\big\| \EE ( \what U^{(t)\tran} X\what V^{(t)}\what V^{(t)\tran}X^{\tran}\what U^{(t)})\big\|=\sigma_{\xi}^2\cdot \frac{r}{d_1d_2}.
\end{align*}
By matrix Bernstein inequality \citep{tropp2012user,koltchinskii2011nuclear}, for any $t\geq 0$, 
\begin{align*}
\PP\Big(\Big\|\frac{d_1d_2}{N_0}\sum_{j\in\mathfrak{D}_{2t}}\xi_j  \what U^{(t)\tran} X_j \what V^{(t)} \Big\|\geq C_2\sigma_{\xi}&\sqrt{\frac{rd_1d_2(t+\log d_1)}{N_0}}+C_3\mumax^2\sigma_{\xi}\frac{\sqrt{r^2d_1d_2}(t+\log d_1)}{N_0}\Big)\leq e^{-t}. 
\end{align*}
By setting $t=2\log d_1$ and the fact $n\geq C_5\mumax^4r \log^2d_1$, we get, with probability at least $1-d_1^{-2}$, that
$$
\Big\|\frac{d_1d_2}{N_0}\sum_{j\in\mathfrak{D}_{2t}}\xi_j  \what U^{(t)\tran} X_j \what V^{(t)} \Big\|\leq C_2\sigma_{\xi}\sqrt{\frac{rd_1d_2\log d_1}{N_0}}
$$
for some absolute constant $C_2>0$. 

We then prove the upper bound for 
\begin{align*}
\Big\|\big(\what G^{(t)}-\what O_U^{(t)\tran}\Lambda&\what{O}_V^{(t)}\big)-\frac{d_1d_2}{N_0}\sum_{j\in\mathfrak{D}_{2t}}\langle\what G^{(t)}-\what O_U^{(t)\tran}\Lambda\what{O}_V^{(t)}, \what U^{(t)\tran} X_j\what V^{(t)}\rangle \what U^{(t)\tran} X_j\what V^{(t)} \Big\|
\end{align*}
where $\what G^{(t)}$ is dependent with $\{(X_j, Y_j)\}_{j\in\mathfrak{D}_{2t}}$. To this end, we write 
\begin{align*}
\Big\|\big(\what G^{(t)}-&\what O_U^{(t)\tran}\Lambda\what{O}_V^{(t)}\big)-\frac{d_1d_2}{N_0}\sum_{j\in\mathfrak{D}_{2t}}\langle\what G^{(t)}-\what O_U^{(t)\tran}\Lambda\what{O}_V^{(t)}, \what U^{(t)\tran} X_j\what V^{(t)}\rangle \what U^{(t)\tran} X_j\what V^{(t)} \Big\|\\
\leq&\|\what G^{(t)}-\what O_U^{(t)\tran}\Lambda\what{O}_V^{(t)}\|
\cdot \sup_{A\in \RR^{r\times r}, \|A\|\leq 1}\Big\|A-\frac{d_1d_2}{N_0}\sum_{j\in\mathfrak{D}_{2t}}\langle A, \what U^{(t)\tran} X_j\what V^{(t)}\rangle \what U^{(t)\tran} X_j\what V^{(t)} \Big\|.
\end{align*}
Denote $\calO_r=\{A\in\RR^{r\times r}, \|A\|\leq 1\}$ and  $\calN_{1/3}(\calO_r)$ the $1/3$-net of $\calO_r$, i.e., for any $A\in\calO_r$, there exists $A_0\in \calN_{1/3}(\calO_r)$ so that $\|A-A_0\|\leq 1/3$. It is well-known by \citep{pajor1998metric, koltchinskii2015optimal} that ${\rm Card}\big(\calN_{1/3}(\calO_r)\big)\leq 3^{C_2r^2}$ for some absolute constants $C_2>0$. By the definition of $\calN_{1/3}(\calO_r)$, 
\begin{align*}
\Big\|\big(\what G^{(t)}-&\what O_U^{(t)\tran}\Lambda\what{O}_V^{(t)}\big)-\frac{d_1d_2}{N_0}\sum_{j\in\mathfrak{D}_{2t}}\langle\what G^{(t)}-\what O_U^{(t)\tran}\Lambda\what{O}_V^{(t)}, \what U^{(t)\tran} X_j\what V^{(t)}\rangle \what U^{(t)\tran} X_j\what V^{(t)} \Big\|\\
\leq&3\|\what G^{(t)}-\what O_U^{(t)\tran}\Lambda\what{O}_V^{(t)}\|
\cdot \max_{A\in \calN_{1/3}(\calO_r)}\Big\|A-\frac{d_1d_2}{N_0}\sum_{j\in\mathfrak{D}_{2t}}\langle A, \what U^{(t)\tran} X_j\what V^{(t)}\rangle \what U^{(t)\tran} X_j\what V^{(t)} \Big\|.
\end{align*}
For each $A\in\calN_{1/3}(\calO_r)$, 
\begin{align*}
\big\| \langle A, \what U^{(t)\tran} X_i&\what V^{(t)}\rangle \what U^{(t)\tran} X_i\what V^{(t)} \big\|
\leq \|\what U^{(t)\tran} X_i\what V^{(t)}\|_{\star}\cdot \| \what U^{(t)\tran} X_i\what V^{(t)} \|
\leq  \|\what U^{(t)}\|_{\submaxx}^2\|\what V^{(t)}\|_{\submaxx}^2\leq \mumax^4\frac{r^2}{d_1d_2}
\end{align*}
where $\|\cdot\|_{\star}$ denotes the matrix nuclear norm. Moreover,
\begin{align*}
\big\|\EE  \langle A, \what U^{(t)\tran}& X_i\what V^{(t)}\rangle^2 \what U^{(t)\tran} X_i\what V^{(t)} \what V^{(t)\tran}X_i^{\tran}\what U^{(t)}\big\|
\leq\mumax^4\frac{r^3}{(d_1d_2)^2}.
\end{align*}
Therefore, for each $A\in\calN_{1/3}(\calO_r)$ and any $t>0$, 
\begin{align*}
\PP\bigg(\Big\|A-\frac{d_1d_2}{N_0}\sum_{j\in\mathfrak{D}_{2t}}\langle A, \what U^{(t)\tran} X_j\what V^{(t)}\rangle \what U^{(t)\tran} X_j\what V^{(t)} \Big\|\geq &C_1\mumax^2\sqrt{\frac{r^3(t+\log d_1)}{N_0}}\\
+&C_2\mumax^4\frac{r^2(t+\log d_1)}{N_0}\bigg)\leq e^{-t}. 
\end{align*}
By setting $t=C_2r^2+2\log d_1$ and the union bound over all $A\in\calN_{1/3}(\calO_r)$, if $n\geq C_3\mumax^4 (r^3+r\log d_1)\log d_1$, then with probability at least $1-d_1^{-2}$, 
$$
 \max_{A\in \calN_{1/3}(\calO_r)}\Big\|A-\frac{d_1d_2}{N_0}\sum_{j\in\mathfrak{D}_{2t}}\langle A, \what U^{(t)\tran} X_j\what V^{(t)}\rangle \what U^{(t)\tran} X_j\what V^{(t)} \Big\|\leq C_1\mumax^2\sqrt{\frac{r^3(r^2+\log d_1)}{N_0}}
$$
implying that 
\begin{align*}
\Big\|\big(\what G^{(t)}-&\what O_U^{(t)\tran}\Lambda\what{O}_V^{(t)}\big)-\frac{d_1d_2}{N_0}\sum_{j\in\mathfrak{D}_{2t}}\langle\what G^{(t)}-\what O_U^{(t)\tran}\Lambda\what{O}_V^{(t)}, \what U^{(t)\tran} X_j\what V^{(t)}\rangle \what U^{(t)\tran} X_j\what V^{(t)} \Big\|\\
\leq&\|\what G^{(t)}-\what O_U^{(t)\tran}\Lambda\what{O}_V^{(t)}\|\cdot C_1\mumax^2\sqrt{\frac{r^3(r^2+\log d_1)}{N_0}}.
\end{align*}
Similarly if $n\geq C_2\mumax^4 r\log^2d_1$, then with probability at least $1-d_1^{-2}$, 
\begin{align*}
\Big\|\frac{d_1d_2}{N_0}\sum_{j\in\mathfrak{D}_{2t}}\langle \big(\what U^{(t)}&\what O_U^{(t)\tran}\Lambda(\what V^{(t)}\what O_V^{(t)\tran})^{\tran}-U\Lambda V^{\tran}\big), X_j\rangle \what U^{(t)\tran} X_j\what V^{(t)} \Big\|\\
\leq&\big\| \what U^{(t)\tran}\big(\what U^{(t)}\what O_U^{(t)\tran}\Lambda(\what V^{(t)}\what O_V^{(t)\tran})^{\tran}-U\Lambda V^{\tran}\big)\what V^{(t)}\big\|\\
&+C_2\big\|\what U^{(t)}\what O_U^{(t)\tran}\Lambda(\what V^{(t)}\what O_V^{(t)\tran})^{\tran}-U\Lambda V^{\tran}\big\|_{\submax}\cdot \sqrt{\frac{rd_1d_2\log d_1}{N_0}}
\end{align*}
where we used the fact 
\begin{align*}
\big\|\langle \big(\what U^{(t)}&\what O_U^{(t)\tran}\Lambda(\what V^{(t)}\what O_V^{(t)\tran})^{\tran}-U\Lambda V^{\tran}\big), X_j\rangle \what U^{(t)\tran} X_j\what V^{(t)} \big\|\\
\leq& \big\|\what U^{(t)}\what O_U^{(t)\tran}\Lambda(\what V^{(t)}\what O_V^{(t)\tran})^{\tran}-U\Lambda V^{\tran}\big\|_{\submax}\cdot \mumax^2\sqrt{\frac{r^2}{d_1d_2}}
\end{align*}
and
\begin{align*}
\big\|\EE\langle \what U^{(t)}\what O_U^{(t)\tran}&\Lambda(\what V^{(t)}\what O_V^{(t)\tran})^{\tran}-U\Lambda V^{\tran}, X_j\rangle ^2\what U^{(t)\tran} X_j\what V^{(t)} \what V^{(t)\tran}X_j^{\tran}\what U^{(t)} \big\|\\
\leq&  \big\|\what U^{(t)}\what O_U^{(t)\tran}\Lambda(\what V^{(t)}\what O_V^{(t)\tran})^{\tran}-U\Lambda V^{\tran}\big\|_{\submax}^2\cdot \frac{r}{d_1d_2}.
\end{align*}
Therefore, we conclude that if $n\geq C_2\mumax^4 r^3(r^2+\log d_1) \log d_1$, then with probability at least $1-3d_1^{-2}$, 
\begin{align*}
\|\what G^{(t)}-\what O_U^{(t)\tran}&\Lambda\what O_V^{(t)}\|\leq\big\|\what U^{(t)\tran}\big(\what U^{(t)}\what O_U^{(t)\tran}\Lambda(\what V^{(t)}\what O_V^{(t)\tran})^{\tran}-U\Lambda V^{\tran}\big)\what V^{(t)} \big\|
+C_6\sigma_{\xi}\sqrt{\frac{rd_1d_2\log d_1}{N_0}}\\
&+C_2\|\what U^{(t)}\what O_U^{(t)\tran}\Lambda(\what V^{(t)}\what O_V^{(t)\tran})^{\tran}-U\Lambda V^{\tran}\|_{\submax}\cdot \sqrt{\frac{rd_1d_2\log d_1}{N_0}}
\end{align*}
Note that 
\begin{align*}
\|\what U^{(t)}\what O_U^{(t)\tran}\Lambda(\what V^{(t)}&\what O_V^{(t)\tran})^{\tran}-U\Lambda V^{\tran}\|_{\submax}\\
&\leq 3\|\Lambda\|\mumax\cdot\Big(\sqrt{\frac{r}{d_2}}\|\what U^{(t)}-U\what O_U^{(t)}\|_{\submaxx}+\sqrt{\frac{r}{d_1}}\|\what V^{(t)}-V\what O_V^{(t)}\|_{\submaxx}\big).
\end{align*}
By the differential property of Grassmannians, see, e.g., \citep{keshavan2010matrix_a, xia2017polynomial, edelman1998geometry}, 
\begin{align*}
\big\|\what U^{(t)\tran}\big(\what U^{(t)}\what O_U^{(t)\tran}&\Lambda(\what V^{(t)}\what O_V^{(t)\tran})^{\tran}-U\Lambda V^{\tran}\big)\what V^{(t)} \big\|\\
\leq& \|\Lambda\|\cdot \|\what U^{(t)\tran}(\what U^{(t)}\what O_U^{(t)\tran}-U)\|+\|\Lambda\|\cdot \|\what V^{(t)\tran}(\what V^{(t)}\what O_V^{(t)\tran}-V)\|\\
\leq& 2\|\Lambda\|\cdot \|\what U^{(t)}-U\what O_U^{(t)}\|^2+2\|\Lambda\|\cdot \|\what V^{(t)}-V\what O_V^{(t)}\|^2.
\end{align*}
Finally, we conclude with probability at least $1-3d_1^{-2}$,
\begin{align*}
\big\|\what G^{(t)}-&\what O_U^{(t)\tran}\Lambda\what O_V^{(t)} \big\|\leq C_5\sigma_{\xi}\sqrt{\frac{rd_1d_2\log d_1}{N_0}}+2\|\Lambda\|\cdot \big(\|\what U^{(t)}-U\what O_U^{(t)}\|^2+\|\what V^{(t)}-V\what O_V^{(t)}\|^2\big)\\
+&C_7\|\Lambda\|\Big(\sqrt{\frac{r}{d_2}}\|\what U^{(t)}-U\what O_U^{(t)}\|_{\submaxx}+\sqrt{\frac{r}{d_1}}\|\what V^{(t)}-V\what O_V^{(t)}\|_{\submaxx}\big)\cdot\mumax\sqrt{\frac{rd_1d_2\log d_1}{N_0}}
\end{align*}

\section{Proof of Lemma~\ref{lem:gd_infty}}
Recall that 
\begin{align}
\what U^{(t+0.5)}-&U\what O_U^{(t)}\what L_G^{(t)}=(\what U^{(t)}\what L_G^{(t)}-U\what O^{(t)}_U\what L_G^{(t)})\big(I-\eta\cdot \what L_G^{(t)\tran}\what O_U^{(t)\tran}\Lambda \what O_V^{(t)}\what R_G^{(t)}(\what \Lambda^{(t)})^{-1}\big)\nonumber\\
-&\eta\cdot \what U^{(t)}\big(\what G^{(t)}-\what O_U^{(t)\tran}\Lambda\what O_V^{(t)}\big)\what R^{(t)}_G(\what \Lambda^{(t)})^{-1}
-\eta\cdot U\Lambda(\what V^{(t)}\what O_V^{(t)\tran}-V)^{\tran}\what V^{(t)}\what R_G^{(t)}(\what \Lambda^{(t)})^{-1}\nonumber\\
&\hspace{2cm}+\what E_V^{(t)}+\what E_{\xi,V}^{(t)}.\label{eq:hatUt0.5_3}
\end{align}
By Lemma~\ref{lem:hatG},
\begin{align}
\big\|\what\Lambda^{(t)}-&\what L_G^{(t)\tran}\what O_U^{(t)\tran}\Lambda\what O_V^{(t)} \what R_G^{(t)}\big\|\leq C_5\sigma_{\xi}\sqrt{\frac{rd_1d_2\log d_1}{N_0}}
+2\|\Lambda\|\cdot \big(\|\what U^{(t)}-U\what O_U^{(t)}\|^2+\|\what V^{(t)}-V\what O_V^{(t)}\|^2\big)\nonumber\\
+&C_7\|\Lambda\|\Big(\sqrt{\frac{r}{d_2}}\|\what U^{(t)}-U\what O_U^{(t)}\|_{\submaxx}+\sqrt{\frac{r}{d_1}}\|\what V^{(t)}-V\what O_V^{(t)}\|_{\submaxx}\big)\cdot\mumax\sqrt{\frac{rd_1d_2\log d_1}{N_0}}\label{eq:hatSigma}
\end{align}
which implies that $\|\what\Lambda^{(t)}-\what L_G^{(t)\tran}\what O_U^{(t)\tran}\Lambda\what O_V^{(t)} \what R_G^{(t)}\|\leq \lambda_r/20$ under Assumption~\ref{assump:noise} and when $\max\{\|\what U^{(t)}-U\what O_U^{(t)}\|, \|\what V^{(t)}-V\what O_V^{(t)}\|\}\leq 1/(80\sqrt{\kappa_0})$, $\|\what U^{(t)}\|_{\submaxx}\leq 2\mumax\sqrt{r/d_1}, \|\what V^{(t)}\|_{\submaxx}\leq 2\mumax\sqrt{r/d_2}$ and $n\geq C_3\kappa_0^2 \mumax^4 r^3\log d_1$.

Since $\eta\leq 0.75$, we have 
\begin{align*}
\big\|&(\what U^{(t)}\what L_G^{(t)}-U\what O^{(t)}_U\what L_G^{(t)})\big(I-\eta\cdot \what L_G^{(t)\tran}\what O_U^{(t)\tran}\Lambda \what O_V^{(t)}\what R_G^{(t)}(\what \Lambda^{(t)})^{-1}\big) \big\|_{\submaxx}\\
\leq&\big\|(\what U^{(t)}\what L_G^{(t)}-U\what O^{(t)}_U\what L_G^{(t)})\big\|_{\submaxx}\cdot \big\| I-\eta\cdot \what L_G^{(t)\tran}\what O_U^{(t)\tran}\Lambda \what O_V^{(t)}\what R_G^{(t)}(\what \Lambda^{(t)})^{-1}\big\|\\
\leq&\|\what U^{(t)}-U\what O_U^{(t)}\|_{\submaxx}\cdot  (1-\eta)+\|\what U^{(t)}-U\what O_U^{(t)}\|_{\submaxx}\cdot \eta\big\|\what \Lambda^{(t)}- \what L_G^{(t)\tran}\what O_U^{(t)\tran}\Lambda \what O_V^{(t)}\what R_G^{(t)}\big\|\cdot \big\|(\what \Lambda^{(t)})^{-1}\big\|\\
\leq& \big(1-\eta\big)\cdot \|\what U^{(t)}-U\what O_U^{(t)}\|_{\submaxx}+2\|\what U^{(t)}-U\what O_U^{(t)}\|_{\submaxx}\cdot \eta\big\|\what \Lambda^{(t)}- \what L_G^{(t)\tran}\what O_U^{(t)\tran}\Lambda \what O_V^{(t)}\what R_G^{(t)}\big\|\cdot \lambda_r^{-1}
\end{align*}
where the last inequality is due to $\lambda_r(\what \Lambda^{(t)})\geq \lambda_r/2$ by (\ref{eq:hatSigma}). Again, by Lemma~\ref{lem:hatG} and Assumption~\ref{assump:noise},
\begin{align*}
2 \big\|\what \Lambda^{(t)}- \what L_G^{(t)\tran}\what O_U^{(t)\tran}\Lambda \what O_V^{(t)}\what R_G^{(t)}\big\|\cdot \lambda_r^{-1}\leq \frac{1}{10}.
\end{align*}
Then, we obtain
\begin{align*}
\big\|(\what U^{(t)}\what L_G^{(t)}-U\what O^{(t)}_U\what L_G^{(t)})\big(I-\eta\cdot \what L_G^{(t)\tran}&\what O_U^{(t)\tran}\Lambda \what O_V^{(t)}\what R_G^{(t)}(\what \Lambda^{(t)})^{-1}\big) \big\|_{\submaxx}
\leq \Big(1-\frac{9\eta}{10}\Big)\cdot  \|\what U^{(t)}-U\what O_U^{(t)}\|_{\submaxx}.
\end{align*}
Since $\|\what U^{(t)}\|_{\submax}\leq 2\mumax\sqrt{r/d_1}$, by (\ref{eq:hatSigma}), we get 
\begin{align*}
\eta\big\|  \what U^{(t)}\big(&\what G^{(t)}-\what O_U^{(t)\tran}\Lambda\what O_V^{(t)}\big)\what R^{(t)}_G(\what \Lambda^{(t)})^{-1}\big\|_{\submaxx}
\leq 2\eta\|\what U^{(t)}\|_{\submaxx}\cdot \|\what G^{(t)}-\what O_U^{(t)\tran}\Lambda\what O_V^{(t)}\|\cdot \lambda_r^{-1}\\
\leq& C_3\frac{\mumax\eta\sigma_{\xi}}{\lambda_r}\cdot \sqrt{\frac{r^2d_2\log d_1}{N_0}}+C_4\eta\kappa_0\mumax\sqrt{\frac{r}{d_1}}\Big(\|\what U^{(t)}-U\what O_U^{(t)}\|^2+\|\what V^{(t)}-V\what O_V^{(t)}\|^2\Big)\\
+C_5&\eta\kappa_0\mumax\cdot \Big(\sqrt{\frac{r}{d_2}}\|\what U^{(t)}-U\what O_U^{(t)}\|_{\submaxx}+\sqrt{\frac{r}{d_1}}\|\what V^{(t)}-V\what O_V^{(t)}\|_{\submaxx}\big)\cdot\mumax\sqrt{\frac{r^2d_2\log d_1}{N_0}}.
\end{align*}
Observe that 
\begin{align*}
C_4\kappa_0\mumax\sqrt{\frac{r}{d_1}} \|\what U^{(t)}-U\what O_U^{(t)}\|^2\leq C_4\kappa_0\mumax \sqrt{r} \|\what U^{(t)}-U\what O_U^{(t)}\|\cdot \frac{ \|\what U^{(t)}-U\what O_U^{(t)}\|}{\sqrt{d_1}}\\
\leq C_4\kappa_0\mumax \sqrt{r} \|\what U^{(t)}-U\what O_U^{(t)}\|\cdot \frac{ \|\what U^{(t)}-U\what O_U^{(t)}\|_{\subf}}{\sqrt{d_1}}\leq \frac{\|\what U^{(t)}-U\what O_U^{(t)}\|_{\submaxx}}{20}
\end{align*}
if $\|\what U^{(t)}-U\what O_U^{(t)}\|\leq 1/(20\kappa_0\mumax\sqrt{r})$. Similarly, if $\|\what V^{(t)}-V\what O_V^{(t)}\|\leq 1/(20\kappa_0\mumax\sqrt{r\alpha_d})$ where $\alpha_d=d_1/d_2$, then
\begin{align*}
C_4\kappa_0\mumax\sqrt{\frac{r}{d_1}} \|\what V^{(t)}-V\what O_V^{(t)}\|^2\leq  \frac{\|\what V^{(t)}-V\what O_V^{(t)}\|_{\submaxx}}{20}.
\end{align*}
Therefore,
\begin{align*}
C_4\kappa_0\mumax\sqrt{\frac{r}{d_1}}\Big(\|\what U^{(t)}-U&\what O_U^{(t)}\|^2+\|\what V^{(t)}-V\what O_V^{(t)}\|^2\Big)\\
&\leq \frac{\|\what U^{(t)}-U\what O_U^{(t)}\|_{\submaxx}+\|\what V^{(t)}-V\what O_V^{(t)}\|_{\submaxx}}{20}.
\end{align*}
Moreover, if $n\geq C_1\alpha_d\kappa_0^2\mumax^4 r^3\log^2 d_1$, then 
\begin{align*}
C_5\kappa_0\mumax\cdot \Big(\sqrt{\frac{r}{d_2}}\|\what U^{(t)}-&U\what O_U^{(t)}\|_{\submaxx}+\sqrt{\frac{r}{d_1}}\|\what V^{(t)}-V\what O_V^{(t)}\|_{\submaxx}\big)\cdot\mumax\sqrt{\frac{r^2d_2\log d_1}{N_0}}\\
\leq&  \frac{\|\what U^{(t)}-U\what O_U^{(t)}\|_{\submaxx}+\|\what V^{(t)}-V\what O_V^{(t)}\|_{\submaxx}}{40}.
\end{align*}
Then, we get
\begin{align*}
\eta\big\|  \what U^{(t)}&\big(\what G^{(t)}-\what O_U^{(t)\tran}\Lambda\what O_V^{(t)}\big)\what R^{(t)}_G(\what \Lambda^{(t)})^{-1}\big\|_{\submaxx}\\
\leq& C_3\frac{\eta\mumax\sigma_{\xi}}{\lambda_r}\cdot \sqrt{\frac{r^2d_2\log d_1}{N_0}}+\frac{3\eta}{40}\cdot \big(\|\what U^{(t)}-U\what O_U^{(t)}\|_{\submaxx}+\|\what V^{(t)}-V\what O_V^{(t)}\|_{\submaxx}\big).
\end{align*}
Since $\|(\what V^{(t)}\what O_V^{(t)\tran}-V)^{\tran}\what V^{(t)}\|\leq \|\what V^{(t)}\what O_V^{(t)\tran}-V\|^2$, we get 
\begin{align*}
\big\|U\Lambda(\what V^{(t)}\what O_V^{(t)\tran}-V)^{\tran}\what V^{(t)}\what R_G^{(t)}(\what \Lambda^{(t)})^{-1}\big\|_{\submaxx}\leq 2\kappa_0\mumax\sqrt{\frac{r}{d_1}}\|\what V^{(t)}\what O_V^{(t)\tran}-V\|^2 
\leq \frac{1}{40}\cdot \|\what V^{(t)}-V\what O_V^{(t)}\|_{\submaxx}
\end{align*}
if $\|\what V^{(t)}\what O_V^{(t)\tran}-V\|\leq 1/(C_2\kappa_0\mumax\sqrt{r})$.  Putting together the above bounds, we obtain 
\begin{align*}
\big\|&\what U^{(t+0.5)}-U\what O_U^{(t)}\what L_G^{(t)} \big\|_{\submaxx}
\leq \Big(1-\frac{9\eta}{10}\Big)\|\what U^{(t)}-U\what O_U^{(t)}\|_{\submaxx}+C_3\eta\mumax\frac{\sigma_{\xi}}{\lambda_r}\sqrt{\frac{r^2d_2\log d_1}{N_0}}\\
&+\frac{\eta}{10}\cdot \big(\|\what U^{(t)}-U\what O_U^{(t)}\|_{\submaxx}+\|\what V^{(t)}-V\what O_V^{(t)}\|_{\submaxx}\big)+\|\what E_V\|_{\submaxx}+\|\what E_{\xi,V}\|_{\submaxx}
\end{align*}
Since $(\what U^{(t)},\what G^{(t)}, \what V^{(t)})$ are independent with $\mathfrak{D}_{2t+1}$ and $\what U^{(t)},\what V^{(t)}$ are incoherent, by Bernstein inequality and an union bound for all rows, with probability at least $1-d_1^{-2}$,
\begin{align*}
\|\what E_V\|_{\submaxx}+\|\what E_{\xi,V}\|_{\submaxx}\leq C_3\eta\frac{\sigma_{\xi}+\|\what M^{(t)}-M\|_{\submax}}{\lambda_r}\sqrt{\frac{rd_1d_2\log d_1}{N_0}}
\end{align*}
where $\what M^{(t)}=\what U^{(t)}\what G^{(t)}\what V^{(t)\tran}$. Note that 
\begin{align*}
\big\|\what M^{(t)}-&M\big\|_{\submax}\leq \|(\what U^{(t)}-U\what O_U^{(t)})\what G^{(t)}\what V^{(t)\tran}\|_{\submax}+\|U\big(\what O_U^{(t)}\what G^{(t)}-\Lambda\what O_V^{(t)}\big)\what V^{(t)\tran}\|_{\submax}\\
&\hspace{3cm}+\big\|U\Lambda\big(\what V^{(t)}\what O_V^{(t)\tran}-V\big)^{\tran} \big\|_{\submaxx}\\
\leq& 2\mumax\|\Lambda\|\cdot\Big(\sqrt{\frac{r}{d_2}}\|\what U^{(t)}-U\what O_U^{(t)}\|_{\submaxx}+\sqrt{\frac{r}{d_1}}\|\what V^{(t)}-V\what O_V^{(t)}\|_{\submaxx}\Big)\\
&\hspace{3cm}+\mumax^2\sqrt{\frac{r^2}{d_1d_2}}\|\what G^{(t)}-\what O_U^{(t)\tran}\Lambda\what O_V^{(t)}\|.
\end{align*}
Together with Lemma~\ref{lem:hatG}, 
\begin{align*}
\mumax^2\sqrt{\frac{r^2}{d_1d_2}}&\|\what G^{(t)}-\what O_U^{(t)\tran}\Lambda\what O_V^{(t)}\|\\
\leq& C_5\mumax^2\sigma_{\xi}\sqrt{\frac{r^3\log d_1}{N_0}}
+4\mumax^2\|\Lambda\|\sqrt{\frac{r^2}{d_1d_2}}\big(\|\what U^{(t)}-U\what O_U^{(t)}\|^2+\|\what V^{(t)}-V\what O_V^{(t)}\|^2\big)\\
&+C_7\mumax^2\|\Lambda\|\Big(\sqrt{\frac{r}{d_2}}\|\what U^{(t)}-U\what O_U^{(t)}\|_{\submaxx}+\sqrt{\frac{r}{d_1}}\|\what V^{(t)}-V\what O_V^{(t)}\|_{\submaxx}\big)\cdot\mumax\sqrt{\frac{r^3\log d_1}{N_0}}.
\end{align*}
If $n\geq C_3\alpha_d\mumax^4r^3\log^2 d_1$ and $\|\what U^{(t)}-U\what O_U^{(t)}\|+\|\what V^{(t)}-V\what O_V^{(t)}\|\leq 1/(10\mumax\sqrt{r})$, then
\begin{align}
\big\|\what M^{(t)}-M\big\|_{\submax}\leq& \sigma_{\xi}+2\mumax\|\Lambda\|\cdot\Big(\sqrt{\frac{r}{d_2}}\|\what U^{(t)}-U\what O_U^{(t)}\|_{\submaxx}+\sqrt{\frac{r}{d_1}}\|\what V^{(t)}-V\what O_V^{(t)}\|_{\submaxx}\Big)\label{eq:hatM_max}
\end{align}
Therefore, if $n\geq C_3\alpha_d\kappa_0^2\mumax^4r^2d_1\log^2 d_1$, with probability at least $1-2d_1^{-2}$, 
\begin{align*}
\|\what E_V^{(t)}&\|_{\submaxx}+\|\what E_{\xi,V}^{(t)}\|_{\submaxx}
\leq& C_3\eta\frac{\sigma_{\xi}}{\lambda_r}\sqrt{\frac{rd_1d_2\log d_1}{N_0}}+\frac{\eta}{40}\cdot \big(\|\what U^{(t)}-U\what O_U^{(t)}\|_{\submaxx}+\|\what V^{(t)}-V\what O_V^{(t)}\|_{\submaxx}\big)
\end{align*}
and as a result
\begin{align*}
\big\|\what U^{(t+0.5)}-U\what O_U^{(t)}&\what L_G^{(t)} \big\|_{\submaxx}\leq \Big(1-\frac{9\eta}{10}\Big)\|\what U^{(t)}-U\what O_U^{(t)}\|_{\submaxx}+C_3\eta\frac{\sigma_{\xi}}{\lambda_r}\sqrt{\frac{rd_1d_2\log d_1}{N_0}}\\
&+\frac{\eta}{8}\cdot \big(\|\what U^{(t)}-U\what O_U^{(t)}\|_{\submaxx}+\|\what V^{(t)}-V\what O_V^{(t)}\|_{\submaxx}\big).
\end{align*}

Next, we investigate the singular values of $\what U^{(t+0.5)}-U\what O_U^{(t)}\what L_G^{(t)}$. Recall
\begin{align}
\what U^{(t+0.5)}=\underbrace{\what U^{(t)}\what L_G^{(t)}-\eta\cdot (\what U^{(t)}\what G^{(t)}\what V^{(t)\tran}-U\Lambda V^{\tran})\what V^{(t)}\what R_G^{(t)}(\what\Lambda^{(t)})^{-1}}_{\calI_1}
+\underbrace{\what{E}_{V}^{(t)}+\what{E}_{\xi,V}^{(t)}}_{\calI_2}.\label{eq:hatU0.5_1}
\end{align}
By the independence between $(\what U^{(t)},\what G^{(t)},\what V^{(t)}, \what L_G^{(t)}, \what R_G^{(t)},\what \Lambda^{(t)})$ and $\mathfrak{D}_{2t+1}$, and matrix Bernstein inequality \citep{tropp2012user,koltchinskii2011nuclear}, with probability at least $1-2d_1^{-2}$, 
\begin{align*}
\big\|\what E_V^{(t)}+\what E_{\xi,V}^{(t)} \big\|\leq& C_4\eta\cdot\frac{\sigma_{\xi}+\|\what M^{(t)}-M\|_{\submax}}{\lambda_r}\sqrt{\frac{d_1^2d_2\log d_1}{N_0}}
\leq C_4\eta\cdot\frac{\sigma_{\xi}}{\lambda_r}\sqrt{\frac{d_1^2d_2\log d_1}{N_0}}\\
+C_5\eta\kappa_0\mumax&\sqrt{\frac{ d_1^2d_2\log d_1 }{N_0}}\cdot\Big(\sqrt{\frac{r}{d_2}}\|\what U^{(t)}-U\what O_U^{(t)}\|_{\submaxx}+\sqrt{\frac{r}{d_1}}\|\what V^{(t)}-V\what O_V^{(t)}\|_{\submaxx}\Big)
\end{align*}
where the last inequality is due to (\ref{eq:hatM_max}). 
Note that the singular values of $\what U^{(t+0.5)}$ are the square root of eigenvalues of $\what U^{(t+0.5)\tran}\what U^{(t+0.5)}$. We write
\begin{align*}
\what U^{(t+0.5)\tran}\what U^{(t+0.5)}=\calI_1^{\tran}\calI_1+\calI_2^{\tran}\calI_2+\calI_1^{\tran}\calI_2+\calI_2^{\tran}\calI_1.
\end{align*}
Since $\|\what U^{(t)\tran}(\what U^{(t)}-U\what O_U^{(t)})\|\leq 2\|\what U^{(t)}-U\what O_U^{(t)}\|^2$ and $\|\what V^{(t)\tran}(\what V^{(t)}-V\what O_V^{(t)})\|\leq 2\|\what V^{(t)}-V\what O_V^{(t)}\|^2$, by Lemma~\ref{lem:hatG}, we get 
\begin{align*}
\big\|\calI_1^{\tran}&\calI_1-I\big\|\\
\leq& 3\kappa_0\eta\cdot\big(\|\what U^{(t)}-U\what O_U^{(t)}\|^2+\|\what V^{(t)}-V\what O_V^{(t)}\|^2\big)+2\lambda_r^{-1}\eta\cdot\|\what O_U^{(t)}\what G^{(t)}-\Lambda\what O_V^{(t)}\|\\
&+2\eta^2\cdot \big(\kappa_0^2\|\what U^{(t)}-U\what O_U^{(t)}\|^2+\kappa_0^2\|\what V^{(t)}-V\what O_V^{(t)}\|^2+\lambda_r^{-2}\|\what O_U^{(t)}\what G^{(t)}-\Lambda\what O_V^{(t)}\|^2\big)\\
\leq&3(\kappa_0^2\eta^2+\kappa_0\eta)\cdot\big(\|\what U^{(t)}-U\what O_U^{(t)}\|^2+\|\what V^{(t)}-V\what O_V^{(t)}\|^2\big)+4\lambda_r^{-1}\eta\cdot\|\what O_U^{(t)}\what G^{(t)}-\Lambda\what O_V^{(t)}\|\\
\leq&3(\kappa_0^2\eta^2+\kappa_0\eta)\cdot\big(\|\what U^{(t)}-U\what O_U^{(t)}\|^2+\|\what V^{(t)}-V\what O_V^{(t)}\|^2\big)+C_2\eta\frac{\sigma_{\xi}}{\lambda_r}\cdot\sqrt{\frac{rd_1d_2\log d_1}{N_0}}\\
+C_3&\eta\kappa_0\cdot\Big(\sqrt{\frac{r}{d_2}}\|\what U^{(t)}-U\what O_U^{(t)}\|_{\submaxx}+\sqrt{\frac{r}{d_1}}\|\what V^{(t)}-V\what O_V^{(t)}\|_{\submaxx}\Big)\cdot\mumax\sqrt{\frac{rd_1d_2\log d_1}{N_0}}.
\end{align*}
Similar as above analysis, we have
\begin{align*}
\big\|\calI_2^{\tran}\calI_2\big\|\leq& C_4\eta^2 \frac{\sigma_{\xi}^2}{\lambda_r^2}\cdot \frac{d_1^2d_2\log d_1}{N_0}\\
+&C_5\eta^2\kappa_0^2\mumax^2\frac{\alpha_drd_1d_2\log d_1}{N_0}\cdot \big(\|\what U^{(t)}-U\what O_U^{(t)}\|_{\submaxx}^2+\|\what V^{(t)}-V\what O_V^{(t)}\|_{\submaxx}^2\big).
\end{align*}
When $\|\what U^{(t)}-U\what O_U^{(t)}\|+\|\what V^{(t)}-V\what O_V^{(t)}\|\leq 1/(C_3\kappa_0\sqrt{r})$ so that $\|\calI_1\|\leq 2$, 
 due to the independence between $\calI_1$ and $\mathfrak{D}_{2t+1}$, by matrix Bernstein inequality, we get with probability at least $1-2d_1^{-2}$,
\begin{align*}
\big\|\calI_1^{\tran}\calI_2+\calI_2^{\tran}\calI_1\big\|\leq C_3\eta \frac{\sigma_{\xi}+\|\what M^{(t)}-M\|_{\submax}}{\lambda_r}\cdot \sqrt{\frac{rd_1d_2\log d_1}{N_0}}
\leq C_3\eta\frac{\sigma_{\xi}}{\lambda_r}\cdot\sqrt{\frac{rd_1d_2\log d_1}{N_0}}\\
+C_4\eta\mumax\kappa_0\cdot \sqrt{\frac{rd_1d_2\log d_1}{N_0}}\cdot\Big(\sqrt{\frac{r}{d_1}}\|\what U^{(t)}-U\what O_U^{(t)}\|_{\submaxx}+\sqrt{\frac{r}{d_1}}\|\what V^{(t)}-V\what O_V^{(t)}\|_{\submaxx}\Big)
\end{align*}
where the last inequality is due to (\ref{eq:hatM_max}). 
Therefore, with probability at least $1-4d_1^{-2}$, 
\begin{align*}
\big\|&\what U^{(t+0.5)\tran}\what U^{(t+0.5)}-I\big\|\\
\leq& C_3\eta\frac{\sigma_{\xi}}{\lambda_r}\cdot\sqrt{\frac{d_1^2d_2\log d_1}{N_0}}+4(\kappa_0\eta+\kappa_0^2\eta^2)\cdot\big(\|\what U^{(t)}-U\what O_U^{(t)}\|^2+\|\what V^{(t)}-V\what O_V^{(t)}\|^2\big)\\
+C_3&\eta\kappa_0\cdot\Big(\sqrt{\frac{r}{d_2}}\|\what U^{(t)}-U\what O_U^{(t)}\|_{\submaxx}+\sqrt{\frac{r}{d_1}}\|\what V^{(t)}-V\what O_V^{(t)}\|_{\submaxx}\Big)\cdot\mumax\sqrt{\frac{rd_1^2d_2}{N_0}}.
\end{align*}
implying that 
\begin{align*}
\big\{|1-\lambda_r&(\what U^{(t+0.5)})|, |1-\lambda_1(\what U^{(t+0.5)})| \big\}
\leq C_3\eta\frac{\sigma_{\xi}}{\lambda_r}\cdot\sqrt{\frac{d_1^2d_2\log d_1}{N_0}}\\
&+C_4(\kappa_0\eta+\kappa_0^2\eta^2)\cdot\big(\|\what U^{(t)}-U\what O_U^{(t)}\|^2+\|\what V^{(t)}-V\what O_V^{(t)}\|^2\big)\\
+C_5&\eta\kappa_0\cdot\Big(\sqrt{\frac{r}{d_2}}\|\what U^{(t)}-U\what O_U^{(t)}\|_{\submaxx}+\sqrt{\frac{r}{d_1}}\|\what V^{(t)}-V\what O_V^{(t)}\|_{\submaxx}\Big)\cdot\mumax\sqrt{\frac{rd_1^2d_2}{N_0}}
\end{align*}
which concludes the proof of Lemma~\ref{lem:gd_infty}.

\end{document}